\documentclass[letterpaper,10pt, oneside]{article} %arxiv-version+new
\usepackage{stmaryrd}
\usepackage{latexsym, amsmath, amssymb, amsfonts, amscd, amsthm}
\usepackage{bbm}
\usepackage{amsxtra}

\usepackage[utf8]{inputenc}
\usepackage[T1]{fontenc}
\usepackage{lmodern}
\usepackage[mathscr]{eucal}
\usepackage{pgf,calc,pgffor,xkeyval}
\usepackage[all]{xy}
\usepackage{tikz}
\usepackage{xcolor}
\usepackage{booktabs} % To thicken table lines
\newtheorem{thm}{Theorem}[section]
\newtheorem{lem}[thm]{Lemma}\newtheorem{lemma}[thm]{Lemma}
\newtheorem{cor}[thm]{Corollary}
\newtheorem{pro}[thm]{Proposition}

\newtheorem{ep}[thm]{Example}
\newtheorem{rmk}[thm]{Remark}

\newtheorem{defi}[thm]{Definition}

\newcommand {\emptycomment}[1]{}

\newcommand{\calc}[1]{} %to hide calculation
\newcommand{\Cat}{\mathcal M}% category M--->manifolds
\newcommand{\Sets}{\mathsf{Sets}}% category of sets
\newcommand{\Mfd}{\mathsf{Mfd}}% category of smooth manifolds
\newcommand{\Gpd}{\mathsf{Gpd}}% category of gpd
\newcommand{\nGpd}{\mathsf{nGpd}}% category of n gpd
\newcommand{\tGpd}{\mathsf{2Gpd}}% category of 2 gpd
\newcommand{\ecalgdp}{\mathsf{EC}^{p}_{c}} %prestack of exact
\newcommand{\tlp}{{\mathsf{TL}^{p}_{c}}} %prestack of Transitive
\newcommand{\tlpp}{{{\mathsf{TL}^{p}_{c}}}^+} %plus construction of prestack of Transitive Lie algebroids with connections
\newcommand{\tl}{\mathsf{TL}}
 \newcommand{\tcalgdp}{{\mathsf{TC}^{p}_{c}}} %prestack of transitive
\newcommand{\tc}{\mathsf{TC}}%  usual category of  Courant algebroids
\newcommand{\ecalgdpp}{{\mathsf{EC}^{p}_{c}}^+} %plus construction
\newcommand{\tcalgdpp}{{\mathsf{TC}^{p}_{c}}^+} %plus construction of prestack of transitive
\newcommand{\bstringnp}{\mathsf{BString}(n)^p_{c}} %prestack of
\newcommand{\bstringnpp}{{\mathsf{BString}(n)^p_{c}}^+} %plus
\newcommand{\String}{\mathsf{String}}%string group of something
\newcommand{\BString}{\mathsf{BString}}%string group of something
\newcommand{\Spin}{\mathsf{Spin}}%string group of something
\newcommand{\BSpin}{\mathsf{BSpin}}%string group of something
\newcommand{\B}{\mathsf{B}}%B of something, the stack version

\newcommand{\bs}{\mathsf{s}}%source
\newcommand{\bt}{\mathsf{t}}%target
\newcommand{\tr}{\mathsf{tr}}

\DeclareMathOperator{\im}{im}
\DeclareMathOperator{\holim}{\mathsf{holim}}

\newcommand{\tangu}{{\mathbb{T} U}}
\newcommand{\tangui}{{\mathbb{T} U_i}}

\newcommand{\tanguij}{{\mathbb{T} U_{ij}}}
\newcommand{\tagu}{{\mathbb{T}_\g U}}
\newcommand{\tagui}{{\mathbb{T}_\g U_i}}
\newcommand{\taguij}{{\mathbb{T}_\g U_{ij}}}
\newcommand{\tagv}{{\mathbb{T}_\g V}}
\newcommand{\tagvi}{{\mathbb{T}_\g V_i}}
\newcommand{\tagof}[1]{\mathbb{T}_\g {#1}}
\newcommand{\taof}[1]{\mathbb{T}{#1}}
\newcommand{\T}{{\mathbb{T}}}

\DeclareMathOperator{\Kan}{Kan}     %<-- Kan condition
\newcommand{\op}{\mathrm{op}}% opposite

\newcommand{\brtl}[1]{   [    #1   ]_{\nabla}^T   }

\newcommand{\brtls}[1]{   [    #1   ]_S^T   }

\tikzset{cd/.style={matrix of math nodes,row sep=2em,column sep=2em, text height=1.5ex, text depth=0.5ex}}
\tikzset{cdar/.style={->,auto}} % parentheses added, fixed by lidu

\tikzset{triar/.style={anchor=mid,->}}%style for 1-arrows in triangles
\tikzset{tridar/.style={anchor=mid,double,double equal sign distance,-implies}}%style for 2-arrows in triangles

\usetikzlibrary{positioning}
\usetikzlibrary{matrix,arrows}
\usetikzlibrary{calc}

\newcommand{\cl}{\mathrm{cl}}

\newcommand{\be }{\begin{equation}}
\newcommand{\ee }{\end{equation}}

\newcommand{\pf}{\noindent{\bf Proof.}\ }
\newcommand{\ii}{\mathbbm i}
\newcommand{\cs}{\mathrm{cs}}
\newcommand{\TM}{\theta_{\scriptscriptstyle \mathrm{MC}}}
\newcommand{\R}{\mathbb R}\newcommand{\Z}{\mathbb Z}

\newcommand{\huaB}{\mathcal{B}}%{{\mathcal{E}}}%{\mathcal{B}}
\newcommand{\huaS}{\mathcal{S}}
\newcommand{\huaF}{\mathcal{F}}
\newcommand{\huaG}{\mathcal{G}}

\newcommand{\huaP}{\mathcal{P}}
\newcommand{\CWM}{C^{\infty}(M)}

\newcommand{\g}{\mathfrak g} \newcommand{\so}{\mathfrak{so}}

\newcommand{\frkk}{\mathfrak k}

\newcommand{\frkp}{\mathfrak p}

\newcommand{\frkB}{\mathfrak B}
\newcommand{\frkC}{\mathfrak C}

\newcommand{\frkX}{\mathfrak X}

\def\qed{\hfill ~\vrule height6pt width6pt depth0pt}

\newcommand{\half}{\frac{1}{2}}

\newcommand{\pair}[1]{\left\langle #1\right\rangle}

\newcommand{\Courant}[1]{\left\llbracket  #1\right\rrbracket }

\newcommand{\Inn}{\mathrm{Inn}}

\newcommand{\Ad}{\mathrm{Ad}}
\newcommand{\Aut}{\mathrm{Aut}}

\newcommand{\pp}{\frkp} %cohomology class like p_1

\newcommand{\gl}{\mathfrak {gl}}

\newcommand{\SO}{\mathsf {SO}}
\newcommand{\BSO}{\mathsf {BSO}}

\newcommand{\BO}{\mathsf {BO}}

\newcommand{\ad}{\mathrm{ad}}
\newcommand{\pr}{\mathrm{pr}}

\newcommand{\h}{\mathrm{h}}
\newcommand{\id}{\mathrm{id}}

\newcommand{\tB}{\widetilde{B}}
\newcommand{\tg}{\widetilde{g}}
\newcommand{\tf}{\widetilde{f}}
\newcommand{\tA}{\widetilde{A}}
\newcommand{\ta}{\widetilde{a}}
\newcommand{\tomega}{\widetilde{\omega}}
\newcommand{\ttheta}{\widetilde{\theta}}
\newcommand{\bg}{\bar{g}}
\newcommand{\bPhi}{\bar{\Phi}}
\newcommand{\bareta}{\bar{\eta}}
\newcommand{\bTheta}{\bar{\Theta}}

\begin{document}
\title{
{String principal bundles and Courant algebroids } }
\author{Yunhe Sheng,
Xiaomeng Xu  and
Chenchang Zhu }
\date{}

\maketitle 


\begin{abstract}
Just like Atiyah Lie algebroids encode the infinitesimal symmetries of
principal bundles, exact Courant algebroids are believed to encode the
infinitesimal symmetries of $S^1$-gerbes. At the same time, transitive Courant
algebroids may be viewed as the higher analogue of
Atiyah Lie algebroids, and the non-commutative analogue of exact
Courant algebroids. In this article, we explore what the ``principal
bundles'' behind transitive Courant algebroids are, and they turn out
to be principal 2-bundles of string groups. First, we construct the stack of principal 2-bundles of string groups with connection data.
We prove a lifting theorem for the stack of string principal bundles with connections and show the multiplicity of the lifts once they exist. This is a differential geometrical refinement of what is known for string structures by Redden, Waldorf and Stolz-Teichner. We also extend the result of Bressler and Chen-Sti\'enon-Xu on extension obstruction involving transitive Courant algebroids to the case of transitive Courant algebroids with connections, as a lifting theorem with the description of multiplicity once liftings exist. At the end, we build a morphism between these two stacks. The morphism turns out to be neither injective nor surjective in general, which shows that the process of associating the ``higher Atiyah algebroid'' loses some information and at the same time, only some special transitive Courant algebroids come from string bundles.

\end{abstract}

\setcounter{tocdepth}{2}
\tableofcontents

\section{Introduction}

Just like Atiyah Lie algebroids encode the infinitesimal symmetries of $G$-principal bundles \cite{Atiyah,MK2}, exact Courant algebroids are believed to encode the infinitesimal symmetries of $U(1)$-gerbes (or equivalently $\B U(1)$-2-principal bundles) \cite{Bressler:ca,Col,Rogers:localobser,sevlet, Rogers:prequan}. Transitive Courant algebroids may be viewed as the higher analogue of Atiyah Lie algebroids, and the non-commutative analogue of exact Courant algebroids. In this article, we explore what the ``principal bundles'' behind transitive Courant algebroids are.

First, we notice that there are topological obstructions for the
existence of such transitive Courant algebroids. In
\cite{Bressler:Pclass}, Bressler discovered that the obstruction to
extend an Atiyah Lie algebroid to a transitive Courant algebroid is
given by the real first Pontryagin class. This is further fully
generalized to the case of any regular Courant algebroid by
Chen-Sti\'enon-Xu \cite{ChenRCA} in a differential geometry  setting, where the authors gave a complete classification result. In fact, \v{S}evera outlined some very nice ideas to classify transitive Courant algebroids in a series of private letter exchanges with Weinstein \cite{sevlet}. The role of the Pontryagin class is further developed in his later works \cite{sev15, sev16}.

We then notice that the first Pontryagin class arises as an obstruction in another domain. Motivated by  Stolz-Teichner's program of topological modular forms \cite{stolz-teichner:elliptic-obj}, Redden \cite{redden:thesis} defined a string class on a $\Spin(n)$-principal bundle $\bar{P}\to M$,  as a class $\xi \in H^3(\bar{P}, \Z)$, such that for every point $p \in \bar{P}$ the associated inclusion $i_p: \Spin(n) \to \bar{P}$ by $g\mapsto g\cdot p$ pulls back $\xi$ to the standard generator of $H^3(\Spin(n),\Z)$. He further proved that the obstruction for a $\Spin(n)$-principal bundle $\bar{P}$ over $M$ to admit string classes is provided by the integer class of half the first Pontryagin class $\half p_1(\bar{P})\in H^4(M, \Z)$.  In \cite{waldorf:string-conn}, Waldorf proved that such string classes on $\bar{P}$ are in one-to-one correspondence with the isomorphism classes of trivializations of the Chern-Simons 2-gerbe over $\bar{P}$, where the Chern-Simons 2-gerbe is again characterised by $\half p_1(\bar{P})$. Topologically, Stolz and Teichner view the above string structure on a $\Spin(n)$-principal bundle $ M\xrightarrow{\bar{P}} \BSpin(n)$ as a lift of the structure group of $\bar{P}$ from $\Spin(n)$ to a certain three-connected extension, the {\bf string group} $\String(n)$,
\[
\xymatrix{ & {\BString(n)} \ar[d] \\
M \ar@{.>}[ur] \ar[r]^{\bar{P}} & \BSpin(n).
}
\]
Thus, if we realise ``$\String(n)$-principal bundles'' in differential geometry, one may interpret $\half p_1$  as the lifting obstruction of a $\Spin(n)$-principal bundle to a $\String(n)$-principal bundle. A more familiar fact in this style is that the  lifting obstruction of a $\SO(n)$-principal bundle $M\xrightarrow{P} \BSO(n)$ to a $\Spin(n)$-principal bundle $M\xrightarrow{\tilde{P}} \BSpin(n)$ is given by $w_2(P)$.  In fact, there is a whole sequence, called the Whitehead tower:
\[\dots \to \BString(n)\to \BSpin(n)\to \BSO(n)\to \BO(n), \]
with obstruction $w_1$, $w_2$ and $\half p_1$ respectively. Here $w_1$ and $w_2$ are the first and second  Stiefel-Whitney classes.

Since the topological obstruction for both transitive Courant algebroids and $\String(n)$-principal bundles is provided by the first Pontryagin class, we naturally believe that the principal bundle behind a transitive Courant algebroid is exactly a $\String(n)$-principal bundle.

This belief is also supported by another observation from T-duality. Let us start with two T-dual torus bundles $X, \hat{X}$ over $M$, and matching T-dual $S^1$-gerbes $\mathcal{G}\to X$ and $\hat{\mathcal{G}} \to \hat{X}$, %where $H$ and $\hat{H}$ are closed 3-forms representing the classifying DD-classes of $S^1$-gerbe $\mathcal{G}$ in $H^3(X, \Z)$  and that of  $\hat{\mathcal{G}}$ in $H^3(\hat{X}, \Z)$ respectively.
 Bouwknegt-Evslin-Mathai \cite{bem2004} and Bunke-Schick \cite{bs2005}
 proved that the twisted K-theory for the T-dual pairs are isomorphic,
 that is,  there is an isomorphism between twisted K-groups $
 K^\bullet(X, \mathcal{G}) \cong K^{\bullet } (\hat{X}, \hat{ \mathcal{G}})$.  On the level of differential geometric objects,  Cavalcanti-Gualtieri \cite{CG} proved that the exact Courant algebroid associated to the T-dual $S^1$-gerbes are the same. Now we extend this story $\Spin(n)$-equivariantly.  We  begin with two T-dual torus bundles $X, \hat{X}$ over $M$, and their matching T-dual string structures $(P, \xi)\to X$ and $(\hat{P}, \hat{\xi}) \to \hat{X}$, where $P$ and $\hat{P}$ are $\Spin(n)$-principal bundles over $X$ and $\hat{X}$ respectively, and $\xi,\hat{\xi}$ are string classes on $P$ and $\hat{P}$ respectively. Leaving alone what the cohomological invariants should be, Baraglia and Hekmati \cite{baraglia-hekmati}  showed that, on the level of differential geometric objects, the transitive Courant algebroids associated to both sides are isomorphic.

In this article, we realise $\String(n)$-principal bundles and their connections as differential geometric objects by describing the entire $(3,1)$-sheaf (or 2-stack) $\bstringnpp$. Then we make the connection between transitive Courant algebroids and string principal bundles explicit and functorial by constructing a morphism between their corresponding stacks.

For this purpose, first we study what a $\String(n)$-principal bundle with connection data really is. As we have seen, $\String(n)$ is a three-connected cover of $\Spin(n)$, and this forces the model of $\String(n)$ to be either infinite-dimensional or finite-dimensional however higher (namely being a Lie 2-group)\footnote{There are also models which are both higher and infinite-dimensional.}. We take the second approach with the model of Schommer-Pries \cite{schommer-pries:string} for $\String(n)$. The advantage of this model is that the spaces it involves are all nice finite dimensional manifolds, thus there is no additional analytic difficulty when solving equations or constructing covers; at the same time, this is paid off by algebraic difficulty of chasing through various pages of spectral sequences of cohomological calculation.

First we  construct a $(3,1)$-presheaf of $\String(n)$-principal bundles with connection data $\bstringnp$ and complete it into a $(3,1)$-sheaf (or a 2-stack) $\bstringnpp$ using the plus construction. This is essentially to build a $\String(n)$-principal bundle with a connection from local data and gluing conditions in the fashion of Breen-Messing. Breen and Messing studied connections for gerbes in their original work  \cite{breen-messing}. We also notice that in a recent work \cite{waldorf:global},  connections for 2-principal bundles of strict 2-groups are studied both locally and globally. However, the finite-dimensional differential geometric model for $\String(n)$ is a non-strict Lie 2-group. This forces us to develop our own formula instead of using existing results in literature.
It turns out that the glued stack involves descent equations of first Pontryagin class. In a recent work \cite{AFXZ}, these equations are further studied in a universal setting and proved to be closely related to Kashiwara-Vergne theory and Drinfeld associators.

To justify our construction, we prove directly the lifting theorem that one
expects for $\String(n)$-principal bundles and provide a comparison
to previous string concepts of Stolz-Teichner, Redden and Waldorf
respectively in Section \ref{sec:string-lift}:

\begin{thm} \label{thm:main}
Given a $\Spin(n)$-principal bundle with a connection $M\xrightarrow{\bar{P}} \BSpin(n)_{c}$,
\begin{enumerate}
    \item[\rm(i)] it lifts to an object in $M\xrightarrow{P_{c}}{\bstringnp}^+$,
\begin{equation}\label{lift}
  \xymatrix{ & {\bstringnp}^+ \ar[d] \\
M \ar@{.>}[ur] \ar[r]^{\bar{P}} & \BSpin(n)_{c},
}
\end{equation}
if and only if $\frac{1}{2} p_1(\bar{P})=0$;

    \item[\rm(ii)] if a lift in \eqref{lift} exists, then the isomorphism classes of different lifts form a torsor of the Deligne cohomology group $H^2(M, \underline{U(1)} \xrightarrow{d\log} \Omega^1\xrightarrow{d}\Omega^2)$ mod out by a certain subspace $I$.
\end{enumerate}
\end{thm}
It is proved in Theorem \ref{thm:lift} and Theorem \ref{thm:torsor}.

After this, we build the $(2,1)$-sheaf (or 1-stack) of transitive Courant algebroids with connections. We benefit much from \cite{ChenRCA} where transitive Courant algebroids and their gauge transformations are well studied.  However, the  gauge transformations which preserve the connection data are still needed to be specified. We thus have additional equations in the definition of 1-morphisms (see Eqs. \eqref{eq:congluetc1}-\eqref{eq:congluetc3}). To make the construction mathematically strict, however at the same time avoiding the routine checking of gluing conditions of stacks over several layers, as before, we first   construct   a $(2,1)$-presheaf $\tcalgdp$ by simply mapping to the category of {\em standard}  transitive Courant algebroids  with connections and their gauge transformations. Then we complete it to a $(2,1)$-sheaf $\tcalgdpp$ using Nikolaus-Schweigert's plus construction \cite{nikolaus-schweigert}.  We prove that the gluing result gives us exactly a transitive Courant algebroid (not necessarily standard) with a connection and its gauge transformations. This in turn justifies our construction of the transitive Courant algebroid stack. There is a subtle difference between our construction  of the transitive Courant algebroid stack and the one in \cite{Bressler:Pclass}. In \cite{Bressler:Pclass}, the stack is directly taken to be a functor mapping to the category of transitive Courant algebroids (not just standard ones), however the checking of gluing conditions seems to be omitted. Also in the  language of stacks,  a recent work \cite{Pym} has studied the relation between twisted Courant algebroids and shifted symplectic Lie algebroids, and has further hinted an even higher correspondence of our type involving fivebrane structures.

In the end, we construct a morphism from the $(3,1)$-sheaf of $\String(n)$-principal bundles with connections to the $(2,1)$-sheaf of transitive Courant algebroids with connections for $\Spin(n)$. To achieve this, we only need to build a morphism on the presheaf level since the plus construction is functorial. It turns out that the difficulty of the construction lies on the level of morphisms, that is, to construct the gauge transformation of transitive Courant algebroids associated to that of $\String(n)$-principal bundles. The formula of the symmetric part of the $(3,1)$-position in the gauge transformation remains rather mysterious. \v{S}evera suggests  us  some connection to Alekseev-Malkin-Meinrenken's theory on group valued moment maps \cite{AMM}. We reserve it for future investigation. We remark in Appendix  \ref{app:innerauto} that these gauge transformations are all inner ones noticed by \v{S}evera. Similar results of these inner automorphisms are also studied in \cite{GRT} in another setting.  We further verify that this morphism from the string stack to the Courant stack is neither injective nor surjective. This tells us that the process of associating a ``higher Atiyah algebroid'' to a string principal 2-bundle loses some information, and, at the same time, not all transitive Courant algebroids come from this process.

\section{Preliminaries on prestacks, stacks and the plus construction}

Recall that an $(n+1,1)$-presheaf  over a category $\Cat$ is a (higher) functor
$\Cat^{\op} \to \nGpd$ to the higher category of $n$-groupoids, where $n \in
\{0, 1, 2, \dots\} \cup \{\infty\}$. Here $0$-groupoids are
interpreted as sets. Therefore, a $(1, 1)$-presheaf (or a
 presheaf) over a category $\Cat$  is a functor $\Cat^{\op} \to \Sets$
to the category of sets;  a $(2, 1)$-presheaf over a
category $\Cat$ is a (higher) functor $\Cat^{\op} \to \Gpd$ to the 2-category of
(1-)groupoids;  and a $(3, 1)$-presheaf over a category $\Cat$ is a (higher)
functor $\Cat^{\op} \to \tGpd$ to the higher category of $2$-groupoids. These are
all the cases that we will use in this paper. Then
we perform a plus construction
(namely a procedure of higher sheafification) to obtain the corresponding
sheaves. Sometimes, $(2, 1)$-sheaves are also called stacks, and $(3,1)$-sheaves
are called 2-stacks. The model we use is as in
\cite[Section 2]{nikolaus-schweigert}. For technical details, we refer
readers to this paper and the references therein.

Here we briefly recall the model we use for 2-groupoids. Our model
for a
2-groupoid (in the sense of Duskin and Glenn \cite{duskin, glenn}) is a simplicial set satisfying Kan conditions $\Kan(n, j)$
for all $n\ge 1$ and $0\le j \le n$ and strict
Kan conditions $\Kan!(n, j)$ for all $n\ge 3$ and $0\le j \le n$. Readers who are not familiar with Kan conditions may equivalently
understand it as a bicategory \cite{ben}, whose 2-morphisms are
invertible and whose 1-morphisms are all invertible up to
2-morphisms. The compositions of 1-morphisms are associative up to an
associator, and the associator in turn satisfies a higher coherence
condition. For the precise definition, we refer to \cite[Definition
5.2]{sheng-zhu3}, where a semi-strict Lie 2-groupoid is defined. If we
equip the object therein with discrete topology, we obtain what a
2-groupoid is. Let us also recall the equivalence between the two
different descriptions: if we start with a simplicial set $X_\bullet$
satisfying the above Kan condition, then we take $C_0=X_0$ on the
object level; $C_1=X_1$ on the 1-morphism level; $C_2= d_0^{-1} (
s_0(X_0))$ on the 2-morphism level, we obtain a bicategory $(C_0, C_1,
C_2)$ satisfying
required conditions; as for the other direction, we take $X_0=C_0$,
$X_1=C_1$ and $X_2= C_1\times_{C_0}C_1$.  For details we refer to
\cite[Section 4]{z:tgpd}.

Now we shortly recall the process of the plus construction in the case when $\Cat=\Mfd$ is
the category of differential manifolds for our application. Given a $(3, 1)$-presheaf $\huaF: \Mfd^\op \to \tGpd$, the plus construction in \cite{nikolaus-schweigert}
gives us a $(3,1)$-sheaf
$\huaF^+: \Mfd^\op \to \tGpd$. To describe this $(3, 1)$-sheaf, we first need to take
the homotopy limit $\holim
\huaF(U(M)_\bullet)$ for an open cover $\{U_i\}$ of $M$ over the \v{C}ech
simplicial manifold
\begin{equation}\label{defi:UM}
\xymatrix{U(M)_\bullet=  \sqcup U_i  &
  \ar@2[l]^{\qquad\partial_1}_{\qquad\partial_0} \sqcup U_{ij} & \sqcup
              U_{ijk} \ar@3[l]_{\partial_0, \partial_1}^{\partial_2} \dots .}
  \end{equation}
%which is the simplicial
%nerve of the \v{C}ech groupoid $ \sqcup U_{ij} \Rightarrow \sqcup
%U_i$.
Let us describe  $\holim
\huaF(U(M)_\bullet)$ explicitly: the result is a
2-groupoid.

\begin{itemize}
\item[$\bullet$]Its object consists of
\begin{enumerate}
\item[Ob0]   an element $\theta=(\theta_i)\in \huaF(\sqcup U_i)_0$;
\item[Ob1]   an element $g=(g_{ij}) \in \huaF(\sqcup U_{ij})_1$, which
  is a 1-morphism $
  \theta_i|_{U_{ij}} \xleftarrow{g_{ij}} \theta_j|_{U_{ij}}$, or
  equivalently, a 1-morphism $ \partial^*_1 \theta \xleftarrow{g} \partial_0^*\theta$;
\item[Ob2]  an element $a=(a_{ijk})\in \huaF(\sqcup U_{ijk} )_2$, which is a 2-morphism $a:
  g_{ij}\circ g_{jk} \Leftarrow g_{ik}$;
\item[Ob3] pentagon condition for $a$, that is, $ (\id\circ_\h \partial^*_0 a) \circ \partial_2^* a
 = \alpha(g_{ij}, g_{jk}, g_{kl}) \circ (\partial_3^* a
  \circ_\h \id) \circ \partial_1^* a$, where $\alpha(g_{ij}, g_{jk}, g_{kl})$ is the
  associator $g_{ij} \circ (g_{jk} \circ g_{kl}) \Leftarrow (g_{ij} \circ g_{jk} ) \circ g_{kl}
  $, where $\circ_\h$ is the
  horizontal composition of 2-morphisms.
\end{enumerate}
\item[$\bullet$] A 1-morphism from $(\ttheta, \tg, \ta)$ to $(\theta, g, a)$ consists of
\begin{enumerate}
\item[1M0] a 1-morphism $A=(A_i): \theta \leftarrow \ttheta$ in $\huaF(\sqcup U_i)$;
\item[1M1] a 2-morphism $f:    g
  \circ \partial^*_0 A \Leftarrow  \partial_1^* A \circ \tg $ in $\huaF(\sqcup U_{ij})$;
\item[1M2] a higher coherence condition, $(a\circ_\h \id)^{-1} \circ
 (\id \circ_\h \partial_0^* f) \circ   (\partial^*_2 f \circ_\h \id )
  = \partial_1^* f \circ (  \id \circ_\h \ta)^{-1} $ of 2-morphisms in
  $\huaF(\sqcup U_{ijk})$.
\end{enumerate}
\item[$\bullet$] A 2-morphism from $(\tA, \tf)$ to $(A, f)$ consists of
\begin{enumerate}
\item[2M0]  a 2-morphism $\omega: A \Leftarrow \tA$ in
  $\huaF(\sqcup U_i)$;
\item[2M1] a higher coherence condition $ f \circ   (\partial_1^* \omega \circ_\h \id)=
  (\id
  \circ_\h \partial_0^* \omega)  \circ \tf $.
\end{enumerate}
\end{itemize}

Then the $(3,1)$-sheaf $\huaF^+$ maps $M \in \Cat$ to the following
2-groupoid:
\begin{itemize}
\item $\huaF^+(M)_0$: an object is a pair $(\{U_i\}, P)$, where
  $\{U_i\}$ is a cover of $M$ and $P$ is an object in $\holim \huaF (U(M)_\bullet)$;
\item $\huaF^+(M)_1$: a 1-morphism from $(\{\widetilde{U}_i\}, \widetilde{P})$ to $(\{U_i\},
  P)$ is a common refinement $\{V_i\}$ of $\{\widetilde{U}_i\}$ and $\{U_i\}$,
  and a
  1-morphism $\phi $ in $\holim \huaF (V(M)_\bullet)$;
\item $\huaF^+(M)_2$: a 2-morphism from $(\{\widetilde{V}_i\}, \widetilde{\phi})$ to $(\{V_i\},
  \phi)$ is a common refinement $\{W_i\}$ of $\{\widetilde{V}_i\}$ and $\{V_i\}$, and
  a 2-morphism $\alpha$ in $\holim \huaF (W(M)_\bullet)$. Moreover,  $(\{\widetilde{W}_i\}, \widetilde{\alpha})$ and $(\{W_i\},
  \alpha)$ are identified if $\alpha$ and $\widetilde{\alpha}$ are
  identified on a further common refinement of $\{\widetilde{W}_i\}$ and $\{W_i\}$.
\end{itemize}

The plus construction for $(2, 1)$-presheaves is then a truncation of that of
$(3, 1)$-sheaves viewing 1-groupoids as 2-groupoids with identity
2-morphisms. Let us explain it with a nice example.

\begin{ep}\label{ep:bgconn}{\rm
Given a Lie group $G$ and its Lie algebra $\g$, there is a $(2, 1)$-presheaf $\B G^p_{c}: \Mfd^\op \to \Gpd$
sending $U\in \Mfd$ to the groupoid whose objects are trivial $G$-principal bundles
$U\times G$ together with $\theta\in \Omega^1(U, \g)$ and whose morphisms from $(U\times G,\ttheta)$ to $(U\times G,\theta)$ are gauge
transformations $g: U\to G$ satisfying $\theta - \Ad_g\ttheta = - g^* \TM$;  and sending a
morphism $U\to V$ to the functor between the corresponding groupoids
induced by pullbacks of principal bundles and differential forms. Here $\TM$ is the right invariant Maurer-Cartan form on $G$. It satisfies the following Maurer-Cartan equation
$$
d\TM-\half[\TM,\TM]_\g=0.
$$Let us form the $\holim
\B G^p_c(U(M)_\bullet)$ with respect to an open cover $\{ U_i\}$ of $M\in
\Mfd$. An object in $\holim
\B G^p_{c}(U(M)_\bullet)$ consists of
\begin{itemize}
\item $ U_i \times G, \, \theta_i \in \Omega^1 (U_i, \g)$;
\item $g_{ij} : U_{ij} \to G$,   $(U_i \times G, \theta_i)
  \xleftarrow{g_{ij}} (U_j \times G, \theta_j)$, with  $\theta_i - \ad_{g_{ij}} \theta_j = - g_{ij}^* \TM$;
\item compatibility condition $g_{ij} \circ g_{jk}=g_{ik}$ on $U_{ijk}$.
\end{itemize}
Thus, we see that such an object gives us exactly the local data of a
$G$-principal bundle with a connection 1-form.
A morphism in $\holim
\B G^p_c(U(M)_\bullet)$ from $(\ttheta_i; \tg_{ij})$ to $(\theta_i; g_{ij})$ consists of
\begin{itemize}
\item $g_i : U_i \to G$, $\theta_i- \ad_{g_i}\ttheta_i = -g_i^* \TM$;
\item compatibility condition $g_{ij} \cdot g_j = g_i \cdot \tg_{ij} $ on $U_{ij}$.
\end{itemize}
This gives us exactly the local data of a gauge transformation  preserving connections between
the corresponding $G$-principal bundle glued by $\tg_{ij}$ and
$g_{ij}$.
}
 \end{ep}

\emptycomment{
Then equivalence between two such functors $X$ and $Y$ is a stalk-wise
weak equivalence. Let's suppose $X: \Cat \to \tGpd$ has level $X_0,
X_1, B^X_2$, where $\Theta(X):=B^X_2\Rightarrow X_1$ is a groupoid. Notice that
$X_i$ becomes then a presheaf values in $\Sets$ on $\Cat$. Then a
stalk-wise weak equivalence $X\xrightarrow{f} Y$ is  to make
\[
\begin{split}
X_0 \times_{Y_0} Y_1 \to Y_0, \quad \text{stalkwise surjective} \\
X_1\times_{Y_1} B^Y_2 \to X_0\times X_0 \times_{Y_0\times Y_0} Y_1, \quad \text{stalkwise surjective} \\
B_2^X \cong B_2^Y\times_{Y_1\times Y_1} X_1 \times X_1
\times_{Y_0\times Y_0} X_0 \times X_0
\end{split}
\]
The last is an isomorphism, and we know that stalkwise isomorphism is
equivalent to objectwise isomorphism (that is the two sheaves are
isomorphic). We notice that objectwise surjectivity implies stalkwise
surjectivity. If $X$ and $Y$ are just $(2,1)$-sheaves, that is, they
take value in $\Gpd$ instead of $\tGpd$, we then have corresponding
condition:
\[
\begin{split}
X_0 \times_{Y_0} Y_1 \to Y_0, \quad \text{stalkwise surjective} \\
X_1\cong Y_1\times_{Y_0\times Y_0}X_0\times X_0.
\end{split}
\]
We notice that it is the usual stack isomorphism: the first condition
is to ask the morphism to be essentially surjective, the second fully
faithful. Thus the last two conditions in the 2-stack case is in fact
to ask $\Theta(X) \to \Theta(Y)|_{X_0}$ to be a weak equivalence for 1-groupoids.
(see also later in the example of principal bundles)
}
\section{$(3,1)$-sheaf of string principal bundles}

\subsection{Finite dimensional model of $\String_\pp(G)$}\label{sec:string-gp}

In this section, $G$ is a finite dimensional compact  Lie
group. Let us first recall the finite dimensional model of the Lie 2-group
$\String_{\pp}(G)$ built in \cite{schommer-pries:string} for a given class $\pp \in H^4(BG_\bullet, \Z)$. The idea is to realise $\String_\pp(G)$ as a $\B U(1)$-central extension of a Lie group
$G$
$$
\B U(1)\longrightarrow \String_\pp(G)\longrightarrow G,
$$
with the extension class  $\pp \in H^4(BG_\bullet, \Z)$. Here $BG_\bullet$ is the simplicial nerve of a Lie
group $G$, and  $H^\bullet(BG_\bullet, \Z)$  denotes the sheaf cohomology of
the sheaf of $\Z$-valued functions on $BG_\bullet$. Similarly, we use $H^\bullet(BG_\bullet,
\underline{U(1)})$ and $H^\bullet(BG_\bullet, \underline{\R})$ to denote the
sheaf cohomology of $U(1)$-valued function and $\R$-valued function
respectively.

Let us explain a bit the terminology here: a
Lie 2-group is a differentiable stack equipped with a
group structure (up to homotopy). For example, $\B U(1)$ is an abelian
Lie 2-group. Here $\B U(1)$ denotes the stack presented by groupoid
$U(1) \Rightarrow pt$. The multiplication $m: \B U(1) \times \B U(1)
\to \B U(1)$ is induced by the multiplication of $U(1)$. Notice that
since $U(1)$ is abelian, thus the $U(1)$-multiplication is a functor
\[
(U(1) \Rightarrow pt) \times (U(1) \Rightarrow pt) \longrightarrow
(U(1) \Rightarrow pt).
\]
For more details in the topic of (Lie) 2-groups and examples see e.g. \cite{baez:2group}\cite[Sect 3.1]{bnz}. In the case that $G=\Spin(n)$, the generator of $H^4(BG_\bullet, \Z)$ is given by $\half p_1$, half
of the Pontryagin class. The sheaf cohomology may be calculated by taking a
hypercover of $BG_\bullet$ and taking the \v{C}ech cohomology, as long as the
hypercover is acyclic. This is a classical result. See \cite[Proposition 2.4]{wz:int} for a concrete
statement. See \cite{friedlander} and \cite[Sect. 2]{wz:int} and references therein for definition and properties of the sheaf cohomology for simplicial objects. The cohomology on $BG_\bullet$ is then equivalent to the group cohomology used in \cite{schommer-pries:string} originally coming from Segal and Brylinski \cite{Seg70, Seg75, Bry00}.  In our case, as long as
the cover of $BG_\bullet$ on each layer $G^{(\bullet)}$ is good, namely intersections are contractible, it is acyclic with respect to sheaves in our study.

The short exact sequence of sheaves $0\to \Z \to \R\to U(1)\to 0$ gives us a long exact sequence of cohomology. Notice that $H^{\ge 1} (BG_\bullet, \underline{\R})=0$ for compact group $G$. Thus $H^{n}(BG_\bullet, \underline{U(1)})\cong H^{n+1}(BG_\bullet, \Z)$ for $n\ge 1$.

To build the finite dimensional model for the Lie 2-group $\String_\pp(G)$, let us
take a good simplicial hypercover $G^{(\bullet)}$ for $BG_\bullet$ and write the
simplicial-\v{C}ech double complex whose total cohomology is
$H^\bullet (BG_\bullet, \underline{U(1)})$.
\begin{equation}\label{db-cx:BG}
 \xymatrix{
 &&&&&&
\\
 C(\sqcup G^{(3)}_{s},\underline{U(1)})\ar[r]^{\check{\delta}}\ar[u]^{\delta}
 &C(\sqcup G^{(3)}_{s,t},\underline{U(1)})\ar[r]^{\check{\delta}}\ar[u]^{\delta}
 &\dots%\Omega^3(U^1)\ar[r]^{\check{\delta}}\ar[u]^{-d}
&%\Omega^3(U^2)\ar[r]^{\check{\delta}}\ar[u]^{d}
&%\Omega^3(U^3)\ar[r]^{\check{\delta}}\ar[u]^{-d}
&
 \\
  C(\sqcup G^{(2)}_{p},\underline{U(1)})\ar[r]^{\check{\delta}}\ar[u]^{\delta}
  &C(\sqcup G^{(2)}_{p,q},\underline{U(1)})\ar[r]^{\check{\delta}}\ar[u]^{\delta}
 &C(\sqcup G^{(2)}_{p,q,r},\underline{U(1)})\ar[r]^{\check{\delta}}\ar[u]^{\delta}
&\dots %\Omega^2(U^2)\ar[r]^{\check{\delta}}\ar[u]^{d}
& %\Omega^2(U^3)\ar[r]^{\check{\delta}}\ar[u]^{-d}
&
 \\
 C(\sqcup G^{(1)}_{\alpha},\underline{U(1)})\ar[r]^{\check{\delta}}\ar[u]^{\delta}
 & C(\sqcup G^{(1)}_{\alpha,\beta},\underline{U(1)})\ar[r]^{\check{\delta}}\ar[u]^{\delta}
 &C(\sqcup G^{(1)}_{\alpha,\beta,\gamma},\underline{U(1)})\ar[r]^{\check{\delta}}\ar[u]^{\delta}
&C(\sqcup G^{(1)}_{\alpha,\beta,\gamma,\delta},\underline{U(1)})\ar[r]^{\qquad\check{\delta}}\ar[u]^{\delta}
&\dots %\Omega^1(U^3)\ar[r]^{\check{\delta}}\ar[u]^{-d}
&
 \\
 C(\cdot,\underline{U(1)})\ar[r]^{0}\ar[u]^{\delta}
 &C(\cdot,\underline{U(1)})\ar[r]^{id} \ar[u]^{\delta}
&C(\cdot,\underline{U(1)})\ar[r]^{0}\ar[u]^{\delta}
&C(\cdot,\underline{U(1)})\ar[r]^{id} \ar[u]^{\delta}
&\dots %C(\cdot,\underline{U(1)})\ar[r]^{\check{\delta}} \ar[u]^{\delta}
&
}
 \end{equation}
We take a representative $(\Theta, \Phi, \eta, 0)$ of $\pp \in H^3(BG_\bullet, \underline{U(1)})=H^4(BG_\bullet, \Z)$, where $$\Theta\in C(\sqcup
 G^{(3)}_{s},\underline{U(1)}), \quad\Phi\in C(\sqcup G^{(2)}_{p,q},\underline{U(1)}),\quad \eta \in
 C(\sqcup G^{(1)}_{\alpha, \beta, \gamma}, \underline{U(1)}).$$ The last entry being
 0 is implied by the closedness.

To build up a  Lie 2-group, we first need to have an underlying
Lie groupoid which presents the stack $\String_\pp(G)$, and then establish a group structure ``up to homotopy''
on top of it. Here we follow the
convention in \cite[Section 2]{wz:int}.

Our underlying Lie groupoid $\Gamma[\eta]$ is a $U(1)$-extension of
the \v{C}ech groupoid with respect to the cover $G^{(1)}$, that is $
\sqcup G^{(1)}_{\alpha,\beta}\times U(1)\Rightarrow \sqcup
G^{(1)}_\alpha$, together with source and target
$\bs(g_{\alpha,\beta},a)=g_\alpha, \bt(g_{\alpha,\beta},a)=g_\beta$,  multiplication
$$(g_{\alpha, \beta},a)(g_{\beta,\gamma},b)=(g_{\alpha,\gamma},a+b-\eta_{\alpha,\beta,\gamma}(g)),$$
identity  $e(g_\alpha)=(g_{\alpha,\alpha},0)$, and inverse
$(g_{\alpha, \beta},a)^{-1}=(g_{\beta,\alpha},-a)$. Here
$\check{\delta}(\eta)=0$ guarantees that the construction gives rise to a Lie groupoid structure.

Now we build the multiplication for the 2-group structure on
$\Gamma[\eta]$, which should be a generalized morphism
$\Gamma[\eta]^{\times 2}\longrightarrow \Gamma(\eta)$. We realize the
generalized morphism by a span of a Morita morphism and a usual
morphism,  $$\Gamma[\eta]^{\times
  2} \stackrel{M.E.}{\longleftarrow}
\Gamma^2[\eta]\stackrel{m}{\longrightarrow}\Gamma(\eta),$$ where
$\Gamma^2[\eta]= \Big( G_{[1]}^{(2)}\times {U(1)}^{\times
  2}\Longrightarrow G_{[0]}^{(2)} \Big)$, is similarly constructed as
$\Gamma[\eta]$, however, by the pullback \v{C}ech
cocycle  $(d_0^*(\eta),d_2^*(\eta))\in C(G^{(2)}_{[2]},{U(1)}^{\times
  2})$. Here, and later, $G^{(j)}_{[i]}$ denotes the disjoint union of
$(i+1)$-fold intersections of the hypercover $G^{(j)}$.

The natural projection $\Gamma^2[\eta]\to \Gamma^{\times
  2}[\eta]$ is a Morita morphism, i.e.,  a morphism that gives arise to a
Morita equivalence of Lie groupoids.

The morphism $\Gamma^2(\eta)\stackrel{m}{\longrightarrow}\Gamma(\eta)$
is given by
$$
(v_0,v_1,a_0,a_1)\stackrel{m}{\longrightarrow} (d_1(v_0),d_1(v_1),a_0+a_1+\Phi(v_0,v_1)).
$$It being a groupoid morphism is equivalent to the fact that
$\delta(\eta)+\check{\delta}(\Phi)=0$. However the multiplication is not
strictly associative, that is, the following diagram of differentiable stacks commutes up to a
2-morphism $a$, which is called an associator,
\[
 \xymatrix{\String_\pp(G)^{\times 3}\ar[r]^{m\times \id}\ar[d]^{\id\times m}&\String_\pp(G)^{\times 2}\ar[d]^{  m}\\
 \String_\pp(G)^{\times 2}\ar[r]^{  m}&\String_\pp(G).
 }
 \]
We now find a suitable Lie groupoid presentation of
$\String_\pp(G)^{\times 3}$ so that certain desired morphisms can be
written as strict morphisms of Lie groupoids. Notice that there are
three maps
$d_0d_0, d_2d_0, d_2d_2:G^{\times 3}\longrightarrow G$. Just like before, we take
$\Gamma^3[\eta]$
to be the Lie groupoid constructed by the pullback cocycle
$\big( (d_0d_0)^*\eta,(d_2d_0)^*\eta,(d_2d_2)^*\eta \big) \in
C(G^{(3)}_{[3]}, {\underline{U(1)}}^{\times 3})$,
that is, $
\Gamma^3[\eta]=\Big( G^{(3)}_{[1]}\times {U(1)}^{\times 3}\Longrightarrow G^{(3)}_{[0]} \Big)
$.
Now the two composed morphism $m_1:  \String_\pp(G)^{\times 3} \xrightarrow{\id\times m}
\String_\pp(G)^{\times 2} \xrightarrow{m} \String_\pp(G) $ and $m_2: \String_\pp(G)^{\times 3} \xrightarrow{m\times \id}
\String_\pp(G)^{\times 2} \xrightarrow{m} \String_\pp(G)$ are given by strict
Lie groupoid morphisms:
%\[
% \xymatrix{\Gamma^3(\eta)\ar[r]^{m\times id}\ar[d]^{id\times m}&\Gamma^2(\eta)\ar[d]^{  m}\\
% \Gamma^2(\eta)\ar[r]^{  m}&\Gamma(\eta)
% }
% \]
 \begin{eqnarray*}
   m_1:((w_s,w_t),a_0,a_1,a_2)&\longmapsto (d_1d_2(w_s),d_1d_2(w_t),a_0+a_1+a_2+d_2^*\Phi(w_s,w_t)+d_0^*\Phi(w_s,w_t)),\\
   m_2:((w_s,w_t),a_0,a_1,a_2)&\longmapsto (d_1d_1(w_s),d_1d_1(w_t),a_0+a_1+a_2+d_1^*\Phi(w_s,w_t)+d_3^*\Phi(w_s,w_t)).
 \end{eqnarray*}
In this model, the associator $a:m_2\Longrightarrow m_1$, is a
map $ \Gamma^3[\eta]_0\longrightarrow \Gamma[\eta]_1$, given by
$$w_0\longmapsto (d_1d_2(w_0),d_1d_2(w_0),\Theta(w_0)). $$
The naturality condition $m_1(r)a(\bs(r))=a(\bt(r))m_2(r)$ is
equivalent to the equation $\delta(\Phi)-\check{\delta}(\Theta)=0$. The pentagon condition for the  associator is equivalent to the equation $\delta(\Theta)=0$.

%Tuesday evening 29.11

\subsection{$(3,1)$-sheaf $\bstringnpp$ of string principal bundles with connection data}
Now we take $G=\Spin(n)$ and  $\pp$ to be $\half p_1 \in H^4(B\Spin(n)_\bullet, \Z)$. We denote by $\String(n)$ the corresponding string group $\String_\pp(G)$.

The $(3,1)$-presheaf of $\String(n)$-principal bundles with connections
$\bstringnp: \Mfd^{\op} \to \tGpd$ is constructed as following: first of
all, $\bstringnp(U)$ is a 2-groupoid made up by the following data:
\begin{itemize}
\item{$\bstringnp(U)_0$}:  an object is a triple $\big( U\times \String(n)\to U, \theta, B \big) $ consisting
 of a locally  trivial $\String(n)$-principal bundle $(U\times \String(n)\to U)$
  together with a 2-form $B\in\Omega^2(U) $ and an $\so(n)$-valued 1-form
  $\theta \in \Omega^1(U, \so(n))$.
\item{$\bstringnp(U)_1$}: a 1-simplex  \[\big( U\times \String(n)\to U,
  \theta_0, B_0 \big) \xleftarrow{(g_{\scriptscriptstyle 01}, A_{\scriptscriptstyle 01}, \omega^2_{\scriptscriptstyle 01})} \big( U\times \String(n)\to U,
  \theta_1, B_1 \big)  \] consists of a generalized morphism $g_{\scriptscriptstyle 01}: U\to \String(n)$ given by
  a bibundle $E_{g_{\scriptscriptstyle 01}}$,  a 1-form $A_{\scriptscriptstyle 01}\in \Omega^1(U)$  and a 2-form
  $\omega^2_{\scriptscriptstyle 01} \in \Omega^2(U)$ such that
\[
B_1-B_0=\omega^2_{\scriptscriptstyle 01} +dA_{\scriptscriptstyle 01}, \quad
\cs_3(\theta_1)-\cs_3(\theta_0) =d\omega^2_{\scriptscriptstyle 01}, \quad \theta_0 -
\ad_{\bar{g}_{\scriptscriptstyle 01}} \theta_1 = -\bar{g}^*_{\scriptscriptstyle 01}\TM,
\]  where $\cs_3(\theta)$ is the Chern-Simon 3-form associated to an $\so(n)$-valued 1-form $\theta$  given by
\begin{equation}\label{eq:cs3-form}
 \cs_3(\theta)=(\theta,d\theta)+\frac{1}{3}(\theta,[\theta,\theta]).
\end{equation}
Here $(-,-)$ is a certain invariant symmetric bilinear form on $\so(n)$,  and $\bar{g}_{\scriptscriptstyle 01}
: U\xrightarrow{g_{\scriptscriptstyle 01}} \String(n) \xrightarrow{\pi} \Spin(n)$ is the composition
of $g_{\scriptscriptstyle 01}$ with the natural projection $\String(n)\xrightarrow{\pi} \Spin(n)$. Composition
is given by the multiplication in $\String(n)$:
\[
(g_{\scriptscriptstyle 01}, A_{\scriptscriptstyle 01}, \omega^2_{\scriptscriptstyle 01}) \circ (g_{\scriptscriptstyle 12}, A_{\scriptscriptstyle 12},
\omega^2_{\scriptscriptstyle 12}) = (g_{\scriptscriptstyle 01} \cdot g_{\scriptscriptstyle 12}, A_{\scriptscriptstyle 01}+A_{\scriptscriptstyle 12},
\omega^2_{\scriptscriptstyle 01}+ \omega^2_{\scriptscriptstyle 12}).
\]

\emptycomment{
Composition of 1-morphisms is given
by multiplication in $\String(n)$ and addition of forms,
\[(g_1, A_1, \omega^2_1) \circ ( g_2, A_2, \omega^2_2) = (g_2\cdot
g_1, A_2+A_1, \omega^2_2+\omega^2_1).
\] }

\item{$\bstringnp(U)_2$}: a 2-simplex with edges $(g_{ij}, A_{ij},
  \omega^2_{ij})$ for $0\le i<j \le 2$, or equivalently (in the model of bicategory), a 2-morphism between $(g_{\scriptscriptstyle 01}, A_{\scriptscriptstyle 01}, \omega^2_{\scriptscriptstyle 01})
  \circ (g_{\scriptscriptstyle 12}, A_{\scriptscriptstyle 12}, \omega^2_{\scriptscriptstyle 12})$ and $(g_{\scriptscriptstyle 02}, A_{\scriptscriptstyle 02}, \omega^2_{\scriptscriptstyle 02}) $,  is given by a pair $(f,
  \omega^1)$ with $f\in C^\infty(U,
  \underline{U(1)})$ and $\omega^1 \in \Omega^1(U)$ such that
\[ A_{\scriptscriptstyle 12}-A_{\scriptscriptstyle 02}+A_{\scriptscriptstyle 01}= \omega^1- d \log f , \quad \omega^2_{\scriptscriptstyle 12} -
\omega^2_{\scriptscriptstyle 02} + \omega^2_{\scriptscriptstyle 01} = -d\omega^1  .
\] Moreover, $f$ gives rise to an isomorphism\footnote{Note that a $U(1)$-valued function on $U$ provides an isomorphism of bibundles from $U$ to $\String(n)$ via the map $\B U(1) \to \String(n)$. See also Lemma \ref{lem:lift-g}.} of bibundles  $g_{\scriptscriptstyle 01} \cdot g_{\scriptscriptstyle 12}
\Leftarrow g_{\scriptscriptstyle 02}$.
\[
\begin{tikzpicture}[>=latex',mydot/.style={draw,circle,inner
    sep=1pt},every label/.style={scale=1},scale=1]

  \foreach \i in {0}{
  \node[mydot, fill=black, label=240:$({\theta_0}{,}\, B_0)$]                at (-1+6*\i, -1.73)    (p\i0) {};
  \node[mydot,fill=black,label=90:$(\theta_1{,}\,B_1)$]      at (+6*\i,0) (p\i1) {};
  \node[mydot,fill=black,label=-60:$(\theta_2{,}\,B_2).$]     at (1+6*\i,-1.73)    (p\i2) {};
}
\begin{scope}[<-]
    %--fig 0
    \draw (p00) --node[left,scale=.8]{$(g_{\scriptscriptstyle 01}, A_{\scriptscriptstyle 01}, \omega^2_{\scriptscriptstyle 01})$} (p01);
    \draw (p00) --node[below,scale=.8]{$(g_{\scriptscriptstyle 02}, A_{\scriptscriptstyle 02}, \omega^2_{\scriptscriptstyle 02})$} (p02);
    \draw (p01) --node[right,scale=.8]{$(g_{\scriptscriptstyle 12}, A_{\scriptscriptstyle 12}, \omega^2_{\scriptscriptstyle 12})$} (p02);
\end{scope}
 \node at (0,-1){\tiny{$(f, \omega^1)$}};
 % \node at (20, -1.8){\footnotesize{$\Kan(2,2)$,}};
\end{tikzpicture}
\emptycomment{
\begin{split}
d \log f &= \omega^1 - A_0+A_1-A_2, \\
d\omega^1 &= - \omega^2_0 +
\omega^2_1 - \omega^2_2 .
\end{split} }
\]
\end{itemize}

Given a morphism $U\xrightarrow{\phi}V$, the associated functor
$\bstringnp(V)\to \bstringnp(U)$ is given by pre-compositions and
pullbacks of forms.

Now let us look at $ \holim \bstringnp ( U(M)_\bullet)$ for a cover
$\{U_i \}$ of $M$, here $U(M)_\bullet$ given in \eqref{defi:UM} is the nerve of the \v{C}ech groupoid associated
to the cover $\{U_i \}$. An object in $ \holim \bstringnp ( U(M)_\bullet)$
consists of
\begin{itemize}
\item $U_i \times \String(n) \to U_i$, $B_i \in \Omega^2(U_i)$,
  $\theta_i \in \Omega^1 (U_i, \so(n))$;
\item $g_{ij}: U_{ij} \to \String(n)$, $A_{ij} \in \Omega^1(U_{ij})$,
  $\omega^2_{ij} \in \Omega^2(U_{ij})$, such
  that
\begin{equation}\label{eq:a-omega3-theta}
B_j - B_i= d A_{ij} +\omega_{ij}^2, \quad
\cs_3(\theta_j) - \cs_3(\theta_i) = d\omega_{ij}^2, \quad
\theta_i - \ad_{\bg_{ij}} \theta_j = - \bg^*_{ij} \TM;
\end{equation}
\item $f_{ijk}: U_{ijk} \to U(1)$, $\omega^1_{ijk} \in
  \Omega^1(U_{ijk})$, such that $f_{ijk}$ is an
  isomorphism $g_{ik}\Rightarrow g_{ij}\cdot g_{jk}$, and
\begin{equation}\label{eq:f-omega1}
(\check{\delta} A)_{ijk} =  \omega^1_{ijk} - d \log f_{ijk} , \quad
 (\check{\delta} \omega^2)_{ijk}=-d\omega^1_{ijk} .
\end{equation}
\item a pentagon condition for 2-morphisms indicated by the following diagram,
\begin{equation}\label{eq:pent}
\xymatrix{
\big( (g_{ij}, A_{ij}, \omega^2_{ij}) \circ (g_{jk}, A_{jk},
\omega^2_{jk} ) \big) \circ (g_{kl}, A_{kl}, \omega^2_{kl})
\ar@{=>}[dd]^{(\id\circ_\h a, 0, 0)} \ar@{<=}[r]^{\qquad\qquad(f_{ijk}, \omega^1_{ijk})}&
(g_{ik}, A_{ik}, \omega^2_{ik}) \circ (g_{kl}, A_{kl}, \omega^2_{kl}) \ar@{<=}[rd]^{(f_{ikl},
  \omega^1_{ikl}) }
\\
& & (g_{il}, A_{il}, \omega^2_{il}) \\
 (g_{ij}, A_{ij}, \omega^2_{ij}) \circ \big( (g_{jk}, A_{jk},
 \omega^2_{jk}) \circ (g_{kl}, A_{kl}, \omega^2_{kl}) \big) \ar@{<=}[r]^{\qquad\qquad(f_{jkl}, \omega^1_{jkl})}  &  (g_{ij}, A_{ij}, \omega^2_{ij}) \circ (g_{jl}, A_{jl},
 \omega^2_{jl}), \ar@{<=}[ru]^{(f_{ijl},
  \omega^1_{ijl}) } }
\end{equation}
where $a$ is the associator of the string group
$\String(n)$, and $\circ_h$ is the horizontal composition of
2-morphisms, noticing that $(g_{ij} \circ g_{jk}) \circ g_{kl}$ and
$g_{ij} \circ (g_{jk} \circ g_{kl} )$ are composed 1-morphisms $U_{ijk}
\xrightarrow{ (g_{ij}, g_{jk}, g_{kl})} \String(n)^{\times 3}\to \String(n)$. According to Lemma \ref{lem:tech}, the 2-morphism
$\id\circ_h a$ is given by $U(1)$-valued functions $F_{ijkl}: U_{ijkl}\to U(1)$, which
converges to a class in the \v{C}ech cohomology group $ \check{H}^3(M, \Z)$ determined by the extension class $\half p_1$. Thus the above diagram
says exactly
\begin{equation}\label{eq:f-F-omega1}
\check{\delta} f =F, \quad \check{\delta} \omega^1 - \check{\delta}
d\log f =0, \quad \check{\delta} d \omega^1 = 0.
\end{equation}
The latter two equations are implied by
\eqref{eq:f-omega1}.
\end{itemize}

A 1-morphism  in $\holim \bstringnp(U(M)_\bullet)$ from $(U_i\times \String(n)\to U_i, \ttheta_i, \widetilde{B}_i;
\widetilde{g}_{ij}, \tA_{ij}, \widetilde{\omega}^2_{ij}; \widetilde{ f}_{ijk}, \widetilde{\omega}^1_{ijk})$ to $(U_i\times \String(n)\to U_i, \theta_i, B_i;
g_{ij}, A_{ij}, \omega^2_{ij}; f_{ijk}, \omega^1_{ijk})$ consists of
\begin{itemize}
\item $g_i: U_i \to \String(n)$, $A_i\in \Omega^1(U_i)$, $\omega^2_{i}
  \in \Omega^2(U_i)$, such that $(g_i, A_i,
  \omega^2_i) \in \bstringnp(\sqcup_i U_i)_1$ is a 1-morphism between
\[(U_i\times \String(n)\to U_i, \theta_i, B_i) \xleftarrow {(g_i, A_i,
  \omega^2_i)} (U_i\times\String(n)\to U_i, \ttheta_i,  \tB_i),  \]
that is,
\begin{equation} \label{eq:1mor0}
    \tB_i-B_i =  \omega^2_i +dA_i, \quad d\omega^2_i =
\cs_3(\ttheta_i)-\cs_3(\theta_i), \quad \theta_i-\ad_{\bar{g}_i}
\ttheta_i = - \bg_i^* \TM;
\end{equation}
\item an element $(f_{ij}, \omega^1_{ij} ) \in \bstringnp(\sqcup_{ij}
  U_{ij})_2$ which provides a 2-morphism making the following diagram 2-commutative
\[
\xymatrix{
( \theta_j,B_j) \ar[d]_{(g_{ij}, A_{ij}, \omega^2_{ij})} & (\ttheta_j,\tB_j
) \ar[l]_{(g_j, A_j, \omega^2_j)} \ar[d]^{(\tg_{ij}, \tA_{ij},
  \tomega^2_{ij})} \\
(\theta_i,B_i ) & (\ttheta_i,\tB_i), \ar[l]^{(g_i, A_i, \omega^2_i)}
}
\]
that is
\begin{equation}\label{eq:1mor1}
\xymatrix{
g_{ij} \cdot g_j & \ar @{=>}[l]_{f_{ij}\qquad\qquad\qquad\qquad\qquad\qquad\qquad\qquad} g_i \cdot \tg_{ij},  \quad
A_{ij} + A_j - (A_i + \tA_{ij}) = \omega^1_{ij} - d\log f_{ij}, \quad
\omega^2_{ij} + \omega^2_j -(\tomega^2_{ij} + \omega^2_i) =
-d\omega^1_{ij}.}
\end{equation}

\item a higher coherence condition between 2-morphisms
\[
\xymatrix{ & ( \theta_k,B_k) \ar[dl] \ar[dd]& & ( \ttheta_k,\tB_k) \ar[dl]
  \ar[dd] \ar[ll] \\
( \theta_j,B_j) \ar[dr] & & (\ttheta_j,\tB_j) \ar[dr] \ar[ll] \\
& (\theta_i,B_i) && ( \ttheta_i,\tB_i) \ar[ll]}
\]
which gives us the following equations\footnote{These two equations
  are deduced from the coherence condition between $g$'s and $A$'s, the one for
  $\omega^2$'s can be implied by the second equation in \eqref{eq:1mor1}.}
\[
f_{ijk}\tf_{ijk} ^{-1} = (\check{\delta} f_{..})_{ijk}, \quad
\tomega_{ijk}^1-d\log \tf_{ijk} -\omega^1_{ijk} + d\log f_{ijk} =
-(\check{\delta}\omega_{..}^1)_{ijk} + (\check{\delta} d \log f_{..})_{ijk},
\]
\end{itemize}

A 2-morphism in $\holim \bstringnp(U(M)_\bullet)$ from $(\tg_i, \tA_i,
\tomega^2_i ; \tf_{ij}, \tomega^1_{ij} )$ to $(g_i, A_i, \omega^2_i;
f_{ij}, \omega_{ij}^1) $ consists of
\begin{itemize}
\item an element $(f_i, \omega^1_i) \in \bstringnp(\sqcup U_i)_2 $
  from $(\tg_i, \tA_i, \tomega^2_i) $ to $(g_i, A_i, \omega^2_i)$, that is, $\xymatrix{g_i &
    \ar@{=>}[l]_{f_i}\tg_i} $, and
\[
A_i -\tA_i = -d\log f_{i} + \omega^1_i, \quad
\omega^2_{i}-\tomega^2_{i} =-d\omega^1_{i};
\]
\item a coherence condition held on $U_{ij}$,
\[
(f_{ij} , \omega^2_{ij}) \circ  (f_i , \omega_i^1) = (f_j, \omega^1_j) \circ (\tf_{ij}, \tomega_{ij}^2),
\]
that is\footnote{The same as before, the two equations
  are deduced from the coherence condition between $g$'s and $A$'s, and the one coming from
  $\omega^2$'s can be implied by the second equation.},
\[
f_{ij} \tf_{ij}^{-1}=  f_j f_i^{-1}, \quad \omega^1_{ij} - \tomega^1_{ij}  = (\check{\delta} \omega_{.}^1)_{ij}.
\]
\end{itemize}

Then $(\bstringnp)^+: \Mfd^{\op} \to \tGpd$ is a $(3,1)$-sheaf consists of
\begin{itemize}
\item[$\bullet$]~ $\bstringnpp(M)_0$: an object is a pair $(\{U_i\},
  P_{c})$, where $\{U_i\}$ is an open cover of $M$ and $P_{c}$
  is an element in $\holim \bstringnp(U(M)_\bullet)_0$;
\item[$\bullet$]~ $\bstringnpp(M)_1$: a 1-morphism between $(\{ U_i \},
  P_{c})$ and $(\{\widetilde{U}_i\},\widetilde{P}_{c})$ is a pair consisting of  a common refinement $\{V_i\}$ of
  $\{U_i\}$ and $\{\widetilde{U}_i\}$ and an element $\phi_{c} \in \holim \bstringnp(V(M)_\bullet)_1$;
\item[$\bullet$]~$\bstringnpp(M)_2 $: a 2-morphism between $(\{V_i \}, \phi_{c})$
  and $(\{\widetilde{V}_i,\widetilde{\phi}_c \})$ consists of a common refinement $\{
  W_i \} $ of $\{V_i \}$
  and $\{ \widetilde{V}_i \}$ and an element $\alpha_{c} \in \holim
  \bstringnp(W(M)_\bullet)_2$. Moreover, $(\{ W_i \}, \alpha_{c})$ and
  $(\{ \widetilde{W}_i \}, \widetilde{\alpha}_{c})$ are identified if  $\alpha_{c}$ and
  $\widetilde{\alpha}_{c}$ agree on a further
  common refinement of $\{W_i \}$ and $\{ \widetilde{W}_i \}$.
\end{itemize}

\begin{rmk}
For simplicity, we call an element $P_c \in \holim \bstringnp(U(M)_\bullet)_0$ a {\bf string data}. The construction of our string (pre)sheaf works for a general compact Lie group $G$ with the extension class $\pp$ a multiple of $p_1$, as long as we adjust  the coefficient in the front of the Chern-Simon $3$-form  by this multiple also.
\end{rmk}

\subsection{Lifting Theorem and Comparison}\label{sec:string-lift}

As we state in the introduction, there are already several ways to
grasp the concept of string structure. Redden's string
class is probably the most accessible and concise, while Waldorf's
method includes connection data and makes it easy and natural to locate the
integrity of the string class. The reason to develop yet another way
here, is to connect with
concepts involving differential forms, such as Courant algebroids,
descent equations, and Deligne cohomology. We found it also much
easier in this language to relate to physics literatures, such as
\cite{fss:characteristic-classes}.  Then, we own readers a
justification.

The direct comparison to previous methods might be a
wrong approach to see the nature of the problem since both sides
(especially our side) involves heavy machinery. We
remark (Remark \ref{rmk:compare-wo-conn}, \ref{rmk:curv},
\ref{rmk:compare-redden}, \ref{rmk:compare-waldorf}) carefully on links of these concepts, and
focus ourselves on the proof of lifting theorem from the viewpoint of Stolz-Teichner for justification.

Let $\mathsf{B}G_{c}$ be the $(2, 1)$-sheaf of $G$-principal bundles with
connections. We take the model  ${\mathsf{B}G^p_c}^+$ for it. Thus comparing the construction of $\bstringnp$ with Example \ref{ep:bgconn}, we see that  there is a natural projection ${\bstringnp}^+
\xrightarrow{\pi} \mathsf{B}\Spin(n)_c$, by forgetting higher data. More precisely, given an object $$P_c =( \sqcup_i
U_{i}; U_i\times \String(n)\to U_i, \theta_i, B_i;
g_{ij}, A_{ij}, \omega^2_{ij}; f_{ijk}, \omega^1_{ijk}) \in  {\bstringnp}^+,$$
where $(\theta_i,B_i;
g_{ij}, A_{ij}, \omega^2_{ij}; f_{ijk}, \omega^1_{ijk})$ satisfies
equations \eqref{eq:a-omega3-theta}, \eqref{eq:f-omega1} and
\eqref{eq:f-F-omega1}, since the isomorphism of bibundles between $g_{ij}\cdot
g_{jk} \Leftarrow g_{ik}$ is given by a $U(1)$-function   $f_{ijk}$,
the projected $\Spin(n)$-valued function, $\bg_{ij}: U_{ij}
\xrightarrow{g_{ij}} \String(n) \to \Spin(n)$, satisfies strict cocycle
condition $\bg_{ij}\cdot \bg_{jk} =\bg_{ik}$. This gives us a
$\Spin(n)$-principal bundle $\bar{P}$. Furthermore, equation \eqref{eq:a-omega3-theta} implies that $\theta_i$
provides a connection on $\bar{P}$.

\emptycomment{
Given an object $M\xrightarrow{\bar{P}_c  } \BSpin(n)_c$, the class $\half p_1(\bar{P}_c) \in \check{H}^3_{U_\bullet}(M, \underline{U(1)})$ constructed in Lemma \ref{lem:tech} converges exactly to $\half p_1(\bar{P} )\in H^3(M, \underline{U(1)})\cong H^4(M, \Z)$ when we refine the cover $\{U_i\}$, where $\bar{P}$ is the $\Spin(n)$-principal bundle that $\bar{P}_c$ glues to. In particular, if $\{U_i \}$ is a good cover, then $\half p_1(\bar{P}_c) = \half p_1(\bar{P}) \in \check{H}^3_{U_\bullet}(M, \underline{U(1)})=H^3(M, \underline{U(1)})$.}

%30.11.afternoon

\begin{thm} \label{thm:lift}
An object $M\xrightarrow{\bar{P}_c} \mathsf{B}\Spin(n)_c$ in
$\mathsf{B}\Spin(n)_c$ over a fine enough good cover $\{U_i\}$ of $M$ lifts to an object in $M\xrightarrow{P_c}{\bstringnp}^+$ in ${\bstringnp}^+$
\[
\xymatrix{ & {\bstringnp}^+ \ar[d] \\
M \ar@{.>}[ur] \ar[r]^{\bar{P}_c} & \mathsf{B}\Spin(n)_c,
}
\] if and only if $\frac{1}{2} p_1(\bar{P})=0$, where $\bar{P}$ is the $\Spin(n)$-principal bundle that $\bar{P}_c$ glues to.
\end{thm}
\begin{rmk}\label{rmk:lift}
This theorem says that as long as a good cover is fine enough, the obstruction for a $\Spin(n)$ data $\bar{P}_c$ to  lift to a $\String(n)$ data is $\half p_1(\bar{P}_c)$. Since we may always find fine enough good covers, this theorem justifies that our construction $\bstringnp$ is indeed
reasonable to be called the $(3,1)$-presheaf of $\String(n)$-principal
bundles.
\end{rmk}

We first prove some technical lemmas:

\begin{lemma}\label{lem:lift-g}
Given a function $\bar{g}: U \to G$, one may always lift it to a
morphism $g:U\to \String(G)$ if $\bar{g}^{-1}(G^{(1)}_\alpha)=U\times_{\bg,
  G, \pr} G^{(1)}_\alpha$'s are contractible, where $\pr: \sqcup G^{(1)}_\alpha \to
G$ is the covering map.
\end{lemma}
\pf
A morphism $g:U\to \String(G)$ is given by a bibundle $E$ which is a
$\sqcup G^{(1)}_{\alpha \beta} \times U(1) \Rightarrow \sqcup G^{(1)}_{\alpha}$
principal bundle over $U$. We know that the underlying morphism $\bg: U \to
G$ of $g$ is given by the bibundle
$U\times_{\bar{g}, G, \pr} \sqcup G^{(1)}_{\alpha}$, and $E$ is an
$U(1)$-bundle over it. Since $\bar{g}^{-1}(G^{(1)}_\alpha)$'s are contractible, $E\cong U\times_{\bar{g}_{ij}, G, \pr}
\sqcup G^{(1)}_\alpha \times U(1)$ as manifolds.      Suppose that the action is given by
\[ (x, g_\alpha, a) \cdot (g_{\alpha\beta}, a') := (x, g_\beta, a+a'
+ \lambda_{\alpha \beta} (x) ), \quad \text{for}\; (x, g_\alpha, a) \in E,\; \text{and} \;
(g_{\alpha\beta}, a') \in G^{(1)}_{\alpha \beta} \times U(1),  \]
for a certain function $\lambda_{\alpha\beta}: U \to U(1) $. The
associativity of the action is equivalent to the fact that $(\check{\delta}
\lambda)_{\alpha\beta\gamma} =
\eta_{\alpha\beta\gamma}(\bg(x))$. But $H^2(
U\times_{\bar{g}, G, \pr} \sqcup G^{(1)}_{\alpha}, \underline{U(1)}) =0$, thus 2-cocycle
$\bar{g}^*\eta$ is always exact. Therefore we may always find such
$\lambda$.
\qed\vspace{3mm}

We endow a simplicial hypercover $G^{(\bullet)}$ of $BG_\bullet$--the nerve of
$G$. Suppose that the extension class $\pp \in H^3(BG_\bullet, \underline{U(1)})\cong
H^4(BG_\bullet, \Z)$ is represented by  the $U(1)$-valued 3-cocycle $(\Theta,
\Phi, \eta, 0)$ supported on this cover. The last entry being $0$ is
implied by the closedness.
Notice that $\bg_{ij}: U_{ij} \to G $ extends to a simplicial morphism $\bg_{(\bullet)}$
from the simplicial nerve $U(M)_\bullet$ of the covering groupoid
$\sqcup U_{ij} \Rightarrow \sqcup U_i$ to
$BG_\bullet$. On each simplicial level $U(M)_k=\sqcup
U_{i_0 i_1\dots i_k} $, we endow it with the pullback cover of the one on
$BG_k$ pulled back by $\bg_{(k)}$. We may always start with a fine enough cover $\{U_i\}$ so that $\bar{g}_{ij}(U_{ij})$ is either entirely in $G_\alpha$ or does not intersect $G_\alpha$. Thus we may assume that both $\{U_i\}$ and the pullback covers are good. Then the
simplicial-\v{C}ech double complex \eqref{eq:db-cx-M} calculates the cohomology $H^\bullet(U(M)_\bullet, \underline{U(1)})\cong \check{H}^\bullet(M, \underline{U(1)})$. We denote the pullback cocycle $\bg_{(\bullet)}^*(\Theta, \Phi, \eta, 0)$ by
$(\bTheta, \bPhi, \bareta, 0)$ and it is a cocycle in double complex \eqref{eq:db-cx-M} representing a class $\pp(\bar{P}_c)\in \check{H}^3(M, \underline{U(1)})$.
% Not to be confused, the simplicial differential comes from the face maps of $U(M)_\bullet$,
%and the \v{C}ech differential comes from the covers on $U(M)_k$.
\begin{equation}\label{eq:db-cx-M}
\xymatrix{
 &&&&&&\\
 C(\sqcup U_{ijkl;s},\underline{U(1)})\ar[r]^{\check{\delta}}\ar[u]^{\delta}
 &C(\sqcup U_{ijkl;s,t},\underline{U(1)})\ar[r]^{\check{\delta}}\ar[u]^{\delta}
 &\dots %\Omega^3(U^1)\ar[r]^{\check{\delta}}\ar[u]^{-d}
&%\Omega^3(U^2)\ar[r]^{\check{\delta}}\ar[u]^{d}
&%\Omega^3(U^3)\ar[r]^{\check{\delta}}\ar[u]^{-d}
&
 \\
  C(\sqcup U_{ijk;p},\underline{U(1)})\ar[r]^{\check{\delta}}\ar[u]^{\delta}
  &C(\sqcup U_{ijk;p,q},\underline{U(1)})\ar[r]^{\check{\delta}}\ar[u]^{\delta}
 &C(\sqcup U_{ijk;p,q,r},\underline{U(1)})\ar[r]^{\check{\delta}}\ar[u]^{\delta}
&\dots %C(U^2)\ar[r]^{\check{\delta}}\ar[u]^{d}
&%\Omega^2(U^3)\ar[r]^{\check{\delta}}\ar[u]^{-d}
&
 \\
 C(\sqcup U_{ij;\alpha},\underline{U(1)})\ar[r]^{\check{\delta}}\ar[u]^{\delta}
 & C(\sqcup U_{ij;\alpha,\beta},\underline{U(1)})\ar[r]^{\check{\delta}}\ar[u]^{\delta}
 &C(\sqcup U_{ij;\alpha,\beta,\gamma},\underline{U(1)})\ar[r]^{\check{\delta}}\ar[u]^{\delta}
&C(\sqcup U_{ij;\alpha,\beta,\gamma,\delta},\underline{U(1)})\ar[r]^{\qquad\check{\delta}}\ar[u]^{\delta}
&\dots %\Omega^1(U^3)\ar[r]^{\check{\delta}}\ar[u]^{-d}
&
 \\
 C(\sqcup U_i ,\underline{U(1)})\ar[r]^{0}\ar[u]^{\delta}
 &C(\sqcup U_i, \underline{U(1)})\ar[r]^{id} \ar[u]^{\delta}
&C(\sqcup U_i, \underline{U(1)})\ar[r]^{0}\ar[u]^{\delta}
&C(\sqcup U_i , \underline{U(1)})\ar[r]^{id} \ar[u]^{\delta}
&\dots %C(\sqcup U_i , \underline{U(1)})\ar[r]^{\check{\delta}} \ar[u]^{-d}
&
}
\end{equation}

\begin{lemma}\label{lem:tech}
The $2$-morphism $(\id
\circ_\h a, 0, 0) $ in diagram \eqref{eq:pent} is given by
$U(1)$-valued functions $F_{ijkl}:U_{ijkl}\to U(1)$. Moreover, if the cover is fine enough,
$F_{ijkl}$'s give rise to a cocycle representing $\half p_1(\bar{P}_c) \in
\check{H}^3(M, \underline{U(1)})$.
\end{lemma}
\pf We continue to use the notation and a fine enough cover $\{U_i\}$ given just before this lemma. For us now $G=\Spin(n)$.
Since $\check{H}^2(\sqcup U_{ij}, \underline{U(1)})=0$,
$(\check{\delta} \bareta)=0$ implies that
$\bareta=-\check{\delta}\lambda$.  We continue this tic-tac-toe procedure, since
$\check{H}^{\ge 0}(\sqcup U_{ijk}, \underline{U(1)})=\check{H}^{\ge 0}(\sqcup U_{ijkl}, \underline{U(1)}) =0$,  we have
\begin{equation}\label{eq:lambda-varphi}
\bareta=-\check{\delta}\lambda, \quad \bPhi=\delta \lambda +
\check{\delta} \varphi, \quad \bTheta=\delta \varphi +F,
\end{equation}
and $F$ is a function $U_{ijkl}\to U(1)$.

A calculation shows that the bibundle of $(g_{ij} \circ g_{jk})\circ
g_{kl}$ is $$U_{ijkl}\times_{G^{\times 3}} \sqcup G^{(3)}_s \times
(U(1))^{\times 3} \times_{\sqcup G_{\alpha} } \sqcup G_{\alpha \beta} \times
U(1)/ \sqcup G^{(3)}_{s, t} \times (U(1))^{\times 3}, $$
 where the action is
given by\footnote{Recall that  $\{U_i
  \}$ is fine enough so that $\bar{g}_{ij}(U_{ij})$ is either entirely in $G_\alpha$
  or does not intersect $G_\alpha$. Then by Lemma
  \ref{lem:lift-g}, one may always take trivial
  bibundles for $g_{ij}$, therefore their various composites.}
\[
\begin{split}
& (x_{ijkl}, w_s, a_1, a_2, a_3, g_{\alpha,\beta}, a) \cdot (w_{s, t},
a'_1, a'_2, a'_3)  \\
=& (x_{ijkl}, w_t, a_1+a'_1 + \lambda_{\alpha_1,\gamma_1}(x_{ij}),
a_2+a'_2+\lambda_{\alpha_2,\gamma_2}(x_{jk}),
a_3+a'_3+\lambda_{\alpha_3,\gamma_3}(x_{kl}), g_{\gamma,\beta}, \\
&a-a'_1-a'_2-a'_3-d_2^*\Phi(w_{s,t})-d^*_0\Phi(w_{s,t})-\bareta_{\alpha,\beta,\gamma}(x_{il})
).
\end{split}
\] This bibundle is isomorphic to the bibundle $U_{ijkl}\times_{G} \sqcup G_\alpha
\times U(1)$ through
\[
[(x_{ijkl}, w_s, a_1, a_2, a_3, g_{\alpha,\beta}, a)] \xrightarrow{\psi}
[(x_{ijkl}, g_\beta,
a+a_1+a_2+a_3+\varphi_{p_0}(v_0)+\varphi_{p_2}(v_2)+\lambda_{\alpha,\beta}(x_{ijkl}) )],
\]
where the right action of $\sqcup G_{\alpha,\beta}\times U(1)
\Rightarrow \sqcup G_{\alpha}$ is given by $$(x_{ijkl}, g_\beta,
a)\cdot (g_{\beta,\beta'}, a' )=(x_{ijkl}, g_{\beta'},
a+a'+\lambda_{\beta,\beta'}(x)).$$ The morphism $\psi$ is well-defined
thanks to the second equation in \eqref{eq:lambda-varphi}, and it is a
bibundle isomorphism thanks to the first equation in
\eqref{eq:lambda-varphi}.

The bibundle of $g_{ij} \circ (g_{jk}\circ
g_{kl})$ is given by exactly the same form but  the quotient is
given by a different action,
\[
\begin{split}
& (x_{ijkl}, w_s, a_1, a_2, a_3, g_{\alpha,\beta}, a) \cdot (w_{s, t},
a'_1, a'_2, a'_3)  \\
=& (x_{ijkl}, w_t, a_1+a'_1 + \lambda_{\alpha_1,\gamma_1}(x_{ij}),
a_2+a'_2+\lambda_{\alpha_2,\gamma_2}(x_{jk}),
a_3+a'_3+\lambda_{\alpha_3,\gamma_3}(x_{kl}), g_{\gamma,\beta}, \\
&a-a'_1-a'_2-a'_3-d_1^*\Phi(w_{s,t})-d^*_3\Phi(w_{s,t})-\bareta_{\alpha, \beta, \gamma}(x_{il})
).
\end{split}
\] Similarly,  this bibundle is also isomorphic to the same bibundle $U_{ijkl}\times_{G} \sqcup G_\alpha
\times U(1)$ through
\[
[(x_{ijkl}, w_s, a_1, a_2, a_3, g_{\alpha,\beta}, a)] \xrightarrow{\psi'}
[(x_{ijkl}, g_\beta,
a+a_1+a_2+a_3+\varphi_{p_1}(v_1)+\varphi_{p_3}(v_3)+\lambda_{\alpha,\beta}(x_{il}) )].
\]

The 2-morphism $\id\circ_\h a$ is then to add $\Theta$ on the last $U(1)$
component, and is
explicitly given by
\[
(x_{ijkl}, w_s, a_1, a_2, a_3, g_{\alpha,\beta}, a) \mapsto (x_{ijkl},
w_s, a_1, a_2, a_3, g_{\alpha,\beta}, a-\Theta(w_s) ),
\] before quotient. This map is
equivariant with respect to the above two actions because
$\delta\Phi=\check{\delta} \Theta$, thus it descends to
the quotients.  Under the isomorphism $\psi$ and $\psi'$, this
2-morphism is then given by
\[
\psi' (\id\circ_\h a (\psi^{-1}(x_{ijkl}, g_\beta, a)))=(x_{ijkl},
g_{\beta}, a+F_{ijkl}(x)),
\]
guaranteed by the last equation of \eqref{eq:lambda-varphi}. Moreover,
$[(F, 0, 0, 0)]=[ (\bTheta, \bPhi, \bareta, 0)] = \half p_1(\bar{P}_c)
\in \check{H}^3_{U_\bullet}(M, \underline{U(1)})$. Thus $F$ is a representative of $
\half p_1(\bar{P}_c)$.
\qed\vspace{3mm}

Now we are ready to prove Theorem \ref{thm:lift}.

%05.12.2016 morning
\pf
If there is a lifting object $P_c \in \bstringnp$ over an object $M \xrightarrow{\bar{P}_c} \BSpin(n)_c$ over $\{ U_i\}$ in $\BSpin_c$, we write $$P_c =( \sqcup_i
U_{i}; U_i\times \String_\pp(n)\to U_i, \theta_i, B_i;
g_{ij}, A_{ij}, \omega^2_{ij}; f_{ijk}, \omega^1_{ijk}),$$ where $(\theta_i,B_i;
g_{ij}, A_{ij}, \omega^2_{ij}; f_{ijk}, \omega^1_{ijk})$ satisfies
equations \eqref{eq:a-omega3-theta}, \eqref{eq:f-omega1} and
\eqref{eq:f-F-omega1}. Then \eqref{eq:f-F-omega1} and Lemma
\ref{lem:tech} implies that $\half p_1(\bar{P}_c)=0$.

\emptycomment{
%proof without tech lemma, can only prove for p1 tensor real =0
Since the 2-morphism between $g_{ij}\circ
g_{jk} \Leftarrow g_{ik}$ is given by an $U(1)$-function   $f_{ijk}$,
the projected $G$-valued function $\bg_{ij}$ satisfies strict cocycle
condition $\bg_{ij}\circ \bg_{jk} =\bg_{ik}$. This gives us a
principal bundle which is exactly the principal bundle $\bar{P}$ that
we start with.  Further equation \eqref{eq:a-omega3-theta} implies that $\theta_i$
provides a connection for $\bar{P}$. Moreover, the curvature 2-form $R_i :=
d\theta_i +*[\theta_i, \theta_i]$ is a global 2-form in $\Omega^2(M,
Ad_\g)$. Moreover,
\[d cs_3(\theta_i) = (R_i, R_i)^\g =:\alpha \]
is a global 4-form in $\Omega^4(M)$.
According to the discussion in Section
\ref{sec:deligne-coh}, the 4-form $\alpha$ represents $\frac{1}{2}
p_1(\bar{P}_c)$ and  \eqref{eq:a-omega3-theta} \eqref{eq:f-omega1}
\eqref{eq:f-F-omega1} implies that $\frac{1}{2}
p_1(\bar{P}_c)=0$.}

For the other direction, we first do some preparation: given an object $M \xrightarrow{\bar{P}_c} \BSpin(n)_c$ over a good cover $\{ U_i\}$ of $M$ in $\BSpin(n)_c$, we
take a $U(1)$-cocycle $F$ representing $\half p_1(\bar{P}_c) \in \check{H}^3_{U_\bullet}(M,
\underline{U(1)})=H^3(M,
\underline{U(1)})$. Let
\begin{equation}\label{eq:Dp}
D_m=\underline{U(1)}\xrightarrow{d\log}\Omega^1\xrightarrow{d}\Omega^2\xrightarrow{d}\dots \xrightarrow{d}\Omega^m
\end{equation}
be the Deligne sheaf of depth $m$. Recall that the Deligne cohomology $H^\bullet(M,
D_3)$ is then the limit of the total cohomology of the following double complex
taking over all covers $\{U_i\}$ of $M$,
\[
\Big( C(U(M)_\bullet, \underline{U(1)}) \xrightarrow{d\log}
\Omega^1(U(M)_\bullet) \xrightarrow{d} \Omega^2(U(M)_\bullet )
\xrightarrow{d} \Omega^3(U(M)_\bullet), \; \;\check{\delta} \Big).
\] Then the general theory of Deligne cohomology tells us that there
is a surjective morphism
$
H^\bullet(M, D_3) \xrightarrow{\pi} H^\bullet(M, \underline{U(1)})$
given by forgetting the part of differential forms of a
Deligne cocycle. Moreover, there is a morphism $H^3(M, D_3)\xrightarrow{\bar{d}}
\Omega^4_{\cl}(M)$ to closed 4-forms by $[(F, \omega^1, \omega^2, \omega^3)] \mapsto
d\omega^3$. Here notice that $\omega^3=(\omega^3_i)$ is made up by
local 3-forms on each $U_i$. However, $d\omega^3_i$'s glue together to a closed
global 4-form, which we denote by $d\omega^3$, and it is independent of the choice of the Deligne
cocycle.    The above two morphisms fit into the following commutative diagram:
\begin{equation}\label{eq:deligne}
\xymatrix{H^3(M, D_3) \ar[d]^\pi \ar[r]^{\bar{d}} &\Omega_{\cl}^4(M) \ar[r] & H^4_{dR}(M, \R) \\
H^3(M, \underline{U(1)}) \ar[r] & H^4(M, \Z) \ar[ur]^{\otimes \R} }
\end{equation}

Since the natural morphism $H^3(M, D_3) \to H^3(M, \underline{U(1)})$ is surjective, we lift $F$ to a Deligne cocycle $(F, \omega^1, \omega^2, \omega^3) $, that is
\begin{equation}\label{eq:F-omegas}
\check{\delta} \omega^1 - d\log F = 0, \quad \check{\delta}\omega^2 + d \omega^1 = 0, \quad \check{\delta}\omega^3- d\omega^2 = 0.
\end{equation}
Now we adjust $\omega^3_i$ to be $\cs_3(\theta_i)$, where $\theta_i$ is the connection data with respect to the cover $\{U_i\}$ in $\bar{P}_c$. Both $d\omega^3_i$ and $d\cs_3(\theta_i)$ give to closed global 4-forms, and both represent the deRham classes $ \half p_1 \otimes \R$. Thus $\omega^3_i - \cs_3(\theta_i) = \gamma+d\beta_i$, where $\gamma\in \Omega^3(M)$ and $\beta_i\in\Omega^2(U_i)$. Then it is easy to verify that the Deligne class $(F, \omega_1, \omega_2+ \check{\delta} \beta, \cs_3(\theta_i) )$ is a lift of $F$. Thus we can begin with a Deligne cocycle $(F, \omega^1, \omega^2, \omega^3) $ with $\omega^3_i=\cs_3(\theta_i)$.

Now we construct a lift of $\bar{P}_c$ with respect to the cover $\{U_i\}$ under the condition $\half
p_1(\bar{P}_c)=0$. Fix a good cover $G^{(\bullet)}$ of $\BSpin(n)_\bullet$ as in the construction of $\String(n)$ in Subsection \ref{sec:string-gp}. Refining $\{U_i\}$ if necessary,  we may assume that $\bar{g}_{ij} (U_{ij})$ either falls entirely into  $G^{(1)}_{\alpha}$ or does not intersect $G^{(1)}_\alpha$. Then the condition in Lemma \ref{lem:lift-g} is naturally fulfilled. Thus, there is no obstruction to lift the transition functions $\bar{g}_{ij}$ for $\bar{P}_c$ to $g_{ij}: U_{ij}\to \String(n)$.  Since  $\half
p_1(\bar{P}_c)=0$, we may take a primitive $f$ of $F$, that is $F=\check{\delta} f$. Since $d\log F = \check{\delta}\omega^1$, $\omega^1 = \check{\delta} A + d\log f$ for some $A=(A_{ij}) \in \Omega^1(\sqcup U_{ij})$. We continue such a tic-tac-toe process, and find
\begin{equation} \label{eq:fAB}
F=\check{\delta} f, \quad \omega^1 = \check{\delta} A + d\log f, \quad \omega^2=\check{\delta}B - dA,
\end{equation}
for $f=(f_{ijk}) \in \underline{U(1)}(\sqcup U_{ijk})$,  $A=(A_{ij}) \in \Omega^1(\sqcup U_{ij})$, and $B=(B_i) \in \Omega^2(\sqcup U_i)$.  Both $g_{ij}\circ g_{jk}$ and $g_{ik}$ are morphisms from $U_{ijk}$ to $\String(n)$. They are presented by isomorphic bibundles from the discrete groupoid $U_{ijk}\Rightarrow U_{ijk}$ to  $\Gamma[\eta]$ with a similar construction to that of $\psi$ in Lemma \ref{lem:tech}. Then as in the proof of Lemma \ref{lem:lift-g}, we see that a $U(1)$-valued function $f_{ijk}$ on $U_{ijk}$ serves as a 2-morphism $g_{ik}\Rightarrow g_{ij} \circ g_{jk}$ because $U(1)\Rightarrow pt$ is a subgroupoid of $\Gamma[\eta]$ and sits in the center of it.  Thus, $( \sqcup_i
U_{i}; U_i\times \String(n)\to U_i, \theta_i, B_i;
g_{ij}, A_{ij}, \omega^2_{ij}; f_{ijk}, \omega^1_{ijk}) $ is a lift of $\bar{P}_c$. \qed\vspace{3mm}

%above controlled 21.10.2016

As shown in \cite{stolz-teichner:elliptic-obj}, if a $\Spin(n)$-principal bundle $\bar{P}$ admits string classes then the possible choices of the string classes form a torsor of $H^3(M, \Z)$. Then later in
\cite{waldorf:string-conn}, the author  further showed that for a fixed Chern-Simon 2-gerbe over $\bar{P}$, the choices of trivialisations modding out isomorphisms correspond exactly to string classes on $\bar{P}$.

We see that inside an object $P_c \in \bstringnp$, the determining information, is a covering $\{U_i\}$ together with $(\theta_i, B_i; g_{ij}, A_{ij}, f_{ijk})$, other terms $(F_{ijkl}, \omega^1_{ijk}, \omega^2_{ij}, \omega^3_{i}=\cs_3(\theta_i))$ are determined by these terms through \eqref{eq:fAB}. These terms $(F_{ijkl}, \omega^1_{ijk}, \omega^2_{ij}, \omega^3_{i}=\cs_3(\theta_i))$ representing a refinement of $\half p_1(\bar{P}_c) \in H^3(M, \underline{U(1)})$, may be viewed as the information for a given Chern-Simon 2-gerbe over $\bar{P}_c$ and its connections. Thus, after adding connection data inside, we may ask ourselves again how many string data lift $\bar{P}_c \in \BSpin(n)_c$  if we fix a choice of the cocycle $(F_{ijkl}, \omega^1_{ijk}, \omega^2_{ij}, \omega^3_{i}=\cs_3(\theta_i))$.

The Deligne cohomology group $H^2(M, D_2)$ may be viewed as a refinement of $H^2(M, \underline{U(1)})\cong H^3(M, \Z)$, where $D_2$ is the Deligne sheaf defined in \eqref{eq:Dp}. We have the following diagram (see also \eqref{eq:deligne})
\[\xymatrix{ \ker \pi_2 \cap \ker \bar{d} = \ker \pi_3  \ar[d]\ar[r] & H^2(M, D_3)=\ker \bar{d} \ar[d]\ar[r]^{\pi_3~\quad} & H^2(M, \underline{U(1)}) \cong H^3(M, \Z) \\
\ker \pi_2 \ar[r]& H^2(M, D_2) \ar[d]^{\bar{d}} \ar[ur]^{\pi_2} & \\
& \Omega_{\rm cl}^3(M) &
}
\]
Since we have
\begin{eqnarray*}H^3(M, \Z)=H^2(M, D_3)/\ker \pi_3
\hookrightarrow H^2(M, D_2) /\ker\pi_3 \to H^3(M, \Z)=H^2(M, D_2)/\ker \pi_2,
\end{eqnarray*}
and $H^3(M, \Z)$ may be viewed naturally both as a subgroup and a quotient of our group, we show that different lifts of $\bar{P}_c$ modding out isomorphisms is a torsor of $H^2(M, D_2)/\ker \pi_3 $  in the next theorem.

%07.12.2016 morning

\begin{thm}\label{thm:torsor}
Given an object $M\xrightarrow{\bar{P}_c} \BSpin(n)_c$ over a good cover $\{ U_i\}$ of $M$ in $\BSpin(n)_c$, and a fixed Deligne cocycle $(F_{ijkl}, \omega^1_{ijk}, \omega^2_{ij}, \cs_3(\theta_i))$ representing a refinement of $\half p_1 (\bar{P}_c) \in H^3(M, \underline{U(1)})$,  let us denote the set of all possible $P_c$'s lifting $\bar{P}_c$ with fixed $(F_{ijkl}, \omega^1_{ijk}, \omega^2_{ij}, \cs_3(\theta_i))$ by $S_{(F_{ijkl}, \omega^1_{ijk}, \omega^2_{ij}, \cs_3(\theta_i))}$. Then the Deligne cohomology group $H^2(M, D_2)$ acts on
\[ S_{(F_{ijkl}, \omega^1_{ijk}, \omega^2_{ij}, \cs_3(\theta_i))}/_{\text{$1$-morphisms}}. \]
This action descends to the quotient $H^2(M, D_2)/\ker \pi_3 $ and makes $ S_{(F_{ijkl}, \omega^1_{ijk}, \omega^2_{ij}, \cs_3(\theta_i))}/_{\text{$1$-morphisms}}$ a  $(H^2(M, D_2)/ \ker \pi_3 )$-torsor.
\end{thm}
\pf
Given a cocycle $(f^h, A^h, B^h)$ representing an element in $H^2(M, D_2)$, it acts on $S_{(F_{ijkl}, \omega^1_{ijk}, \omega^2_{ij}, \cs_3(\theta_i))}$ by
\begin{equation}\label{eq:fhahbh-action}
(U_i, \theta_i, B_i; g_{ij}, A_{ij}, f_{ijk}) \xrightarrow{\cdot (f^h, A^h, B^h)} (U_i, \theta_i, B_i+B_i^h; g_{ij}, A_{ij}+A^h_{ij}, f_{ijk}+f^h_{ijk}).
\end{equation}
If $(f^h, A^h, B^h)=D(\varphi, \alpha)$, then the above two elements in $S_{(F_{ijkl}, \omega^1_{ijk}, \omega^2_{ij}, \cs_3(\theta_i))}$ are connected by an isomorphism $(1, -\alpha, 0; \varphi, 0)$. Thus \eqref{eq:fhahbh-action} gives rise to an action of $H^2(M, D_2)$ on $S_{(F_{ijkl}, \omega^1_{ijk}, \omega^2_{ij}, \cs_3(\theta_i))}/_{\text{1-morphisms}}$.

Now we prove that the action of a cocycle $(f^h, A^h, B^h)$ induces an isomorphism if and only if  $[(f^h, A^h, B^h)] \in \ker \pi_2 \cap \ker \bar{d}$. This then will complete the proof of our statement.

First of all, the isomorphism induced by $(f^h, A^h, B^h)$ is possibly given through another  finer cover $\{V_i\}$. However as we may always pull back our cocycle to $V_i$, we might as well assume that $V_i=U_i$.
By a direct calculation,  $(f^h, A^h, B^h)$ induces an isomorphism $(1, A_i, \omega^2_i; f_{ij}, \omega^1_{ij})$,  if and only if
\begin{equation}\label{eq:iso-h}
\begin{split}
B^h_i = \omega^2_i+ dA_i, \quad d\omega^2_i=0, \quad A^h_{ij}=\check{\delta} (A.)_{ij} - \omega^1_{ij} + d\log f_{ij}, \quad \check{\delta}\omega^2. = - d\omega^1_{ij}, \quad (\check{\delta} f..)_{ijk} = (f^h_{ijk})^{-1}.
\end{split}
\end{equation}
Thus one direction is clear. If $(f^h, A^h, B^h)$ induces an isomorphism, then $dB^h_i=0$ and $f^h$ is a coboundary, which exactly shows that $d[(f^h, A^h, B^h)]=0$ and $\pi_2([(f^h, A^h, B^h)])=0$ respectively.

For the other direction,  if $dB^h_i=0$ and $f^h_{ijk}=(\check{\delta} f..)_{ijk}$ for some $f..$, then $(1, 0, B^h_i; f_{ij}, -A^h_{ij}+d\log f_{ij})$ gives us a desired isomorphism.
\qed
\begin{rmk}\label{rmk:compare-wo-conn}
If we forget the connection data inside a string data $P_c$ and only remember $(U_i; g_{ij}; f_{ijk})$, then the action of $H^2(M, D_2)$ simplifies to that of $H^3(M, \Z)$ through $\pi_2$. Thus we recover the structure of the torsor in \cite{stolz-teichner:elliptic-obj} for the string structures over a $\Spin(n)$-principal bundle via the projection $H^2(M, D_2)/\ker \pi_3\to H^3(M, \Z)$.
\end{rmk}

For a string data $P_c= (U_i\times \String(n), \theta_i, B_i; g_{ij}, A_{ij}; f_{ijk})$, we see that $\check{\delta}(\cs_3(\theta_i) - dB_i)=0$ by \eqref{eq:a-omega3-theta}, thus $\{\cs_3(\theta_i) - dB_i\}$ give rise to a global 3-form $H$ on $M$.  We define $H$ to be the {\bf curvature} of $P_c$. By \eqref{eq:1mor0}, we see that isomorphic string data over the same $\bar{P}_c$ has the same curvature. Notice that when $(f^h, A^h, B^h)$ acts on $P_c$, the curvature is changed by $dB^h_i$ which glues to a global closed 3-form $\bar{d}(f^h, A^h, B^h)\in \Omega_{\rm cl}^3(M)$. Thus we have the following corollary,

\begin{cor}\label{cor:curv}
Given an object $M\xrightarrow{\bar{P}_c} \BSpin(n)_c$ over $M$, and a Deligne cocycle $(F_{ijkl}, \omega^1_{ijk}, \omega^2_{ij}, \cs_3(\theta_i))$, the curvatures of all possible $P_c$'s lifting $\bar{P}_c$ with fixed $(F_{ijkl}, \omega^1_{ijk}, \omega^2_{ij}, \cs_3(\theta_i))$ form a torsor of $\im \bar{d}$.
\end{cor}
\begin{rmk}\label{rmk:curv} %we may put it afterwards with surjectivity
We notice the following commutative diagram
\begin{equation}\label{eq:d2-z}
\xymatrix{
H^2(M, D_2) \ar[r]_{\pi_2} \ar[d]_{\bar{d}} & H^3(M, \Z) \ar[d]^{\otimes \R} \\
\Omega^3_{\rm cl}(M) \ar@{->>}[r] & H^3(M, \R).
}
\end{equation}
One may interpret $\im \bar{d}$ as ``integral'' forms. Certainly, $\bar{d}$ is not always surjective.
\end{rmk}

\emptycomment{
this proposition may be viewed as a refinement of the combination of the above two statements. However, the freeness of the action may not be achieved. The reason is following: if two objects are connected by an isomorphism $(1, A_i, \omega^2_i; f_{ij}, \omega_{ij}^1)$, then they must be connected by the cocycle $(f^h, A^h, B^h)$ satisfying the following equations
\[
B_i^h = \omega_i^2 + dA_i, \quad A^h_{ij}= A_j-A_i-\omega_{ij}^1 + d\log f_{ij}, \quad f^h_{ijk} = (\check{\delta} f_{..})_{ijk}.
\]
If $(f^h, A^h, B^h)=D(\varphi, \alpha)$, combining with these equations, we have
\[
(\omega^2_i, -\omega_{ij}^2, 1) = D(\alpha_i-A_i, \varphi_{ij} \cdot f_{ij}^{-1}).
\]
From the condition of 1-morphism, it is not hard to see that $D(\omega^2_i, -\omega_{ij}^1, 1)=0$. However, we still have $H^1(M, D_2)$ which surjectively maps to $H^2(M, \Z)$ as an obstruction of   $(\omega^2_i, -\omega_{ij}^1, 1)$ being exact. ***signs and comparing to the cases before.
}

\begin{rmk}\label{rmk:compare-redden}
We conjecture that when we glue the local data of a string data, we
will obtain an $S^1$-gerbe over the underlying $\Spin(n)$-principal bundle
$\bar{P}$. This gives us the access to Redden's string class. For this
problem, a possible way is to apply the descent for $n$-bundles of
Wolfson \cite[Theorem 5.7]{Wolfson} to realize the gluing process. The difficulty would
lie further on gluing differential forms to obtain an
integral form which presents the string class. Notice that already in Redden's
thesis, there was no explicit formula to adjust to an
integral form for the string class. We thus leave it for future
investigation.
\end{rmk}
\begin{rmk}\label{rmk:compare-waldorf}
To compare thoroughly with Waldorf's method using trivialization of Chern-Simons
$2$-gerbe is neither a simple task. Both definitions use heavy
machineries, one with bundle gerbe theory, and the other with stack
theory. Nevertheless, some traces of equivalence are rather
visible\footnote{Here we thank Konrad Waldorf for very
  helpful conversations.}. In \cite[Definition 2.2.1]{waldorf:string-conn}, if we
take $Y$ to be $\sqcup_{i}U_i$, then bundle gerbe $\mathcal{S}$
characterised by a closed $3$-form,  over
each $U_i$,
corresponds to $2$-form $B_i$ in a string data, since a closed form is locally
exact. Similarly the bundle gerbe $\mathcal{P}$ corresponds to
$\omega^2_{ij}$,  isomorphism $\mathcal{A}$ corresponds to $A_{ij}$,
$\mathcal{M}$ corresponds to $\omega^1_{ijk}$, and $\sigma$ corresponds to $f_{ijk}$. Then
various coherence conditions correspond to our descent equations for
differential forms.
\end{rmk}

\emptycomment{
representing $\half p_1 (\bar{P}_c) \in H^3(M, U(1)\xrightarrow{d\log} \Omega^1
\xrightarrow{d} \Omega^2 \xrightarrow \Omega^3 )$, and $\omega^3$ is
exactly the Chern-Simon 3-form for the connection $\theta_i$ of
$\bar{P}_c$. ****(here I think
Deligne cohomology $H^3(M, U(1)\xrightarrow{d\log} \Omega^1
\xrightarrow{d} \Omega^2 \xrightarrow \Omega^3 )=H^3(M, \underline{U(1)})$ because
we can use tic-tac-toe, but need to check it for sure... in Chris's
article, he only mentions that it's a surjection. but I think for
manifold it's the same). ***(also why can we choose $\omega^3$ exactly
Chern-Simon 3-from?? because Chern-Simon 3-form and other choices
differ by a closed 3-form (or exact locally) because both of gives
arise to $\alpha=\langle R,  R \rangle$). *** here we need the section
on Deligne cohomology. *** Since  $\half
p_1(\bar{P}_c)=0$, there is a primitive Deligne cochain $(f, A, B
)$ such that $D(f, A, B) = (F, \omega^1, \omega^2, \omega^3)$.  We
thus construct a lift accordingly. ***Notice that there is no
obstruction to lift $g$ from $\bg$. }

\section{$(2,1)$-sheaf $\tcalgdpp$ of transitive Courant algebroids with connections}

The notion of a Courant algebroid was introduced in \cite{LWXmani}. See also \cite{roy:thesis,royt,rw,Severa:3-form} for various other aspects of Courant algebroids.

\begin{defi}\label{defi:ca}
A {\bf Courant algebroid} is a vector  bundle $C$ together with a bundle map $\rho:C\longrightarrow TM$, a nondegenerate symmetric  bilinear form  $\pair{-,-}$,
and an operation $\Courant{-,-}:\Gamma(C)\times \Gamma(C)\longrightarrow \Gamma(C)$ such that for all $e_1,e_2,e_3\in{\Gamma(C)}$, the following axioms hold:
\begin{itemize}
\item[\rm(i)] $(\Gamma(C),\Courant{-,-})$ is a Leibniz algebra;
\item[\rm(ii)] $\pair{\Courant{e_{1},e_{1}}, e_2}  = \half\rho(e_2)\pair{e_{1},e_{1}}$;
\item[\rm(iii)]
$\rho(e_{1})\pair{e_{2},e_{3}}=\pair{\Courant{e_{1},e_{2}},e_{3}}+\pair{e_{2},\Courant{e_{1}, e_{3}}}$.
  \end{itemize}
  \end{defi}

A Courant algebroid $(C,\Courant{-,-},\pair{-,-},\rho)$ is called {\bf transitive} if $\rho$ is surjective, that is, $\im \rho=TM$. A transitive Courant algebroid is an extension of a transitive Lie algebroid. However, not every transitive Lie algebroid $A$ admits such a Courant extension. The obstruction is given by the first Pontryagin class \cite{Bressler:Pclass,ChenRCA,sev15}.  See also \cite{MK2,Vai05} for more details about transitive Courant algebroids. In this section, we introduce the $(2,1)$-presheaf $\tcalgdp$ of transitive Courant algebroids with
connections, and we use the plus construction to sheafify it to a $(2,1)$-sheaf ${\tcalgdp}^+$. We then reinterpret the extension obstruction as the lifting obstruction,
\begin{equation}
\xymatrix{
&\tcalgdp^+\ar[d]\\
M \ar@{.>}[ur]\ar[r]&\tlp^+.
}
\end{equation}
Here $\tlp$ and $\tlp^+$ are the $(2,1)$-presheaf and $(2,1)$-sheaf of transitive Lie algebroids with connections respectively (see appendix \ref{app:tl}).
Notice how similar it is to the lifting story on the string side.

Given a transitive Courant algebroid, we have the following two short exact sequences:
\begin{eqnarray}
  \label{seq:tc1}&0\longrightarrow\ker\rho\longrightarrow C\stackrel{\rho}{\longrightarrow}TM\longrightarrow 0,\\
 \label{seq:tc2} &0\longrightarrow(\ker\rho)^\bot\longrightarrow \ker\rho\stackrel{\rho'}{\longrightarrow}\huaG\longrightarrow 0,
\end{eqnarray}
where  $\huaG=\ker\rho/(\ker\rho)^\bot$ is a Lie algebra bundle, whose fiber is isomorphic to a quadratic Lie algebra $(\g,(-,-)^\g)$. We will also use $(-,-)^\g$ to denote the fiberwise metric on $\huaG$.   A {\bf connection} of a transitive Courant algebroid $C$ consists of the following data:
 \begin{itemize}
   \item[$\bullet$] an isotropic splitting $s:TM\longrightarrow C$ of the short exact sequence \eqref{seq:tc1};
   \item[$\bullet$] a splitting $\sigma_s:\huaG\longrightarrow\ker\rho$ of the short exact sequence \eqref{seq:tc2} that is orthogonal to $s(TM)$ in $C$, i.e. $\pair{s(X),\sigma_s(a)}=0$ for all $X\in \Gamma(TM)$ and $a\in \Gamma(\huaG)$.
 \end{itemize}
 In \cite{ChenRCA}, the authors show that splittings $s$ and $\sigma_s$ always exist. A connection gives rise to an isomorphism $C\cong TM\oplus \huaG\oplus T^*M$ between vector bundles. Transferring the Courant algebroid structure on $C$ to $TM\oplus \huaG\oplus T^*M$, we obtain  the transitive Courant algebroid  $(TM\oplus \huaG\oplus T^*M,\Courant{-,-}^T_{\nabla,R,H},\pair{-,-}^T,\pr_{TM})$, which is determined by a connection $\nabla$ on $\huaG$, a  2-form $R \in \Omega^2(M, \huaG)$ and a 3-form $H\in \Omega^3(M)$, which obey a set of conditions given in \cite[Propositoin 2.2]{ChenRCA}. Here the bracket $\Courant{-,-}^T_{\nabla,R,H}$ and the pairing $\pair{-,-}^T$ are defined by
\begin{eqnarray}
    \nonumber\Courant{X+a+\xi,Y+b+\eta}^T_{\nabla,R,H}&=&[X,Y]+\nabla_Xb-\nabla_Ya+[a,b]_\huaG+R(X,Y)\\
  &&\label{eq:Tbracket}+L_X\eta-\ii_Yd\xi+P(a,b)-2Q(X,b)+2Q(Y,a)+H(X,Y),\\
  \label{eq:Tpair} \pair{X+a+\xi,Y+b+\eta}^T&=&\half\big(\xi(Y)+\eta(X)\big)+(a,b)^\g,
\end{eqnarray}
where $P:\Gamma(\huaG)\otimes \Gamma(\huaG)\longrightarrow\Omega^1(M)$ and $Q:\frkX(M)\otimes \Gamma(\huaG)\longrightarrow\Omega^1(M)$ are given by
\begin{eqnarray}
  P(a,b)(Y)&=&2(b,\nabla_Ya)^\g,\\
  Q(X,a)(Y)&=&(a,R(X,Y))^\g.
\end{eqnarray}

%21.11 Monday, evening

In particular, if $\huaG$ is the trivial bundle $M\times\g$ and the connection is given by $\nabla_Xa=X(a)$, we obtain the {\bf standard transitive Courant algebroid} structure on $TM\oplus (M\times \g)\oplus T^*M$ with the Courant bracket given by
\begin{equation}\label{eq:TSb}
 \Courant{X+a+\xi,Y+b+\eta}^T_S=[X,Y]+ X(b)-Y(a)+[a,b]_\g+L_X\eta-\ii_Yd\xi+\huaP(a,b),
\end{equation}
where $\huaP:\Gamma(M\times \g)\otimes \Gamma(M\times \g)\longrightarrow\Omega^1(M)$ is given by
$$
 \huaP(a,b)(Y)=2(b,Y(a))^\g.
$$For simplicity, for an object $U\in \Mfd$, we write
\[\tagu := TU\oplus(U\times \g)\oplus T^*U.\]
According to \cite[Proposition 2.7]{ChenRCA}, automorphisms of the standard transitive Courant algebroid are given as follows.

\begin{cor}\label{pro:auto}
  An automorphism of the standard transitive Courant algebroid $(\tagof{M},\Courant{-,-}^T_S,\pair{-,-}^T,\pr_{TM})$, where $\Courant{-,-}^T_S$ and $\pair{-,-}^T$ are given by \eqref{eq:TSb} and \eqref{eq:Tpair} respectively,  is of the form $\left(\begin{array}{ccc}
1&0&0\\
\phi&\tau&0\\
\beta&-2\phi^\star\tau&1
\end{array}\right)$, where $\tau$ is an orthogonal automorphism of the bundle of quadratic Lie algebras $M\times \g$ and $\phi:TM\longrightarrow M\times \g$ and $\beta:TM\longrightarrow T^*M$ are bundle maps satisfying the following compatibility conditions:
\begin{eqnarray}
  \label{eq:auto1}\half\big(\beta(X)(Y)+\beta(Y)(X)\big)+(\phi(X),\phi(Y))^\g&=&0,\\
   \label{eq:auto2}\tau(X(b))-X(\tau(b))-[\phi(X),\tau(b)]_\g&=&0,\\
   \label{eq:auto3}d\phi+\half[\phi,\phi]_\g&=&0,\\
   \label{eq:auto4}L_X(\beta(Y))-\ii_Yd(\beta(X))-\beta([X,Y])+\huaP(\phi(X),\phi(Y))&=&0.
\end{eqnarray}
Here $\phi^\star:M\times\g\longrightarrow T^*M$ is defined by
\begin{equation}\label{eq:star}
    \phi^\star(a)(X)=(a,\phi(X))^\g,\quad \forall a\in\Gamma(M\times \g), X\in\Gamma(TM).
\end{equation}
\end{cor}

%\subsection{$(2,1)$-(pre)sheaf of transitive Courant algebroids with connections}\label{SheafCA}
%\begin{pro}\label{$(2,1)$-presheafCA}
There is a $(2,1)$-presheaf of transitive Courant algebroids with
connections $\tcalgdp: \Mfd^{\op} \to \Gpd$, where $\Mfd^{\op}$ is the
opposite  category of
$\Mfd$, and $\Gpd$ is the
$2$-category of (discrete) groupoids and groupoid morphisms.

For an object $U\in \Mfd$, the groupoid $\tcalgdp(U)$ is made up
by the following data:
\begin{itemize}
\item{$\tcalgdp(U)_0$}: an object is a $6$-tuple %I checked online, 6 elements should be called sextuple.
$(\tagu, \Courant{-,-}^T_S, \pair{-,-}^T, \pr_{TU}, \theta, B)$, where $\theta\in\Omega^1(U,\g)$ is a $\g$-valued $1$-form, $B\in \Omega^2(U)$ is a $2$-form, and $\Courant{-,-}^T_S$ and
  $\pair{-,-}^T$ are given by \eqref{eq:TSb} and \eqref{eq:Tpair} respectively. We will simply denote an object by $(\tagu,\theta,B)$ in the sequel.

\item{$\tcalgdp(U)_1$}:  a $1$-morphism from  $(\tagu, \widetilde{\theta},  \widetilde{B} )$ to   $(\tagu,  \theta, B )$   is an automorphism of the standard transitive Courant algebroid $(\tagu,\Courant{-,-}^T_S,  \pair{-,-}^T,\pr_{TU})$ given by the
  matrix
$\left(\begin{array}{ccc}
1&0&0\\
\phi&\tau&0\\
\beta&-2\phi^\star\tau&1
\end{array}\right)$
such that \begin{eqnarray}
\label{eq:congluetc1}\theta(X)&=&\tau\widetilde{\theta}(X)+\phi(X),\\
\label{eq:congluetc2} \ii_X(\widetilde{B}-B)&=&\beta(X)+(\theta,\theta(X))^\g-(\widetilde{\theta},\widetilde{\theta}(X))^\g-2\phi^\star\tau(\widetilde{\theta}(X)),\\
\label{eq:congluetc3}\cs_3(\widetilde{\theta})-\cs_3(\theta)&=&d(\widetilde{B}-B).
 \end{eqnarray}
 Here the Chern-Simon 3-form $\cs_3(\theta)$ of $\theta$ is a $3$-form on $U$ defined by
 \eqref{eq:cs3-form} using bilinear form $(-, -)^\g$.

 The
composition of $1$-morphisms is simply the matrix multiplication.
\end{itemize}
Then for a morphism $\varphi: U\to V$ in $\Mfd$,  the associated
functor $\tcalgdp(\varphi): \tcalgdp(V) \to \tcalgdp(U)$ is induced by
pulling back forms. \calc{ More precisely, on the level of objects,
\begin{eqnarray*}
\tcalgdp(\varphi)(\tagv, \Courant{-,-}^T_S,  \pair{-,-}^T,\pr_{TV}, \theta,B)
= (\tagu, \Courant{-,-}^T_S,  \pair{-,-}^T,\pr_{TU},\varphi^*\theta, \varphi^*B),
\end{eqnarray*}
and on the level of morphisms,
\[
\tcalgdp(\varphi)\left(\begin{array}{ccc}
1&0&0\\
\phi&\tau&0\\
\beta&-2\phi^\star\tau&1
\end{array}\right)= \left(\begin{array}{ccc}
1&0&0\\
\varphi^*\phi&\varphi^*\tau&0\\
\varphi^*\beta&-2(\varphi^*\phi)^\star(\varphi^*\tau)&1
\end{array}\right).
\]
One verifies that $\tcalgdp(\varphi)$ is indeed a functor between
desired categories.}
%  \end{pro}

\vspace{3mm}

Take an open cover $\{U_i\}$ of $M\in
\Mfd$. An object in $\holim \tcalgdp(U(M)_\bullet)$ consists of
\begin{itemize}
\item  an object $\sqcup (\tagui,
  \theta_i,B_i)$ in $\tcalgdp(\sqcup U_i)_0$,

\item a $1$-morphism $\Lambda_{ij}=\left(\begin{array}{ccc}
1&0&0\\
\phi_{ij}&\tau_{ij}&0\\
\beta_{ij}&-2\phi_{ij}^\star\tau_{ij}&1
\end{array}\right)$ in $\tcalgdp(\sqcup U_{ij})_1$ from $(\taguij,\theta_j|_{U_{ij}}, B_j|_{U_{ij}})$ to $(\taguij,    \theta_i|_{U_{ij}}, B_i|_{U_{ij}}).
$
This implies that
\begin{eqnarray}
\label{eq:morcon11} \theta_i(X)&=&\tau_{ij}\theta_j(X)+\phi_{ij}(X),\\
\label{eq:morcon22} \ii_X(B_j-B_i)&=&\beta_{ij}(X)+(\theta_i,\theta_i(X))^\g-(\theta_j,\theta_j(X))^\g-2\phi_{ij}^\star\tau_{ij}(\theta_j(X)),\\
\label{eq:morcon33}\cs_3(\theta_j)-\cs_3(\theta_i)&=&d(B_j-B_i).
\end{eqnarray}
\item compatibility conditions $\Lambda_{ij}\Lambda_{jk}=\Lambda_{ik}$ on $U_{ijk}$, which are equivalent to the following equations
     \begin{eqnarray*}
      \phi_{ij}+\tau_{ij}\phi_{jk}&=&\phi_{ik},\\
      \tau_{ij}\tau_{jk}&=&\tau_{ik},\\
      \beta_{ij}-2\phi^\star \tau_{ij}\phi_{jk}+\beta_{jk}&=&\beta_{ik},\\
      -2\phi_{ij}^\star\tau_{ij}\tau_{jk}-2\phi_{jk}^\star\tau_{jk}&=&-2\phi_{ik}^\star\tau_{ik}.
    \end{eqnarray*}
\end{itemize}

\begin{defi}
An object $C_c=\Big(\sqcup(\tagui, \Courant{-,-}^T_S,  \pair{-,-}^T,\pr_{TU_i}\theta_i,B_i);\phi_{ij},\tau_{ij},\beta_{ij} \Big)$ in $\holim \tcalgdp(U(M)_\bullet)$ is called a {\bf transitive Courant data}. For simplicity, we also denote $C_c$ by  $(\theta_i, B_i;\phi_{ij},\tau_{ij},\beta_{ij})$ when there is no confusion.
\end{defi}

Note that $\holim \tcalgdp(U(M)_\bullet)$ of the $(2,1)$-presheaf $\tcalgdp$ might be empty. When it is not, we may describe the objects and morphisms by the following two propositions. We then describe the condition for  $\holim \tcalgdp(U(M)_\bullet)$ to be non-empty.

%In the following, we will always consider objects over fine enough good cover $\{U_i\}$. One can see that this is necessary to solve the decent equations used in this paper.
%Note that the local standard Courant algebroid $(\tagui, \Courant{-,-}^T_S,  \pair{-,-}^T,\pr_{TU_i})$ is canonical. Thus we will simply denote a Courant data by $(\theta_i, B_i;\phi_{ij},\tau_{ij},\beta_{ij})$.

\begin{pro}\label{thm:tco}\label{pro:tco}
A transitive Courant data  $C_c=(\theta_i, B_i;\phi_{ij},\tau_{ij},\beta_{ij})$ gives rise to a  transitive Courant algebroid $(C,\Courant{-,-},\pair{-,-},\rho)$ with a connection.% $s:TM\longrightarrow C$ and $\sigma_s:\huaG\longrightarrow\ker(\rho)$.
\end{pro}
%\cc{I feel that we need good covers here, no? because if your cover is not fine enough, how can you glue to a non-trivial principal bundle, in the principal bundle case. It should be also like this on the morphism level. And is it like this for all such propositions in Courant section? }
\pf
Since $\Lambda_{ij}\Lambda_{jk}=\Lambda_{ik}$, $\tagof{U_i}$'s glue to  a vector bundle $C$. Since $\Lambda_{ij}$ preserves the standard bracket $\Courant{-,-}^T_S$ and the standard pairing $\pair{-,-}^T$ on $\tagof{U_i}$, we have a well-defined bracket $\Courant{-,-}$ and a  nondegenerate symmetric bilinear form $\pair{-,-}$ on $\Gamma(C)$. Clearly, we have the following exact sequence of vector bundles:
$$
0\stackrel{}{\longrightarrow}\ker\rho\stackrel{}{\longrightarrow}C\stackrel{\rho}{\longrightarrow}TM\stackrel{}{\longrightarrow}0,
$$
where $\rho$ is induced by the projection $\tagof{U_i}\longrightarrow TU_i$.
\calc{ $$
C=\coprod TU_i\oplus(U_i\times\g)\oplus T^*U_i/\thicksim,
$$
where the equivalence relation $\thicksim$ is given by
$$
X+a+\xi\thicksim Y+b+\eta \Longleftrightarrow \left(\begin{array}{c}Y\\b\\\eta\end{array}\right)=\Lambda_{ij}\left(\begin{array}{c}X\\a\\\xi\end{array}\right),
$$
for all $X+a+\xi\in \tagof{U_j},~Y+b+\eta\in \tagof{U_i}$.}

The fact that $\Lambda_{ij}$ preserves the standard bracket $\Courant{-,-}^T_S$ and the standard pairing $\pair{-,-}^T$ also implies that Axioms (i)-(iii) in Definition \ref{defi:ca} are satisfied. Therefore, $(C,\Courant{-,-},\pair{-,-},\rho)$ is a transitive Courant algebroid.

On $U_i$, consider the splitting $s_i:TU_i\longrightarrow C|_{U_i}$ given by
\begin{equation}\label{eq:tcs1}
s_i(X)=X+\theta_i(X)-(\theta_i,\theta_i(X))^\g-\ii_XB_i.
\end{equation}
Straightforward calculation shows that $\pair{s_i(X),s_i(Y)}^T=0$. Thus, the splitting $s_i$ is isotropic.
\calc{
\begin{eqnarray*}
  \pair{s_i(X),s_i(Y)}^T&=& \pair{X+\theta_i(X)-(\theta_i,\theta_i(X))^\g-\ii_XB_i,Y+\theta_i(Y)-(\theta_i,\theta_i(Y))^\g-\ii_YB_i}^T\\
  &=&\half\Big(-(\theta_i(Y),\theta_i(X))^\g-B_i(X,Y)-(\theta_i(X),\theta_i(Y))^\g-B_i(Y,X)\Big)+(\theta_i(X),\theta_i(Y))^\g\\
  &=&0.
\end{eqnarray*}}
Eqs. \eqref{eq:morcon11} and \eqref{eq:morcon22} implies that $\Lambda_{ij}s_j(X)=s_i(X)$. Thus, we have a globally well-defined isotropic splitting
$s:TM\longrightarrow C$.
\calc{\begin{eqnarray*}
 \Lambda_{ij}s_j(X)&=&\left(\begin{array}{ccc}
1&0&0\\
\phi_{ij}&\tau_{ij}&0\\
\beta_{ij}&-2\phi_{ij}^\star\tau_{ij}&1
\end{array}\right)\left(\begin{array}{c}
X\\
\theta_j(X)\\
-(\theta_j,\theta_j(X))^\g-\ii_XB_j
\end{array}\right)\\
&=&\left(\begin{array}{c}
X\\
\phi_{ij}(X)+\tau_{ij}\theta_j(X)\\
\beta_{ij}(X)-2\phi_{ij}^\star\tau_{ij}\theta_j(X))^\g-(\theta_j,\theta_j(X))^\g-\ii_XB_j
\end{array}\right)\\
&=&\left(\begin{array}{c}
X\\
\theta_i(X)\\
-(\theta_i,\theta_i(X))^\g-\ii_XB_i
\end{array}\right)\\
&=&s_i(X).
\end{eqnarray*}}

Furthermore, $\coprod U_i\times \g$ and the transition function $\tau_{ij}$ give us a Lie algebra bundle $\huaG$, and there is a short exact sequence
$$
0\stackrel{}{\longrightarrow}T^*M\stackrel{}{\longrightarrow}\ker\rho\stackrel{\rho'}{\longrightarrow}\huaG\stackrel{}{\longrightarrow}0,
$$
where $\rho'$ is induced by the projection $ (U_i\times\g)\oplus T^*U_i\longrightarrow U_i\times \g$.
Consider the splitting $\sigma_{s_i}:\huaG|_{U_i}\longrightarrow (\ker\rho)|_{U_i}$ given by
\begin{equation}\label{eq:tcs2}
\sigma_{s_i}(a)=a-2(\theta_i,a)^\g.
\end{equation}
Then $\pair{s_i(X),\sigma_{s_i}(a)}^T=0$.
\calc{
It is straightforward to obtain that
\begin{eqnarray*}
 \pair{s_i(X),\sigma_{s_i}(a)}^T&=& \pair{X+\theta_i(X)-(\theta_i,\theta_i(X))^\g-\ii_XB_i,a-2(\theta_i,a)^\g}^T\\
  &=&\half\Big(-2(\theta_i(X),a)^\g \Big)+(\theta_i(X),a)^\g\\
  &=&0.
\end{eqnarray*}}
Thus, $\sigma_{s_i}$ is orthogonal to $s$. By \eqref{eq:morcon11}, we have
$$
\left(\begin{array}{cc}
 \tau_{ij}&0\\
 -2\phi_{ij}^\star\tau_{ij} &1
\end{array}\right)\left(\begin{array}{c}
 a\\
 -2(\theta_j,a)^\g
\end{array}\right)=\left(\begin{array}{c}
 \tau_{ij}a\\
 -2\phi_{ij}^\star\tau_{ij}(a)-2(\theta_j,a)^\g
\end{array}\right)=\left(\begin{array}{c}
 \tau_{ij}(a)\\
 -2(\theta_i,\tau_{ij}(a))^\g
\end{array}\right),
$$
which implies that we have a globally well-defined splitting $\sigma_{s}:\huaG\longrightarrow \ker\rho$ that orthogonal to the splitting $s$.\qed

\begin{rmk}
In Appendix \ref{GluingCA}, we write down the explicit formula for the glued transitive Courant algebroid in the form provided in \cite{ChenRCA}.
\end{rmk}

A 1-morphism in $\holim \tcalgdp(U(M)_\bullet)$ from  $(\widetilde{\theta}_i, \widetilde{B}_i; \widetilde{\phi}_{ij}, \widetilde{\tau}_{ij}, \widetilde{\beta}_{ij}) $  to  $(\theta_i, B_i;\phi_{ij},\tau_{ij},\beta_{ij}) $ consists of
  a 1-morphism $\left(\begin{array}{ccc}
1&0&0\\
\phi_i&\tau_i&0\\
\beta_i&-2\phi_i^\star\tau_i&1
\end{array}\right)$  in $\tcalgdp(\cup U_i)_1$ from $\sqcup (\tagui,
  \widetilde{\theta}_i, \widetilde{B}_i) $ to $\sqcup (\tagui,
  \theta_i, B_i)$, which satisfies
  \begin{equation}\label{eq:gluetcam}
    \Lambda_{ij}\left(\begin{array}{ccc}
1&0&0\\
\phi_j&\tau&0\\
\beta_j&-2\phi_j^\star\tau_j&1
\end{array}\right)=\left(\begin{array}{ccc}
1&0&0\\
\phi_i&\tau&0\\
\beta_i&-2\phi_i^\star\tau_i&1
\end{array}\right)\widetilde{\Lambda}_{ij}.
  \end{equation}

\begin{pro}\label{pro:tcm}
A $1$-morphism  in $\holim \tcalgdp(U(M)_\bullet)$ gives rise to a Courant algebroid isomorphism preserving connections.
  \end{pro}
\pf
The proof is similar to that of Proposition \ref{pro:tlm}. Eq. \eqref{eq:gluetcam} implies that the local morphisms glue together to a global morphism $\frkB$ between the gluing results of two Courant data. Eqs. \eqref{eq:congluetc1} and \eqref{eq:congluetc2} imply that $ \frkB \widetilde{s}_i(X) = s_i(X)$ and $\frkB\sigma_{\widetilde{s}}=\sigma_s$. Thus $\frkB$ preserves connections.
\calc{By \eqref{eq:gluetcam},   $\{ \left(\begin{array}{ccc}
1&0&0\\
\phi_i&\tau&0\\
\beta_i&-2\phi_i^\star\tau_i&1
\end{array}\right)\}$ glue to a bundle map $\frkB:C\longrightarrow C$. Obviously, $\frkB$ preserves the pairings and anchor maps. Furthermore, also by \eqref{eq:gluetcam}, $\frkB$ also preserves the brackets. Thus, $\frkB$ is a homomorphism between Courant algebroids.

 By \eqref{eq:congluetc1} and \eqref{eq:congluetc2},  we have
\begin{eqnarray*}
  \frkB \widetilde{s}_i(X)&=&\left(\begin{array}{ccc}
1&0&0\\
\phi_{i}&\tau_{i}&0\\
\beta_{i}&-2\phi_{i}^\star\tau_{i}&1
\end{array}\right)\left(\begin{array}{c}
X\\
\widetilde{\theta}_i(X)\\
-(\widetilde{\theta}_i,\theta_i(X))^\g-\ii_X\widetilde{B}_i
\end{array}\right)\\
&=&\left(\begin{array}{c}
X\\
\phi_{i}(X)+\tau_{i}\widetilde{\theta}_i(X)\\
\beta_{i}(X)-2\phi_{i}^\star\tau_{i}\widetilde{\theta}_i(X))^\g-(\widetilde{\theta}_i,\widetilde{\theta}_i(X))^\g-\ii_X\widetilde{B}_i
\end{array}\right)\\
&=&\left(\begin{array}{c}
X\\
\theta_i(X)\\
-(\theta_i,\theta_i(X))^\g-\ii_XB_i
\end{array}\right)\\
&=&s_i(X).
\end{eqnarray*}
Similarly, one can show that   $\frkB\sigma_{\widetilde{s}}=\sigma_s$.
}\qed\\

%22.11 Tuesday afternoon

By Proposition \ref{pro:tco} and Proposition \ref{pro:tcm},  after the
plus construction, we arrive at the $(2,1)$-sheaf ${\tcalgdp}^+$ of
transitive Courant algebroids with connections, where the 1-morphism
are the isomorphisms of Courant algebroid preserving connections. See
\cite{LBM} for the notion of morphisms (not necessary isomorphisms) of
Courant algebroids. However, in our case, all our sheaves (stacks) are
functors to (higher) groupoids, that is, 1-morphisms are always isomorphisms.

Obviously, there is a projection $\pr$ from the $(2,1)$-presheaf $\tcalgdp$ to the $(2,1)$-presheaf $\tlp$, which sends a transitive Courant data $\Big(\sqcup(\tagui, \Courant{-,-}^T_S,  \pair{-,-}^T,\pr_{TU_i},\theta_i,B_i);\phi_{ij},\tau_{ij},\beta_{ij} \Big)$ to the transitive Lie data $\Big(\sqcup (TU_i\oplus(U_i\times\g),
  [-,-]^T_S,\pr_{TU},\theta_i);\phi_{ij},\tau_{ij} \Big)$ (see Definition \ref{Liedata}), and behaves in a similar obvious way on the level of morphisms. After plus construction, we arrive at a projection $\tcalgdp^+\xrightarrow{\pr}\tlp^+$. Here $\tlp$ and $ \tlp^+$ are the $(2,1)$-presheaf and $(2,1)$-sheaf of transitive Lie algebroids respectively. See Appendix \ref{app:tl} for more details.

%%We know that $\holim \tlp(U(M)_\bullet)$ is not empty and its object $\Big(\sqcup (TU_i\oplus(U_i\times\g),
%  [-,-]^T_S,\pr_{TU_i},\theta_i),\phi_{ij},\tau_{ij}\Big)$ glue to a transitive Lie algebroid. \cc{why it's non-empty here? is it because over $M$ there is always a transitive Lie algebroid $A$ for a given $G$, and we may always decompose it to local data?} Thus the existence of Courant data $\Big(\sqcup (\tagui,
%  \Courant{-,-}^T_S, \pair{-,-}^T,\pr_{TU_i},\theta_i, B_i),\phi_{ij},\tau_{ij},\beta_{ij}\Big)$ amounts to the question of the following lifting problem. It turns out that the lift has an obstruction which is the first Pontryagin class of the underlining transitive Lie data. \cc{we can't say this obstruction is topological, can we? there isn't diff. top. here, or?}

Now we fix a transitive Lie data $A_c=\Big(\sqcup (TU_i\oplus(U_i\times\g),
[-,-]^T_S,\pr_{TU_i},\theta_i);\phi_{ij},\tau_{ij}\Big)$. We simply denote it by $(\theta_i;\phi_{ij},\tau_{ij})$.
\begin{lem}\label{lem:gglueR}
Define  $R_i:\wedge^2\Gamma(TU_i)\longrightarrow \Gamma(U_i\times \g)$ by
\begin{equation}
  R_i=d\theta_i+\half[\theta_i,\theta_i]_\g.
\end{equation}  Then $R_i$'s glue to a globally well-defined curvature $R:\wedge^2TM\longrightarrow\huaG$.   We call $R$ the {\bf curvature} of our transitive Lie data.
\end{lem}
\pf We need to show $\tau_{ij}R_j=R_i$. We bring \eqref{eq:morcon11} inside the expression. Then the result follows from Eqs. \eqref{eq:auto2} and \eqref{eq:auto3}.\qed
\calc{
By \eqref{eq:auto2}, \eqref{eq:auto3} and \eqref{eq:morcon11}, we have
\begin{eqnarray*}
  R_i(X,Y)-\tau_{ij}R_j(X,Y)&=&[\theta_i(X),\theta_i(Y)]_\g+d\theta_i(X,Y)-\tau_{ij}\big([\theta_j(X),\theta_j(Y)]_\g+d\theta_j(X,Y)\big)\\
  &=&[\tau_{ij}\theta_j(X)+\phi_{ij}(X),\tau_{ij}\theta_j(Y)+\phi_{ij}(Y)]_\g\\
  &&+X\big(\tau_{ij}\theta_j(Y)+\phi_{ij}(Y)\big)-Y\big(\tau_{ij}\theta_j(X)+\phi_{ij}(X)\big)\\
  &&-\tau_{ij}\theta_j([X,Y])-\phi_{ij}([X,Y])-\tau_{ij}[\theta_j(X),\theta_j(Y)]_\g\\
    &&-\tau_{ij}X(\theta_j(Y))+\tau_{ij}Y(\theta_j(X))+\tau_{ij}(\theta_j([X,Y]))\\
    &=&X(\phi_{ij}(Y))-Y(\phi_{ij}(X))-\phi_{ij}([X,Y])+[ \phi_{ij}(X), \phi_{ij}(Y)]_\g\\
    &&+X (\tau_{ij}\theta_j(Y))-\tau_{ij}X(\theta_j(Y))+[\phi_{ij}(X),\tau_{ij}\theta_j(Y)]_\g\\
    &&-Y (\tau_{ij}\theta_j(X))+\tau_{ij}Y(\theta_j(X))-[\phi_{ij}(Y),\tau_{ij}\theta_j(X)]_\g\\
    &=&0.
\end{eqnarray*}}
\begin{rmk}
 It turns out that the curvature form $R$ in this lemma is the same $R$ appearing in the Courant bracket \eqref{eq:Tbracket} for the gluing result (see Proposition \ref{GlobalCA}).
\end{rmk}

Then the {\bf first Pontryagin class} of the transitive Lie data $A_c$ is defined through the curvature $R$,
\begin{equation}
  p_1(A_c):=(R,R)^\g.
\end{equation}

\begin{thm} \label{thm:lift-courant}
  Given an object $M\xrightarrow{A_c} \tlp^+$ in $\tlp^+$ over a good cover $\{U_i\}$ of $M$,
  \begin{enumerate}
      \item[\rm(i)] there exists a lift $C_c$ that fits the diagram
  \begin{equation} \label{eq:lifting-courant}
\xymatrix{
&\tcalgdp^+\ar[d]\\
M \ar@{.>}[ur]^{C_c}\ar[r]^{A_c}&\tlp^+,
}
\end{equation}
if and only if $p_1(A_c)=0$;
\item[\rm(ii)] if a lift exists, then the space of lifts up to isomorphisms forms a torsor of $\Omega^3_{\cl}(M)$.
  \end{enumerate}
\end{thm}
The proof of this theorem is given in the following two lemmas.
\begin{lem}\label{lem:Pexact}
There exists a lift $C_c$ of $A_c$ in diagram \eqref{eq:lifting-courant} if and only if
  $$p_1(A_c)=0.$$
\end{lem}
\pf
First, assume that $C_c=\Big(\sqcup(\tagui, \Courant{-,-}^T_S,  \pair{-,-}^T,\pr_{TU_i},\theta_i,B_i);\phi_{ij},\tau_{ij},\beta_{ij} \Big) $ is a lift of $A_c$, then set $H_i:=\cs_3(\theta_i)-dB_i$ on each $U_i$. By \eqref{eq:morcon33}, we have $\cs_3(\theta_i)-dB_i=\cs_3(\theta_j)-dB_j$, which implies that $H_i=H_j$ on $U_{ij}$. Thus, $H_i$'s glue to a global $3$-form $H$. We call $H$ the {\bf curvature} of $C_c$.
Since $H|_{U_i}=\cs_3(\theta_i)-dB_i$, and $\cs_3(\theta_i)=(\theta_i, d\theta_i)^\g+\frac{1}{3} (\theta_i, [\theta_i, \theta_i])^\g$, we have
\begin{eqnarray*}
  dH|_{U_i} &=&d((\theta_i,d\theta_i)^\g+\frac{1}{3}(\theta_i,[\theta_i,\theta_i])^\g)
 % &=&(d\theta_i,d\theta_i)^\g+\frac{1}{3}\Big((d\theta_i,[\theta_i,\theta_i])^\g-(\theta_i,d[\theta_i,\theta_i])^\g\Big)\\
 % &=&(d\theta_i,d\theta_i)^\g+\frac{1}{3}\Big((d\theta_i,[\theta_i,\theta_i])^\g-2(\theta_i,[d\theta_i,\theta_i])^\g\Big)\\
  =(d\theta_i,d\theta_i)^\g+ (d\theta_i,[\theta_i,\theta_i])^\g.
\end{eqnarray*}
Here we use the fact that $(-,-)^\g$ is adjoint invariant. Another direct computation gives
\begin{eqnarray*}
  (R_i,R_i)^\g&=&(d\theta_i+\half[\theta_i,\theta_i]_\g,d\theta_i+\half[\theta_i,\theta_i]_\g)^\g\\
  &=&(d\theta_i,d\theta_i)^\g+ (d\theta^i,[\theta_i,\theta_i])^\g.
\end{eqnarray*}
%+\frac{1}{4}( [\theta_i,\theta_i]_\g, [\theta_i,\theta_i]_\g)^\g.
%Furthermore, using the Jacobi identity and the fact that $(-,-)^\g$ is adjoint invariant, we can obtain  $$( [\theta_i,\theta_i]_\g, [\theta_i,\theta_i]_\g)^\g=0.$$
Therefore, we have $(R,R)^\g=dH,$ i.e. $p_1(A_c)=0.$

On the other hand, given $A_c=\Big(\sqcup (TU_i\oplus(U_i\times\g),
[-,-]^T_S,\pr_{TU_i},\theta_i);\phi_{ij},\tau_{ij}\Big) $ and assuming $p_1(A_c)=0$, let $H\in \Omega^3(M)$ be such that $(R,R)^\g=dH$. Since the cover is good, we may come up with 2-forms $B_i\in \Omega^2(U_i)$ satisfying $$H|_{U_i}=\cs_3(\theta_i)-dB_i.$$
Let $\beta_{ij}:TU_{ij}\longrightarrow T^*U_{ij}$ be the bundle map uniquely determined by \eqref{eq:morcon22}. Then one checks that $\Big(\sqcup(\tagui, \Courant{-,-}^T_S,  \pair{-,-}^T,\pr_{TU_i},\theta_i,B_i);\phi_{ij},\tau_{ij},\beta_{ij} \Big) $ is a Courant data whose image under the projection $\tcalgdpp \xrightarrow{\pr} \tlpp$ coincides with $A_c$.
\qed\\

Now assume that there exists a lift of $A_c=\Big(\sqcup (TU_i\oplus(U_i\times\g),
[-,-]^T_S,\pr_{TU_i},\theta_i);\phi_{ij},\tau_{ij}\Big)\in\tlp^+(M)$, where again $A_c$ is over a good cover of $M$. We denote the fiber category of $\pr$ over $A_c$ by $S_{A_c}$.%This is a category with the morphisms induced from the 1-morphism of $\tcalgdp^+(M)$.
Then the space of lifts up to isomorphisms is the set of equivalent classes $S_{A_c}/_{1\mbox{-morphisms}}$.
We define an action of $\Omega^3_{\rm cl}(M)$ on this space as follows.  For all $h\in\Omega^3_{\rm cl}(M)$, assume that $h|_{U_i}=dB^h_i$. Define the action of $h$ on  $S_{A_c}/_{1\mbox{-morphisms}}$  by
\begin{equation}\label{eq:h-action}
h\cdot [ (\theta_i, B_i;\phi_{ij},\tau_{ij},\beta_{ij})] =[ ( \theta_i,
B_i+B^h_i;\phi_{ij},\tau_{ij},\beta_{ij}+B^h_j-B^h_i )].
\end{equation}
Here $[ \ \ ]$ denotes the isomorphism class in the quotient $S_{A_c}/_{1\mbox{-morphisms}}$.
Now we prove that the above action is well defined.
\begin{itemize}
\item It does not depend on the choice of $\{B^h_i\}$. In fact, if $h|_{U_i}=d\bar{B}^h_i$, we let  $\huaB_i=\bar{B}^h_i-B^h_i$, which is closed. Then $\left(\begin{array}{ccc}
1&0&0\\
0&1&0\\
\huaB_i&0&1
\end{array}\right)$ gives rise to an isomorphism from $(\theta_i,
B_i+\bar{B}^h_i ;
\phi_{ij},\tau_{ij},\beta_{ij}+\bar{B}^h_j-\bar{B}^h_i )$ to
$(\theta_i, B_i+B^h_i;\phi_{ij},\tau_{ij},\beta_{ij}+B^h_j-B^h_i )$.

\item It does not depend on the choice of a representative $(\theta_i, B_i;\phi_{ij},\tau_{ij},\beta_{ij})$ of an isomorphism class.  The isomorphism is possibly given through another finer cover $\{V_i\}$. However as we may always pull back our data to $V_i$, we might as well assume that $V_i=U_i$.  Assume that we have a 1-morphism $\Lambda$ from a transitive Courant data $(\theta_i, \widetilde{B}_i;\phi_{ij},\tau_{ij},\widetilde{\beta}_{ij})\in S_{A_c}$ to a transitive Courant data $(\theta_i, B_i;\phi_{ij},\tau_{ij},\beta_{ij})\in S_{A_c}$.
  Then by \eqref{eq:congluetc3}, one checks that $\Lambda$ is a 1-morphism from $(\theta_i, \widetilde{B}_i+B^h_i;\phi_{ij},\tau_{ij},\widetilde{\beta}_{ij}+B^h_j-B^h_i)$ to $(\theta_i, B_i+B^h_i;\phi_{ij},\tau_{ij},\beta_{ij}+B^h_j-B^h_i)$.
\end{itemize}

  \begin{lem}\label{thm:catorsor}
   With the above notations, $\Omega^3_{\rm cl}(M)$ acts on  $S_{A_c}/_{1\mbox{-morphisms}}$ freely and transitively. Thus, $S_{A_c}/_{1\mbox{-morphisms}}$ is an $\Omega^3_{\rm cl}(M)$-torsor.
  \end{lem}
\pf We first show that the action of $\Omega^3_{\cl}(M)$ is free.
If $h\cdot [(\theta_i, B_i;\phi_{ij},\tau_{ij},\beta_{ij})]$ is isomorphic to $[(\theta_i, B_i;\phi_{ij},\tau_{ij},\beta_{ij})]$, by \eqref{eq:congluetc3}, we deduce that $dB^h_i=0$, which implies that $h=0.$

Then we show that the action is transitive. For  two objects $ (\theta_i, B_i;\phi_{ij},\tau_{ij},\beta_{ij})$ and $ (\theta_i, B_i';\phi_{ij},\tau_{ij},\beta_{ij}')$ in $S_{A_c}$, let $B^h_i=B_i'-B_i$. Since both $\{\cs_3(\theta_i)-dB_i\}$ and $\{\cs_3(\theta_i)-dB_i'\}$ can be glued to a global 3-form, we deduce that $\{dB^h_i\}$ can be glued to a global closed 3-form $h$. By \eqref{eq:morcon22}, we deduce that
\begin{eqnarray*}
\ii_X(B_j'-B_i')&=&\beta_{ij}'+(\theta_i,\theta_i(X))^\g-(\theta_j,\theta_j(X))^\g-2\phi_{ij}^\star\tau_{ij}(\theta_j(X)),\\
\ii_X(B_j-B_i)&=&\beta_{ij}+(\theta_i,\theta_i(X))^\g-(\theta_j,\theta_j(X))^\g-2\phi_{ij}^\star\tau_{ij}(\theta_j(X)),
\end{eqnarray*}
which implies that $$\beta_{ij}'=\beta_{ij}+B^h_j-B^h_i.$$
Thus, $h\cdot (\theta_i, B_i;\phi_{ij},\tau_{ij},\beta_{ij})= (\theta_i, B_i';\phi_{ij},\tau_{ij},\beta_{ij}')$. It finishes the proof.\qed

%Given a transitive Courant data $C_c=(\theta_i, B_i;\phi_{ij},\tau_{ij},\beta_{ij})$, by \eqref{eq:morcon33}, $\{\cs(\theta_i)-dB_i\}$ give rise to %a global 3-form $H$ on $M$, which we define to be the {\bf curvature} of the $C_c$.

\begin{cor}\label{cor:curv-courant}
  Given an object $M\stackrel{A_c}{\longrightarrow}\tlp^+$ with $p_1(A_c)=0$, the curvatures of all possible lifts form an $\Omega^3_{\rm cl}(M)$-torsor.
\end{cor}

\begin{rmk}
 Note that objects in the category  ${\tcalgdp}^+(M)$ are transitive Courant algebroids with connections and $1$-morphisms in  ${\tcalgdp}^+(M)$ are isomorphisms between transitive Courant algebroids that preserve connections.
 Thus, obviously there is a forgetful functor $\huaF$ from the category  ${\tcalgdp}^+(M)$ to the usual category of transitive Courant algebroids $\tc(M)$. Similarly, there is also a forgetful functor from the category $\tlp^+(M)$ to the usual category of transitive Lie algebroids $\tl(M)$ over $M$.  %For an object  $A_c=\Big(\sqcup (TU_i\oplus(U_i\times\g),
%  [-,-]^T_S,\pr_{TU_i},\theta_i),\phi_{ij},\tau_{ij}\Big)$ in
%  $\tlp^+(M)$, forget the connection, we obtain a transitive Lie
%  algebroid $A\in \tl(M)$.
Let $S_A$ denote the fiber of the projection $\pr:\tc(M)\longrightarrow \tl(M)$. Then we have the following diagram:
  \[
 \xymatrix{S_{A_c}\ar[r]\ar[d]^{\huaF}& {\tcalgdp}^+(M)   \ar[r]^{\pr}\ar[d]^{\huaF}&\tlp^+(M)\ar[d]^{\huaF}\\
S_A\ar[r] &\tc(M)\ar[r]^{ \pr}&\tl(M).
 }
 \]
 In \cite{ChenRCA}, the isomorphism classes in $S_{A}$ are classified
 for a quadratic transitive Lie algebroid $A\in\tl(M)$ with vanishing
 first Pontryagin class. Compared to Lemma \ref{thm:catorsor}, the
 classification results are not the same since our morphisms need to
 preserve connections. More precisely, our space is a torsor of
 $\Omega_{\rm cl}^3(M)$ which surjectively maps to the group $H^3(M,
 \R)/I$ of which their space is a torsor. Here $I$ is a certain
 subspace of $H^3(M, \R)$ which comes from automorphisms of the Lie
 algebroid $A$. In our case, such automorphisms do not show up because
 our $\Omega^3_{\cl}(M)$-action \eqref{eq:h-action} fixes the underlying
 transitive Lie data.
\end{rmk}

\begin{ep}[action Courant algebroids]{\rm\label{ep:action-courant}
A quadratic Lie algebra $(\frkk, [-,-]_\frkk,(-, -)^\frkk)$ gives rise to a string Lie 2-algebra \cite[Sect. 2]{sheng-zhu} and may be viewed as a Courant algebroid over a point. Extending this idea to a general base manifold $M$, the authors in \cite{LBM} construct a natural example of Courant algebroids as follows: given an action $\rho: \frkk \to \frkX(M)$ whose stabilizer at each point on $M$ is a coisotropic subspace of $\frkk$,  there is an {\bf action Courant algebroid} $(M\times\frkk,\Courant{-,-},\pair{-, -},\rho)$, where the anchor is given by the action $\rho$, the bilinear form $\pair{-, -}$ is the pointwise pairing induced by $(-, -)^\frkk$ on  $\frkk$, and the Courant bracket on the sections of the vector bundle $M\times \frkk \to M$ is given by
\[ \Courant{X, Y}=[X, Y]_\frkk+L_{\rho(X)}Y-L_{\rho(Y)}X+ \rho^*\pair{dX, Y}, \quad \forall X, Y \in C^\infty(M, \frkk).\]

Now we study this example in the special case of homogeneous spaces.  Take $\frkk={\rm gl}_n(\mathbb{C})$ equipped with the nondegenerate bilinear form $(A, B)^\frkk=\tr(AB)$, and let $\frkk_{\ge 0}$, $\frkk_0$, and $\frkk_+$ be the Lie subalgebras of non-strict upper triangular, diagonal, and strict upper triangular matrices respectively.  Let $K={\rm GL}_n(\mathbb{C})$, and  $K_{\ge 0}$, $K_0$, and $K_+$ be the matrix groups corresponding to $\frkk_{\ge 0}$, $\frkk_0$, and $\frkk_+$ respectively. Take $M$ to be the homogeneous space $K/K_{\ge 0}$. In this case, the anchor $\rho$ is surjective and it follows that the action Courant algebroid $M\times \frkk$ is a transitive Courant algebroid. Thus, after choosing a connection, we have a split of Courant algebroid $M\times \frkk \cong TM\oplus \mathcal{G} \oplus T^*M$, with the Courant bracket and the pairing defined in \eqref{eq:Tbracket} and \eqref{eq:Tpair}.

We now prove that there exists a suitable connection such that the split form on the right hand side is a standard transitive Courant algebroid. Firstly, following \cite[Proposition 3.9]{Xu} the underlying transitive Lie algebroid $A=M\times \frkk / (\ker\rho)^\perp$ is the Atiyah Lie algebroid associated to the $K_0=(\mathbb{C^\times})^n$ principal bundle $K/K_+\to M$. However, this principal bundle is trivial. This can be seen as follows. Note that $M$ is isomorphic to the flag variety $F_n := \{E_\bullet=(E_1\subset E_2\subset\cdot\cdot\cdot\subset E_n = \mathbb{C}^n)~|~ {\rm dim} E_i= i\}$. Let $U_i$ be the tautological $i$-dimensional vector bundle over $M$, whose fiber over a point (a flag $E_\bullet$) is the vector space $E_i$ of the flag. These bundles form a filtration
$0 = U_0 \subset U_1\subset \cdot\cdot\cdot \subset U_n=F_n\times \mathbb{C}^n.$ If we consider the standard representation of $K_0\subset {\rm GL}_n(\mathbb{C})$ on $\mathbb{C}^n$, then the associated vector bundle of the $K_0$ principal bundle $K/K_+\to M$ is isomorphic to $\oplus_{i=1}^n L_i$. Here $L_i:=U_i/U_{i-1}$ is a line bundle for each $0<i\le n$. Therefore the triviality of $K/K_+\to M$ follows from the fact that the associated bundle $\displaystyle{\oplus_{i=1}^n} L_i\cong U_n$ is trivial.

Upon choosing a global trivialization, we can take the natural
connection on the principal bundle $K/K_+\to M$, which in turn induces
a split of the Atiyah Lie algebroid $A\cong TM\oplus (M\times\frkk_0)$
with connection $\nabla_X(a)=X(a)$. Let $TM\oplus  (M\times\frkk_0)
\oplus T^*M$ be the corresponding split of $M\times\frkk $. By
Proposition \ref{GlobalCA} and the fact that the curvature $R$ of the
natural connection is 0,  the Courant bracket on the split form is the
standard Courant bracket up to a $3$-form $H\in \Omega^3(M)$
satisfying $dH=(R, R)^{\frkk_0}=0$. Following from Borel's result, as
a complete flag variety, $M$ has vanishing odd cohomology. Thus we may
assume $H=dB$ for a certain 2-form $B$, and perform a $B$-field
transformation for $TM\oplus  (M\times\frkk_0) \oplus T^*M$,
$X+a+\xi\mapsto X +a+\xi+ \ii_X B$, and arrive at the standard
bracket. Thus composing these two steps, we obtain an isomorphism from
$ M\times\frkk$ to the standard transitive Courant algebroid.

Notice that the $K_0$ principal bundle $K/K_+\to M$ is not necessarily
trivial for  general $K$. For example, when $K={\rm
  SL}_n(\mathbb{C})$, the similar construction will give us nontrivial
$K_0=(\mathbb{C^\times})^{n-1}$ principal bundle\footnote{We thank
very much  Eckhard Meinrenken for sharing this example with
us.}. Since $K_0$ is abelian, the Cartan 3-form on it is
0, thus the basic gerbe on it is a trivial gerbe. Nevertheless, we
shall not expect string groups to be trivial. Therefore, there
will be different features for
abelian counterpart of string structure. We leave it for future discussion.
}
\end{ep}

\section{Morphism from the string sheaf to the transitive Courant sheaf}

Having constructed the sheaves of $\String(n)$-principal bundles and
transitive Courant algebroids with connections, we show that there is
a canonical morphism between them. On the level of objects without
connections, one could build the
correspondence between String structures and transitive Courant
algebroids using the reduction method as in
\cite{baraglia-hekmati,BCG}. Nevertheless, to obtain a functor, it is
convenient to use our language. Throughout this section, we take
the presheaf of transitive Courant algebroids with connections for
$G=\Spin(n)$ and $\g=\so(n)$ with the bilinear form $(-,-)^\g$ the one appeared in \eqref{eq:cs3-form}. We still denote this presheaf by $\tcalgdp$.

\subsection{Construction of the morphism $\Phi: \bstringnpp \to \tcalgdpp$}

%First, we take the $(3,1)$-presheaf of $ \String(n)$-principal bundles with connections $\bstringnp$ and the presheaf of transitive Courant algebroids with connections for $G=\Spin(n)$, which we still denote by $\tcalgdp$.
%Let $G$ be a Lie group with Lie algebra $\g$, $\TM=(dg)g^{-1}$ the right invariant Maurer-Cartan 1-form on $G$.
%Denote by $\overline{\TM}=g^{-1}dg=\Ad_{g^{-1}}\TM$.

\begin{thm}\label{Naturaltransformation}
There is a canonical morphism $\Phi$ from $\bstringnp$ to $\tcalgdp$, where
for any $U\in\Mfd$, the morphism $\Phi(U):\bstringnp(U)\rightarrow\tcalgdp(U)$ is given on $0$-, $1$- and $2$-simplices respectively as follows
%{\bf we think of 1-cat as a 2-cate with identities as the 2-morphisms? Because there is no 2-morphism, natural transformations are identities.}
\begin{itemize}
\item for an object $(U\times \String(n)\to U, \theta, B)$ in $\bstringnp(U)$, we have
$$\Phi(U\times \String(n)\to U, \theta, B)=(\tagu,
   \theta, B);$$

\item for a $1$-simplex $(g_{\scriptscriptstyle 01}, A_{\scriptscriptstyle 01}, \omega^2_{\scriptscriptstyle 01}): (U\times \String(n)\to U, \theta_1, B_1) \rightarrow (U\times \String(n)\to U, \theta_0, B_0),$ we have
$$\Phi(g_{\scriptscriptstyle 01}, A_{\scriptscriptstyle 01}, \omega^2_{\scriptscriptstyle 01})=\Lambda_{\scriptscriptstyle 01} := \left(\begin{array}{ccc}
1&0&0\\
-\bar{g}_{\scriptscriptstyle 01}^*\TM&\ad_{\bar{g}_{\scriptscriptstyle 01}}&0\\
\beta_{\scriptscriptstyle 01}&2(\bar{g}_{\scriptscriptstyle 01}^*{\TM})^\star\circ\ad_{\bar{g}_{\scriptscriptstyle 01}}&1
\end{array}\right),$$
where   $\bar{g}_{\scriptscriptstyle 01}:U\rightarrow G$ is the underlying morphism of $g_{\scriptscriptstyle 01}$ and $$\beta_{\scriptscriptstyle 01}=-(\bar{g}_{\scriptscriptstyle 01}^*\TM,\theta_0)^\g+dA_{\scriptscriptstyle 01}+\omega^2_{\scriptscriptstyle 01}-(\bar{g}_{\scriptscriptstyle 01}^*\TM)^\star \circ \bar{g}_{\scriptscriptstyle 01}^*\TM.$$

\item for a $2$-simplex $(f, \omega^1)$, we have $\Phi(f,\omega^1)=1$.

\end{itemize}
\end{thm}

\begin{rmk}
Note that the symmetric part of $\beta_{\scriptscriptstyle 01} $ is given by $-(\bar{g}_{\scriptscriptstyle 01}^*\TM)^\star \circ \bar{g}_{\scriptscriptstyle 01}^*\TM$, which is the same as the symmetric part of $\Psi$ given in \eqref{eq:severainner} in Appendix A.4. Thus, $\Lambda_{\scriptscriptstyle 01}$ is an inner automorphism in the sense of  ${\rm\check{S}}$evera \cite{sevlet}.
\end{rmk}

To prove this theorem, we need the following lemmas.
\begin{lem}For any $G$-valued function $g:U\longrightarrow G$, and $a,b\in\Gamma(U\times\g)$, $X,Y\in \Gamma(TU)$, we have
  \begin{eqnarray}
 \label{eq:Adauto} [\ad_{g}a,\ad_{g}b]_\g&=&\ad_{g}[a,b]_\g,\\
\label{eq:con1} X( g^*\TM(Y))-Y( g^*\TM(X))-[ g^*\TM(X), g^*\TM(Y)]_\g&=&g^*\TM([X,Y]),\\
 \label{eq:con2}X(\ad_{g}b )-
 [ g^*\TM(X),\ad_{g}b ]_\g&=&\ad_{g}(X(b)).
\end{eqnarray}
\end{lem}
\pf \eqref{eq:Adauto} is obvious. \eqref{eq:con1} follows from the Maurer-Cartan equation $d\TM-\half[\TM,\TM]_\g=0$.  For any $m\in U,$ let $\gamma(s) $ be the integration curve of $X$ through $m$, i.e. $X_m=\frac{d}{ds}|_{s=0}\gamma(s)$.  Then we have
\begin{eqnarray*}
  X_m(\ad_gb)&=&X_m(\frac{d}{dt}|_{t=0}\Ad_g\exp (tb) )\\
%  &=&\frac{d}{ds}|_{s=0}\frac{d}{dt}|_{t=0}\Ad_{g(\gamma(s))}\exp(tb(\gamma(s))) \\
  &=&\frac{d}{dt}|_{t=0}\Big(\frac{d}{ds}|_{s=0}\Ad_{g(\gamma(s))}\exp(tb(\gamma(s))) \Big)\\
 &=&\frac{d}{dt}|_{t=0}\Big( \Ad_{g(m)}\exp(tX_m(b)) +\frac{d}{ds}|_{s=0}\Ad_{g(\gamma(s))}\exp(tb(m))\Big)\\
 &=&\ad_{g(m)} X_m(b)+\frac{d}{ds}|_{s=0}\ad_{g(\gamma(s))\cdot g(m)^{-1}}\ad_{g(m)}b(m)\\
% &=&\ad_{g_m} X_m(b)+[\frac{d}{ds}|_{s=0}g(\gamma(s))(g(m))^{-1},\ad_{g(m)}b(m)]_\g\\
  &=&\ad_{g(m)} X_m(b)+[g(m)^*\TM(X_m),\ad_{g(m)}b(m)]_\g,
\end{eqnarray*}
which implies that \eqref{eq:con2} holds. The proof is finished. \qed

\begin{lem}\label{1morphism}  Let $(g_{\scriptscriptstyle 01}, A_{\scriptscriptstyle 01}, \omega^2_{\scriptscriptstyle 01})$ be a $1$-morphism from $(U\times \String(n)\to U, \theta_1, B_1)$ to $(U\times \String(n)\to U, \theta_0, B_0)$.
Then $\Lambda_{\scriptscriptstyle 01}$ given in Theorem \ref{Naturaltransformation} is a $1$-morphism in $\tcalgdp$ from $(\tagu,  \theta_1, B_1)$ to $ (\tagu, \theta_0, B_0)$.
\end{lem}
\pf
By definition, we first need to show that $\Lambda_{\scriptscriptstyle 01}$ is indeed an automorphism of
the standard transitive Courant algebroid $(\tagu, \Courant{-,-}^T_S,  \pair{-,-}^T,\pr_{TU})$. That is to prove the entries of the vector bundle map $\Lambda_{\scriptscriptstyle 01}=\left(\begin{array}{ccc}
1&0&0\\
-\bar{g}_{\scriptscriptstyle 01}^*\TM&\ad_{\bar{g}_{\scriptscriptstyle 01}}&0\\
\beta_{\scriptscriptstyle 01}&2(\bar{g}_{\scriptscriptstyle 01}^*\overline{\TM},\cdot)^\g&1
\end{array}\right)$  satisfy the identities \eqref{eq:auto1}-\eqref{eq:auto4}. Note that $(g_{\scriptscriptstyle 01}: U \to \String(n)$, $A_{\scriptscriptstyle 01}\in \Omega^1(U)$, $\omega^2_{\scriptscriptstyle 01}
  \in \Omega^2(U))$ gives rise to a 1-morphism,
\[(U\times \String(n)\to U, \theta_0, B_0) \xleftarrow {(g_{\scriptscriptstyle 01}, A_{\scriptscriptstyle 01},
  \omega^2_{\scriptscriptstyle 01})} (U\times\String(n)\to U, \theta_1,  B_1),\] if and only if
\begin{equation}\label{eq:diffcs}
dA_{\scriptscriptstyle 01} = \omega^2 +B_1-B_0, \quad d\omega^2_{\scriptscriptstyle 01} =
\cs_3(\theta_1)-\cs_3(\theta_0), \quad \theta_0-\ad_{\bar{g}_{\scriptscriptstyle 01}}
\theta_1 = - \bg^* \TM,
\end{equation}
where $\bar{g}_{\scriptscriptstyle 01}:U\rightarrow G$ is the underlying morphism of $g_{\scriptscriptstyle 01}:U\rightarrow \String(n)$.

The symmetric part of $\beta_{\scriptscriptstyle 01}:TU\longrightarrow T^*U$ is $-(\bar{g}_{\scriptscriptstyle 01}^*\TM)^\star \circ \bar{g}_{\scriptscriptstyle 01}^*\TM$, which we denote by $\beta_{\scriptscriptstyle 01}^{sym}$. Therefore,
\begin{eqnarray*}
&& \half\Big( \beta_{\scriptscriptstyle 01}(X)(Y)+\beta_{\scriptscriptstyle 01}(Y)(X)\Big)=\beta^{sym}_{\scriptscriptstyle 01}(X)(Y)=-(-\bar{g}_{\scriptscriptstyle 01}^*\TM(X),-\bar{g}_{\scriptscriptstyle 01}^*\TM(Y))^\g,
\end{eqnarray*}
which implies that \eqref{eq:auto1} holds.

 By \eqref{eq:con1} and \eqref{eq:con2}, we deduce that   \eqref{eq:auto2} and \eqref{eq:auto3} hold.

 The skewsymmetric part of $\beta_{\scriptscriptstyle 01}$ is $-(\bar{g}_{\scriptscriptstyle 01}^*\TM,\theta_i)^\g+dA_{\scriptscriptstyle 01}+\omega^2_{\scriptscriptstyle 01}$, which we denote by $\beta_{\scriptscriptstyle 01}^{skew}$. Obviously, we have
 $$\langle L_X\beta_{\scriptscriptstyle 01}^{skew}(Y)-\ii_Y d\beta_{\scriptscriptstyle 01}^{skew}(X)-\beta_{\scriptscriptstyle 01}^{skew}([X,Y]),Z\rangle=d\beta_{\scriptscriptstyle 01}^{skew}(X,Y,Z).$$
Furthermore, we have
\begin{equation}\label{eq:csji}
    \cs_3(\theta_1)- \cs_3(\theta_0)=d(\theta_0,g_{01}^*\TM)^\g-\frac{1}{6}\bar{g}_{01}^*\frkC,
\end{equation}
where $\frkC\in\Omega^3(G)$ defined by
$$
\frac{1}{6}\frkC(\hat{a},\hat{b},\hat{c})=(\TM(\hat{a}),[\TM(\hat{b}),\TM(\hat{c})]_\g)^\g=(a,[b,c]_\g)^\g,\quad \forall a, b, c\in\g.
$$
Here $\hat{a},\hat{b},\hat{c}$ are right invariant vector fields on  $G$.
%The proof of \eqref{eq:csji} is given in Appendix \ref{app:csji}.

 By \eqref{eq:diffcs} and \eqref{eq:csji}, we obtain
 \begin{equation}
 d\beta_{\scriptscriptstyle 01}^{skew}=\cs_3(\theta_1)-\cs_3(\theta_0)-d(\bar{g}_{\scriptscriptstyle 01}^*\TM,\theta_i)^\g=-\frac{1}{6}\bar{g}_{\scriptscriptstyle 01}^*\frkC.
\end{equation}

 On the other hand, by straightforward computations, we have
 $$
 \langle L_X\beta_{\scriptscriptstyle 01}^{sym}(Y)-\ii_Y d\beta_{\scriptscriptstyle 01}^{sym}(X)-\beta_{\scriptscriptstyle 01}^{sym}([X,Y])+\huaP(\bar{g}_{\scriptscriptstyle 01}^*\TM,\bar{g}_{\scriptscriptstyle 01}^*\TM),Z\rangle =\bar{g}_{\scriptscriptstyle 01}^*\frkC(X,Y,Z),
 $$
 which implies that
 \eqref{eq:auto4}  holds.

Finally, \eqref{eq:diffcs} implies that the conditions in \eqref{eq:congluetc1}-\eqref{eq:congluetc3} are satisfied. Thus $\Lambda_{\scriptscriptstyle 01}$ is a $1$-morphism in $\tcalgdp$.
It finishes the proof. \qed
\vspace{3mm}

Recall that the model we use
for a
$2$-groupoid is a simplicial set satisfying Kan conditions $\Kan(n, j)$
for all $n\ge 1$ and $0\le j \le n$ and strict
Kan conditions $\Kan!(n, j)$ for all $n\ge 3$ and $0\le j \le n$.
Lemma
\ref{1morphism} verifies that
$\Phi(U):\bstringnp(U)_1\rightarrow\tcalgdp(U)_1$ is
well-defined. The following lemma will verify that the map
$\Phi(U):\bstringnp(U)_2\rightarrow\tcalgdp(U)_2$ on $2$-simplices commutes with the face maps.
\begin{lem}\label{Catmorphism}
Given any $2$-morphism $(f,\omega^1)$ between $T_{\scriptscriptstyle 01}\circ T_{\scriptscriptstyle 12}$ and $T_{\scriptscriptstyle 02}$, where $T_{ij}:=(\bg_{ij}, A_{ij}, \omega^2_{ij})_{\{0\le i<j\le 2\}}\in \bstringnp(U)_1$ are $1$-morphisms between the objects $(U\times \String(n)\to U, \theta_i, B_i)_{i=1,2,3}\in \bstringnp(U)_0$, the corresponding images $\Lambda_{ij}:=\Phi(T_{ij})=\Phi(\bg_{ij}, A_{ij}, \omega^2_{ij})_{\{0\le i<j\le 2\}}$ under the morphism $\Phi$ satisfy the condition $\Lambda_{\scriptscriptstyle 01}\Lambda_{\scriptscriptstyle 12}=\Lambda_{\scriptscriptstyle 02}$, i.e. we have the following commutative diagram:
\[\begin{tikzpicture}[>=latex',mydot/.style={draw,circle,inner
    sep=1pt},every label/.style={scale=1},scale=1]

  \foreach \i in {0}{
  \node[mydot, fill=black, label=240:$({\theta_0}{,}{B_0})$]                at (-1+6*\i, -1.73)    (p\i0) {};
  \node[mydot,fill=black,label=90:$({\theta_1}{,}{B_1})$]      at (+6*\i,0) (p\i1) {};
  \node[mydot,fill=black,label=-60:$({\theta_2}{,}{B_2}).$]     at (1+6*\i,-1.73)    (p\i2) {};
}
\begin{scope}[<-]
    %--fig 0
    \draw (p00) --node[left,scale=.8]{$\Lambda_{01}$} (p01);
    \draw (p00) --node[below,scale=.8]{$\Lambda_{02}$} (p02);
    \draw (p01) --node[right,scale=.8]{$\Lambda_{12}$} (p02);
\end{scope}
 \node at (0,-1){\tiny{${\rm Id}$}};
 % \node at (20, -1.8){\footnotesize{$\Kan(2,2)$,}};
\end{tikzpicture}\]

\end{lem}
\pf By $\bg_{\scriptscriptstyle 01}\bg_{\scriptscriptstyle 12}=\bg_{\scriptscriptstyle 02}$, we have
\begin{eqnarray*}
  \Lambda_{\scriptscriptstyle 01}\Lambda_{\scriptscriptstyle 12}&=&
  \left(\begin{array}{ccc}
1&0&0\\
-\bg_{\scriptscriptstyle 01}^*\TM&\ad_{\bg_{\scriptscriptstyle 01}}&0\\
\beta_{\scriptscriptstyle 01}&2(\bg_{\scriptscriptstyle 01}^*\TM)^\star\circ \ad_{\bg_{\scriptscriptstyle 01}}&1
\end{array}\right)
\left(\begin{array}{ccc}
1&0&0\\
-\bg_{\scriptscriptstyle 12}^*\TM&\ad_{\bg_{\scriptscriptstyle 12}}&0\\
\beta_{\scriptscriptstyle 12}&2(\bg_{\scriptscriptstyle 12}^*\TM)^\star\circ \ad_{\bg_{\scriptscriptstyle 12}}&1
\end{array}\right) \\
&=&\left(\begin{array}{ccc}
1&0&0\\
-\bg_{\scriptscriptstyle 01}^*\TM-\ad_{\bg_{\scriptscriptstyle 01}}\bg_{\scriptscriptstyle 12}^*\TM&\ad_{\bg_{\scriptscriptstyle 01}}\ad_{\bg_{\scriptscriptstyle 12}}&0\\
D_{\scriptscriptstyle 31}&D_{\scriptscriptstyle 32} &1
\end{array}\right),
\end{eqnarray*}
where
\begin{eqnarray*}
  D_{\scriptscriptstyle 32}&=&2(\bg_{\scriptscriptstyle 01}^*\TM)^\star\circ \ad_{\bg_{\scriptscriptstyle 01}}\circ\ad_{\bg_{\scriptscriptstyle 12}}+2(\bg_{\scriptscriptstyle 12}^*\TM)^\star\circ \ad_{\bg_{\scriptscriptstyle 12}}\\
  &=&2(\bg_{\scriptscriptstyle 01}^*\TM+\ad_{\bg_{\scriptscriptstyle 01}}\bg_{\scriptscriptstyle 12}^*\TM)^\star\circ\ad_{\bg_{\scriptscriptstyle 02}} \\
  &=&2(\bg_{\scriptscriptstyle 02}^*\TM)^\star\circ\ad_{\bg_{\scriptscriptstyle 02}},
\end{eqnarray*}
and
\begin{eqnarray*}
  D_{\scriptscriptstyle 31}&=& \beta_{\scriptscriptstyle 01}-2(\bg_{\scriptscriptstyle 01}^*\TM)^\star\circ \ad_{\bg_{\scriptscriptstyle 01}}\circ\bg_{\scriptscriptstyle 12}^*\TM+\beta_{\scriptscriptstyle 12}\\
  &=&-(\bar{g}_{\scriptscriptstyle 01}^*\TM,\theta_0)^\g+dA_{\scriptscriptstyle 01}+\omega^2_{\scriptscriptstyle 01}-(\bar{g}_{\scriptscriptstyle 01}^*\TM)^\star \circ \bar{g}_{\scriptscriptstyle 01}^*\TM
  -2(\bg_{\scriptscriptstyle 01}^*\TM)^\star\circ \ad_{\bg_{\scriptscriptstyle 01}}\circ\bg_{\scriptscriptstyle 12}^*\TM\\
  &&
  -(\bar{g}_{\scriptscriptstyle 12}^*\TM,\theta_1)^\g+dA_{\scriptscriptstyle 12}+\omega^2_{\scriptscriptstyle 12}-(\bar{g}_{\scriptscriptstyle 12}^*\TM)^\star \circ \bar{g}_{\scriptscriptstyle 12}^*\TM\\
  &=&-(\bg_{\scriptscriptstyle 01}^*\TM+\ad_{\bg_{\scriptscriptstyle 01}}\bg_{\scriptscriptstyle 12}^*\TM,\theta_0)^\g+dA_{\scriptscriptstyle 02}+\omega^2_{\scriptscriptstyle 02}-(\bar{g}_{\scriptscriptstyle 01}^*\TM)^\star \circ (\bar{g}_{\scriptscriptstyle 01}^*\TM
  + \ad_{\bg_{\scriptscriptstyle 01}}\circ\bg_{\scriptscriptstyle 12}^*\TM)\\
  &&-(\ad_{\bg_{\scriptscriptstyle 01}}\circ\bg_{\scriptscriptstyle 12}^*\TM)^\star\circ \bg_{\scriptscriptstyle 01}^*\TM-(\ad_{\bg_{\scriptscriptstyle 01}}\circ\bar{g}_{\scriptscriptstyle 12}^*\TM)^\star \circ\ad_{\bg_{\scriptscriptstyle 01}}\circ \bar{g}_{\scriptscriptstyle 12}^*\TM\\
  &=&-(\bg_{\scriptscriptstyle 02}^*\TM,\theta_0)^\g+dA_{\scriptscriptstyle 02}+\omega^2_{\scriptscriptstyle 02}-(\bar{g}_{\scriptscriptstyle 01}^*\TM)^\star \circ \bar{g}_{\scriptscriptstyle 02}^*\TM-(\ad_{\bg_{\scriptscriptstyle 01}}\circ\bg_{\scriptscriptstyle 12}^*\TM)^\star\circ \bg_{\scriptscriptstyle 02}^*\TM\\
  &=&\beta_{\scriptscriptstyle 02}.
\end{eqnarray*}
Therefore, we have $\Lambda_{01}\Lambda_{12}=\Lambda_{02}$.\qed
\\
\\
{\bf Proof of Theorem \ref{Naturaltransformation}}. Recall that given any manifold $U$, a morphism between the $2$-groupoids $\bstringnp(U)$ and $\tcalgdp(U)$ is a simplicial morphism of the underlying simplicial sets.
Lemma
\ref{1morphism} and \ref{Catmorphism} verify that
$\Phi(U):\bstringnp(U)\rightarrow\tcalgdp(U)$ is indeed a morphism of the underlying simplicial sets. Hence we only need to prove the naturality of the map $\Phi$. To do this, let us assume that $r_{V,U}:V\rightarrow U$ is a smooth map between two manifolds $V$ and $U$. Then $r_{V,U}$ induces maps
\begin{itemize}
\item $\tcalgdp(r_{V,U}):\tcalgdp(U)\rightarrow \tcalgdp(V)$ given by
\[(\tagu,\theta, B)\mapsto (\tagv, r_{V,U}^*(\theta), r_{V,U}^*(B)), \]
\item $\bstringnp(r_{V,U}):\bstringnp(U)\rightarrow \bstringnp(V)$  given by
\[((U, \theta, B);
g, A, \omega^2; f, \omega^1)\mapsto (V, r_{V,U}^*(\theta), r_{V,U}^*(B));
r_{V,U}^*(g), r_{V,U}^*(A), r_{V,U}^*(\omega^2); r_{V,U}^*(f), r_{V,U}^*(\omega^1))
.\]
\end{itemize}

Here $(U, \theta, B)$ stands for an element $(U\times \String(n)\to U,\theta,B)\in \bstringnp(U)_1.$
It is then straightforward to obtain the following commutative diagram
$$
\begin{CD}
\bstringnp(U) @> \Phi(U) >> \tcalgdp(U) \\
@V \bstringnp(r_{V,U}) VV @V \tcalgdp(r_{V,U}) VV \\
\bstringnp(V) @> \Phi(V) >> \tcalgdp(V),
\end{CD}
$$
which shows the naturality of the map $\Phi$. Thus $\Phi$ is morphism from $\bstringnp$ to $\tcalgdp$.
\qed
\vspace{3mm}

As a corollary, we have proven the desired result,

\begin{cor}\label{Functor}
There is a natural morphism $\Phi$ from the $(3,1)$-sheaf
${\bstringnp}^+$ to the $(2,1)$-sheaf ${\tcalgdp}^+$.
\end{cor}

By the discussion above, the morphism $\Phi$ can be described explicitly as follows. Given any manifold $M$, the morphism $\Phi:{\bstringnp}^+(M)\rightarrow {\tcalgdp}^+(M)$ is given on the $0$-, $1$- and $2$-simplices respectively by
\begin{itemize}
    \item
on $0$-simplices
\[\Phi(\{U_i\},P_{c})=\Big(\{U_i\}, (\sqcup (\tagui,  \theta_i, B_i);-\bar{g}_{ij}^*\TM,\ad_{\bar{g}_{ij}},\beta_{ij} ) \Big)\]
where $\{U_i\}$ is an open cover of $M$ and $P_{c}:=(\sqcup (U_i\times \String(n)\to U_i, \theta_i, B_i;g_{ij}, A_{ij}, \omega^2_{ij}; f_{ijk}, \omega^1_{ijk})$ is string data, $\bar{g}_{ij}:U_{ij} \to G$ is the underlying morphism of $g_{ij}$,  and $\beta_{ij}$ is given by $$\beta_{ij}=-(\bar{g}_{ij}^*\TM,\theta_i)^\g+dA_{ij}+\omega^2_{ij}-(\bar{g}_{ij}^*\TM)^\star \circ \bar{g}_{ij}^*\TM.$$
\item on $1$-simplices
\[\Phi(\{V_i\},\phi_{c})=(\tagvi, \Lambda_i),\]
where $\{V_i\}$ is a common refinement of $\{U_i\}$, $\{\widetilde{U}_i\}$, and $\phi_{c}:=(g_i:V_i\rightarrow \String(n), A_i,\omega^2_i)$ provides a $1$-morphism between $(\{U_i\},P_{c})$ and $(\{\widetilde{U}_i\},\widetilde{P}_{c})$, $\Lambda_i$ is the $1$-morphism from $(\tagvi,\widetilde{\theta}_i,\widetilde{B}_i)$ to  $(\tagvi,\theta_i,B_i)$ defined as before
$$\Lambda_{i}=\left(\begin{array}{ccc}
1&0&0\\
-{\bar{g}_{i}}^*\TM&\ad_{\bar{g}_{i}}&0\\
\beta_{i}&2(\bar{g}_{i}^*\overline{\TM})^\star&1
\end{array}\right),$$
where $\bar{g}_{i}:V_{i}\rightarrow G$ is the underlying morphism of $g_i$ and $\beta_i:TV_{i}\longrightarrow T^*V_{i}$ is given by
$$\beta_i=-(\bar{g}_{i}^*\TM,\theta_i)^\g+dA_{i}+\omega^2_{i}-(\bar{g}_{i}^*\TM)^\star \circ \bar{g}_{i}^*\TM;$$
\item on $2$-simplices
\[\Phi(\{W_i\},\alpha_{c})=1,\]
where $\{W_i\}$ is a common refinement of $\{V_i\}$, $\{\widetilde{V}_i\}$, and $\alpha_{c}$ provides a $2$-morphism between $(\{V_i\},\phi_{c})$ and $(\{\widetilde{V}_i\},\widetilde{\phi}_{c})$.
\end{itemize}

As an object in ${\tcalgdpp(M)}$ glues to a Courant algebroid by the  discussion in Appendix \ref{GluingCA},  let us describe explicitly the Courant algebroid with a connection associated to a $\String(n)$-principal bundle with  connection data on a manifold $M$.

\emptycomment{
By the discussion in Appendix \ref{GluingCA}, any object in ${\tcalgdpp(M)}$ takes a
form of $$(TM\oplus \huaG\oplus
T^*M, \Courant{-,-}^T_{\nabla,R,H}, \pair{-,-}^T, \pr_{TM}).$$ Furthermore, Corollary \ref{Functor} directly leads to the following fact:} Given a $\String(n)$ data
$$P_{c}=(\sqcup U_i\times \String(n)\to \sqcup U_i, \theta_i, B_i;g_{ij}, A_{ij}, \omega^2_{ij}; f_{ijk}, \omega^1_{ijk})$$
over a cover $\{U_i\}$ of $M$ in $\bstringnp$, the corresponding transitive Courant algebroid with a connection under the morphism $\Phi:{\bstringnp}^+(M)\rightarrow {\tcalgdp(M)}^+$ is
$$E\cong (TM\oplus \huaG\oplus T^*M, \Courant{-,-}^T_{\nabla,R,H}, \pair{-,-}^T, \pr_{TM}),$$
equipped with  bracket $\Courant{-,-}^T_{\nabla,R,H}$ and pairing $\pair{-,-}^T$ as in formula \eqref{eq:Tbracket}, \eqref{eq:Tpair} with  the global 2-form $R\in \Omega^2(M, \huaG)$, the connection $\nabla:\Gamma(TM)\otimes \Gamma(\huaG)\longrightarrow \Gamma(\huaG)$  and the global $3$-form $H\in \Omega^3(M)$ defined on each $U_i$ respectively by
\emptycomment{
\begin{eqnarray*}
    \nonumber\Courant{X+a+\xi,Y+b+\eta}^T_{\nabla,R,H}&=&[X,Y]+\nabla_Xb-\nabla_Ya+[a,b]_\huaG+R(X,Y)\\
  &&+L_X\eta-i_Yd\xi+P(a,b)-2Q(X,b)+2Q(y,a)+H(X,Y),\\
  \pair{X+a+\xi,Y+b+\eta}^T&=&\half\big(\xi(Y)+\eta(X)\big)+(a,b)^\huaG.
\end{eqnarray*}}
\begin{eqnarray*}
\nabla^i_Xa &:=&X(a)+[\theta_i(X),a]_\g, \quad \forall X\in\Gamma(TU_i), a\in\Gamma(U_i\times \g), \\
R_i &:=& d\theta_i + \half [\theta_i, \theta_i]_\g, \\
 H_i&:=&\cs_3(\theta_i)-dB_i.
\end{eqnarray*}
\emptycomment{
and $P:\Gamma(\huaG)\otimes \Gamma(\huaG)\longrightarrow\Omega^1(M)$ and $Q:\frkX(M)\otimes \Gamma(\huaG)\longrightarrow\Omega^1(M)$ are given by
\begin{eqnarray*}
  P(a,b)(Y)&=&2(b,\nabla_Ya)^\huaG,\\
  Q(X,a)(Y)&=&(a,R(X,Y))^\huaG.
\end{eqnarray*}}

\subsection{Property of the morphism $\Phi: \bstringnpp \to \tcalgdpp$}

We denote the $(3,1)$-presheaf of $U(1)$-gerbes (or $\B U(1)$-principal bundles) with connection data by $\B \B U(1)^p_c$. Then, induced by the morphism  $\B U(1) \xrightarrow{\iota} \String(n)$, there is a morphism $\B \B U(1)^p_c \xrightarrow{\B\iota } \bstringnp$ given by
\[
\begin{split}
(U\times \B U(1), B\in \Omega^2(U); L: U\to \B U(1), A\in \Omega^1(U); a: U\to U(1)) \mapsto \\
(U\times \String(n), B, \theta=0; g=\iota \circ L, A, \omega^2=0; a, \omega^1=0).
\end{split}
\]
On the side of algebroids, we have similar constructions. Let us denote by $\ecalgdp$ the $(2,1)$-presheaf of exact Courant algebroids (see Appendix \ref{app:ec} for the definition and notations). There is a morphism on the presheaf level, $\ecalgdp \to \tcalgdp$,  given by
\[
(\T U, \Courant{-,- }^E_S, \langle-, - \rangle^E, \pr_{TU}, B; \Big( \begin{matrix} 1 & 0 \\ \huaB & 1 \end{matrix} \Big) ) \mapsto (\tagu, \Courant{-,-}^T_S, \langle-, - \rangle^T,  \pr_{TU},0, B; \left( \begin{matrix} 1 & 0 & 0 \\ 0  & 1 &0 \\ \huaB & 0 & 1 \end{matrix} \right)).
\]
Let us denote by $\tcalgdp$ the $(2,1)$-presheaf of transitive Lie algebroids (see Appendix \ref{app:tl}). Similar to the morphism $\bstringnp \xrightarrow{\pi}\mathsf{B}\Spin(n)_{c}^p$ described before Theorem \ref{thm:lift}, we have a morphism $\tcalgdp \to \tlp$ defined also by forgetting some data,
\[
(\tagu, \Courant{-, -}^T_S, \langle-, -\rangle^T, \pr_{TU}; \theta, B;  \left( \begin{matrix} 1 & 0 & 0 \\ \phi & \tau &0 \\ \beta & -2\phi^\star\circ \tau & 1 \end{matrix} \right) ) \mapsto (TU\oplus (U\times{\g}), [-,-]^T_S,  \pr_{TU},\theta; \Big( \begin{matrix} 1 & 0 \\ \phi & \tau  \end{matrix} \Big) ).
\]
Then we can also construct morphisms $\Phi_{U(1)}: \B \B U(1)_c \to \ecalgdpp$ and $\Phi_{\Spin(n)}: \B \Spin(n)_c \to \tlpp$, for $\g=\so(n)$. The constructions are essentially given in \cite{sevlet,MK2}, and here we  spell it out in our setting. On the level of objects, they are given by
\begin{equation}
\begin{split}
\Phi_{U(1)}: & (U_i\times \B U(1), B_i; L_{ij}, A_{ij}; a_{ijk} ) \mapsto (TU_i \oplus T^*U_i, B_i; \Big( \begin{matrix} 1 & 0 \\ dA_{ij} & 1 \end{matrix} \Big)  ), \\
\Phi_{\Spin(n)}: & (U_i \times \Spin(n), \theta_i; g_{ij} ) \mapsto (TU_i\oplus (U_i\times{\so(n)}), [-,-]^T_S,  \pr_{TU},\theta_i; \Big( \begin{matrix} 1 & 0 \\ -g_{ij}^*\TM & \ad_{g_{ij}}  \end{matrix} \Big) ),
\end{split}
\end{equation}
and on the level of morphisms, they are given by corresponding pullbacks and pre-compositions.

Thus, we have the following commutative diagram to connect the principal bundle side and the algebroid side,
\begin{equation}
\xymatrix{
\B \B U(1)_c \ar[r]^{\Phi_{U(1)}} \ar[d]^{\B\iota}& \ecalgdpp \ar[d]\\
\bstringnpp \ar[r]^\Phi \ar[d]^{\pi} & \tcalgdpp \ar[d]\\
\mathsf{B}\Spin(n)_{c} \ar[r]^{\Phi_{\Spin(n)}} & \tlpp.
}
\end{equation}

Now we show that $\Phi_{U(1)}$ is not injective in general and this implies that $\Phi$ is not injective. And $\Phi_{\Spin(n)}$ is not surjective thus $\Phi$ can not be surjective. Moreover, even on the fibre of the image of $\Phi_{\Spin(n)}$, $\Phi$ can not be surjective in general.

\begin{lem}\label{lem:gerbes}
When $H^3(M, \Z)$ has torsion, $\Phi_{U(1)}$ is not injective on the level of objects and not fully faithful.
\end{lem}
\pf
We take a cocycle $a^0_{ijk}$ representing a torsion element in $H^3(M, \Z)$, lifting it to a Deligne cocycle $(a^0_{ijk}, A^0_{ij}, B^0_i)$, then $dB^0_i$ glues to an exact 3-form. We now show that such a Deligne cocycle may always be  adjusted by an exact one to $(a^0_{ijk}, A^0_{ij}, B^0|_{U_i})$ for a global 2-form $B^0$ and some closed $A^0_{ij}$. Since $dB^0_i$ glues to an exact 3-form, there is a global 2-form $B^0$ such that $B^0_i=B|_{U_i} + dA^0_i$. Then adjusting the original Deligne cocycle by $D(1, A^0_i)$ will fulfill our aim.

Therefore we might as well assume that we lift $a^0_{ijk}$ to a Deligne cocycle $(a^0_{ijk}, A^0_{ij}, B^0|{U_i})$ satisfying $dA^0_{ij}=0$. Then clearly the image of the two objects $(U_i \times \B U(1), B_i;L_{ij},  A_{ij}; a_{ijk} ) \in \B \B U(1)$ and $(U_i \times \B U(1), B_i+B^0_i;L_{ij},  A_{ij}+A^0_{ij}; a_{ijk}+a^0_{ijk} ) \in \B \B U(1)$ are the same. However, it is clear that there exists no morphism between these two objects because $a^0_{ijk}$ is not exact.
\qed

\begin{lem}
In general, the map $\Phi$ is not injective on the level of objects and not fully faithful .
\end{lem}
\pf
We take the two different gerbes with connection data constructed in Lemma \ref{lem:gerbes},  $\huaG_1$ and $\huaG_2$,  which maps to the same object under $\Phi_{U(1)}$. Then we see that $\B\iota (\huaG_1)$ and $\B\iota (\huaG_2)$ are  non-isomorphic string data  but mapping to the same Courant data on the right hand side.
\qed

\begin{lem}\label{lem:curvature}
The map $\Phi$ preserves curvatures.
\end{lem}
\pf
It is clear from the definition of curvatures on both sides.
\qed

\begin{lem}
The map $\Phi$ is not essentially surjective in general\footnote{We thank very much Pavol \v{S}evera for pointing out this example to the third author.}.
\end{lem}
\pf
The map $\Phi_{\Spin(n)}$ is not essentially surjective because there are non-integrable transitive Lie algebroids. To show $\Phi$ is not essentially surjective, we need to find a non-integrable transitive Lie algebroid $A$ whose $p_1(A)=0$. Notice that integrability is a property preserved by isomorphisms of Lie algebroids.

We take $M=\R^3-\{p_1\}-\{p_2\}$ where $p_1, p_2\in \R^3$ are two different points and $A=TM\times \R$, with the following Lie bracket
\[
[(X, f), (Y, g)]=[X, Y]+ X(f)-Y(g)+\omega(X, Y).
\] Here $\omega = \iota_1^* \pr_1^* \omega_a + \sqrt{2} \iota_2^* \pr_2^* \omega_a$ with $\iota_j :  M\to \R^3- p_j  $,  $\pr_j: \R^3- p_j  \to S^2$ for $j=1,2$,  and $\omega_a$ the area form on $S^2$. Then the period $\{ \int_{\gamma}\omega, \gamma \in \pi_2(M)\} $ of $\omega$ is dense in $\R$. Therefore, $A$ is not integrable by \cite{cf}. On the other hand, it is clear $p_1(A)=0$ because there is no non-trivial 4-forms on $M$.
\qed

\begin{lem}
On the fibre of the image of $\Phi_{\Spin(n)}$, $\Phi$ is in general not essentially surjective.
\end{lem}
\pf As pointed out in Lemma \ref{lem:curvature}, $\Phi$ preserves the curvature. If we fix an object $M\xrightarrow{\bar{P}_c} \B \Spin(n)_c$ and look at all possible lifts of $\bar{P}_c$, we see that curvatures for these lifts form a torsor of $\im \bar{d}$ as in Corollary \ref{cor:curv}. Now on the Courant side, fixing the underlying $\Spin(n)$-principal bundle and its connection data means that we consider all possible Courant lifts over a fixed Atiyah Lie algebroid and its connection data. By Corollary \ref{cor:curv-courant}, the set of curvatures for such Courant lifts is a torsor of $\Omega^3_{\cl}(M)$. As pointed out in Remark \ref{rmk:curv}, $H^2(M, D_2) \xrightarrow{\bar{d}} \Omega_{\cl}^3(M)$ is not surjective in general. Since curvatures are preserved under isomorphisms, we see that $\Phi$ can not be essentially surjective in general even on the fibre of the image of $\Phi_{\Spin(n)}$.
\qed

\emptycomment{
\begin{rmk}\label{}
We may extend our construction of $\Phi$ from $\Spin(n)$ to a general Lie group. Then the action Courant algebroid $M \times \gl(n)$ over the complete flag $M$,  studied in Example \ref{ep:action-courant},  stays in the image of $\Phi$. This is because $M\times \gl(n)$ is isomorphic to a standard transitive Courant algebroid, whose preimage under $\Phi$ can be a trivial string principal bundle with all connection forms ($\theta_i$, $A_i$, $B_{ij}$, and all $\omega^l$'s) 0, $g_{ij}=1$, and $f_{ijk}=1$.
\end{rmk}}

\appendix

\section{Appendix}\label{sec:app}

\subsection{$(2,1)$-sheaf $\tlp^+$ of transitive Lie algebroids with connections}\label{app:tl}
\begin{defi}
A {\bf Lie algebroid} structure on a vector bundle $ A\longrightarrow M$ is
a pair that consists of a Lie algebra structure $[-,-] $ on
the section space $\Gamma( A)$ and a  bundle map
$\rho: A\longrightarrow TM$, called the anchor, such that the
following relation is satisfied:
$$~[X,fY] =f[X,Y]+\rho(X)(f)Y,\quad \forall~X,Y\in\Gamma(A),~f\in
\CWM.$$
\end{defi}

A Lie algebroid $A$ is called transitive if $\rho$ is surjective, i.e. $\im\rho=TM$. Denote by $\huaG=\ker\rho$. Then $\huaG$ is a bundle of Lie algebras, whose fibre is isomorphic to a Lie algebra $(\g,[-,-]_\g)$. We have the following short exact sequence:
$$
0\longrightarrow\huaG\longrightarrow A\stackrel{\rho}{\longrightarrow}TM\longrightarrow 0.
$$
 A splitting $s:TM\longrightarrow A$   gives rise to a connection $\nabla$ on $\huaG$ by
 $$
 \nabla_Xa=[s(X),a], \quad \forall~X\in\Gamma(A), a\in\Gamma(\huaG).
 $$ Thus we call such a splitting $s: TM\longrightarrow A$ a {\bf connection} of $A$. Connections always exist by partition of unity. Thus, after picking a connection, we have
 $A \cong TM\oplus \huaG$, and the induced bracket on $TM\oplus \huaG$ is
\begin{equation}\label{br:tlc}
\brtl{X+a,Y+b}=[X,Y]+\nabla_Xb-\nabla_Ya+[a,b]_\g+R(X,Y), \quad \forall~X,Y\in\Gamma(TM),~a, b\in\Gamma(\huaG),
\end{equation}
where $R(X,Y)=[s(X),s(Y)]-s([X,Y])$ is the {\bf curvature} of the connection  $s$.  In other words, a transitive Lie algebroid with a connection is always isomorphic to $(TM\oplus \huaG,\brtl{\cdot,\cdot},\rho=\pr_{TM})$ and the isomorphism depends on the choice of the connection.

In particular, if $\huaG=M\times \g$ is a trivial bundle and the connection $\nabla$ is given by the flat connection $\nabla_Xb=X(b)$, we obtain the standard bracket
\begin{equation}\label{br:tls}
\brtls{X+a,Y+b}=[X,Y]+ X(b)-Y(a)+[a,b]_\g.
\end{equation}

An automorphism of the standard transitive Lie algebroid is given by a matrix
$
\left(\begin{array}{cc}
1&0\\
\phi&\tau
\end{array} \right),
$
where $\tau:M\longrightarrow\Aut(\g)$ and $\phi\in \Omega^1(M,\g)$ satisfy
\begin{eqnarray*}
  \phi([X,Y])&=&X(\phi(Y))-Y(\phi(X))+[\phi(X),\phi(Y)]_\g,\\
  \tau([a,b]_\g)&=&[\tau(a),\tau(b)]_\g,\\
  \tau(X(b))&=& X(\tau(b))+[\phi(X),\tau(b)]_\g.
\end{eqnarray*}

There is a $(2,1)$-presheaf of transitive Lie algebroids with
connections $\tlp: \Mfd^{\op} \to \Gpd$, where $\Mfd^{\op}$ is the
opposite  category of
$\Mfd$, and $\Gpd$ is the
$2$-category of (discrete) groupoids and groupoid morphisms.

For an object $U\in \Mfd$, the groupoid $\tlp(U)$ is made up
by the following data:
\begin{itemize}
\item{$\tlp(U)_0$}: an object is a quadruple $(TU\oplus (U\times\g), \brtls{-,-},\pr_{TU}, \theta)$, where $\theta\in \Omega^1(U,\g)$ is a $\g$-valued $1$-form and $\brtls{-,-}$
  is the standard   bracket   given by \eqref{br:tls}. We will simply denote an object by $(TU\oplus (U\times\g),\theta)$.

\item{$\tlp(U)_1$}:  a $1$-morphism from $(TU\oplus (U\times\g),   \widetilde{\theta})$ to $(TU\oplus (U\times\g),
   \theta)$ is an automorphism of the standard transitive Lie algebroid $(TU\oplus
  (U\times\g),\brtls{-,-},\pr_{TU})$ given by the
  matrix
$
\left(\begin{array}{cc}
1&0\\
\phi&\tau
\end{array} \right),
$ such that
$$
\theta(X)-\tau(\widetilde{\theta}(X))=\phi(X).
$$

 The
composition of $1$-morphisms is simply the matrix multiplication.

\end{itemize}

Then for a morphism $\varphi: U\to V$ in $\Mfd$, the associated
functor $\tlp(\varphi): \tlp(V) \to \tlp(U)$ is induced by
pulling back forms.
\calc{
more precisely, on the level of objects,
\[\tlp(\varphi)(TV\oplus (V\times\g),   \brtls{-,-},\theta_V)
= (TU\oplus (U\times\g), \brtls{-,-}, \varphi^*\theta_V)\]
and on the level of morphisms,
\[
\tlp(\varphi) \left(\begin{array}{cc}
1&0\\
\phi_V&\tau_V
\end{array} \right) = \left(\begin{array}{cc}
1&0\\
\varphi^*\phi_V&\tau_V\circ \varphi
\end{array} \right) .
\]
One may verify that $\tlp(\varphi)$ is indeed a functor between
desired categories.}
%  \end{pro}

Take an open cover $\{U_i\}$ of $M\in
\Mfd$. An object in  $\holim \tlp (U(M)_\bullet)$ consists of
\begin{itemize}
\item  an object $\sqcup (TU_i\oplus (U_i\times\g),
   \theta_i)$ in $\tlp(\sqcup U_i)_0$,

\item  $\Lambda_{ij}=\left(\begin{array}{cc}
1&0\\
\phi_{ij}&\tau_{ij}
\end{array} \right) \in \tlp(\sqcup U_{ij})_1$, which is a $1$-morphism from
$
(TU_{ij}\oplus (U_{ij}\times \g),
    \theta_j|_{U_{ij}})$ to
 $ (TU_{ij}\oplus (U_{ij}\times \g),  \theta_i|_{U_{ij}}),$
therefore satisfying
\begin{eqnarray}
\label{eq:tl1}\phi_{ij}([X,Y])&=&X(\phi_{ij}(Y))-Y(\phi_{ij}(X))+[\phi_{ij}(X),\phi_{ij}(Y)]_\g,\\
 \label{eq:tl2} \tau_{ij}([a,b]_\g)&=&[\tau_{ij}(a),\tau_{ij}(b)]_\g,\\
 \label{eq:tl3} \tau_{ij}(X(b))&=& X(\tau_{ij}(b))+[\phi_{ij}(X),\tau_{ij}(b)]_\g.\\
  \label{eq:generalcontheta}\theta_i|_{U_{ij}}-\tau_{ij}\theta_j|_{U_{ij}}&=&\phi_{ij}.
\end{eqnarray}
\item compatibility condition $\Lambda_{ij}\Lambda_{jk}=\Lambda_{ik}$ on $U_{ijk}$, which unpacks itself to the following two equations
    \begin{eqnarray*}
      \phi_{ij}+\tau_{ij}\phi_{jk}&=&\phi_{ik},\\
      \tau_{ij}\tau_{jk}&=&\tau_{ik}.
    \end{eqnarray*}
\end{itemize}

\begin{defi}\label{Liedata}
We call an object $\big(\sqcup(TU_i\oplus(U_i\times\g),\theta_i),\phi_{ij},\tau_{ij}\big)$ in  $\holim \tlp (U(M)_\bullet)$ a {\bf transitive Lie data}.
\end{defi}

 Given a transitive Lie data  $\big(\sqcup(TU_i\oplus(U_i\times\g),\theta_i),\phi_{ij},\tau_{ij}\big)$, since $\Lambda_{ij}=\left(\begin{array}{cc}
1&0\\
\phi_{ij}&\tau_{ij}
\end{array} \right) $ satisfies the cocycle condition $\Lambda_{ij}\Lambda_{jk}=\Lambda_{ik}$, we can glue $TU_i\oplus (U_i\times\g)$'s and obtain a  vector bundle
\begin{equation} \label{eq:glue-A}
A=\coprod TU_i\oplus (U_i\times\g)/\thicksim,
\end{equation}
where the equivalence relation $\thicksim$ is given by
$$
X+a\thicksim Y+b \Longleftrightarrow \left(\begin{array}{c}Y\\b\end{array}\right)=\Lambda_{ij}\left(\begin{array}{c}X\\a\end{array}\right),\quad\forall ~X+a\in TU_j\oplus (U_j\times\g),~Y+b\in TU_i\oplus  (U_i\times \g).
$$
 \begin{pro}\label{pro:tlo} \label{thm:tlo}
   A transitive Lie data  $\big(\sqcup(TU_i\oplus(U_i\times\g),\theta_i),\phi_{ij},\tau_{ij}\big)$  gives rise to a transitive Lie algebroid $(A,[-,-],\rho)$ with a connection $s:TM\longrightarrow A$.
\end{pro}
 \pf Obviously, the vector bundle $A$
 fits the following short exact sequence:
$$
0\stackrel{}{\longrightarrow}\huaG\stackrel{}{\longrightarrow}A\stackrel{\rho}{\longrightarrow}TM\stackrel{}{\longrightarrow}0,
$$
where $\huaG$ denotes the Lie algebra bundle obtained from the transition function $\tau_{ij}$ and $\rho$ is induced by the projection $TU_i\oplus (U_i\times\g)\longrightarrow TU_i$.

Since $\Lambda_{ij}$ preserves the Lie bracket $\brtls{-,-},$ there is a well-defined Lie bracket $[-,-]$ on $\Gamma(A)$. Then we obtain a Lie algebroid $(A,[-,-],\rho)$.

On $U_i$, consider the splitting $s_i:TU_i\longrightarrow A|_{U_i}$ given by
$$
s_i(X)=X+\theta_i(X).
$$
By \eqref{eq:generalcontheta}, we have
$\Lambda_{ij}s_j(X)=s_i(X)$, which implies that we have a global splitting $s:TM\longrightarrow A.$  \qed\vspace{3mm}

A $1$-morphism from  $\Big(\sqcup (TU_i\oplus (U_i\times\g),
  \brtls{-,-},  \widetilde{\theta}_i),\widetilde{\phi}_{ij},\widetilde{\tau}_{ij}\Big) $ to  $\Big(\sqcup (TU_i\oplus (U_i\times\g),
  \brtls{-,-},  \theta_i),\phi_{ij},\tau_{ij}\Big) $ in   $\holim \tlp (U(M)_\bullet)$  consists of
  a 1-morphism $\left(\begin{array}{cc}
1&0\\
\phi_{i}&\tau_{i}
\end{array} \right)$ from $\sqcup (TU_i\oplus (U_i\times\g),
   \widetilde{\theta}_i) $ to $\sqcup (TU_i\oplus (U_i\times\g),
   \theta_i) $  in $\tlp(\cup U_i)_1$, which satisfies
  \begin{equation}\label{eq:gluetlm}
   \Lambda_{ij}\left(\begin{array}{cc}
1&0\\
\phi_{j}&\tau_{j}
\end{array} \right)=\left(\begin{array}{cc}
1&0\\
\phi_{i}&\tau_{i}
\end{array} \right)\widetilde{\Lambda}_{ij}.
  \end{equation}
 We have

  \begin{pro}\label{pro:tlm}
  A  $1$-morphism  in    $\holim \tlp (U(M)_\bullet)$ gives rise to a Lie algebroid isomorphism
  preserving connections.
  \end{pro}
\pf  Denote by $(\widetilde{A},[-,-\widetilde{]},\widetilde{\rho})$ (respectively  $(A,[-,-],\rho)$) the transitive Lie algebroid with the connection $\widetilde{s}:TM\longrightarrow \widetilde{A}$ (respectively $s: TM \longrightarrow A$) obtained from the object $\Big(\sqcup (TU_i\oplus (U_i\times\g),
  \widetilde{\theta}_i),\widetilde{\phi}_{ij},\widetilde{\tau}_{ij}\Big) $ (respectively $\Big(\sqcup (TU_i\oplus (U_i\times\g),
     \theta_i),\phi_{ij},\tau_{ij}\Big) $)  in  $\holim \tlp (U(M)_\bullet)$. Thanks to   \eqref{eq:gluetlm},  a 1-morphism $\{\left(\begin{array}{cc}
1&0\\
\phi_{i}&\tau_{i}
\end{array} \right)\}$ from  $\Big(\sqcup (TU_i\oplus (U_i\times\g),
  \widetilde{\theta}_i),\widetilde{\phi}_{ij},\widetilde{\tau}_{ij}\Big) $ to $\Big(\sqcup (TU_i\oplus (U_i\times\g),
  \brtls{-,-},  \theta_i),\phi_{ij},\tau_{ij}\Big) $ glues to a bundle map which gives rise to a Lie algebroid isomorphism between $(\widetilde{A},[-,-\widetilde{]},\widetilde{\rho})$ and $(A,[-,-],\rho)$.

Furthermore, we have
\begin{eqnarray*}
   \left(\begin{array}{cc}
1&0\\
\phi_{i}&\tau_{i}
\end{array} \right)\widetilde{s}_i(X)=X+\phi_i(X)+\tau_i\widetilde{\theta}_i(X)=X+\theta_i(X)=s_i(X),
\end{eqnarray*}
which implies that the connections are also preserved.\qed

\begin{rmk}
 The above way to glue a transitive Lie algebroid is essentially the same as the one given by Mackenzie \cite{MK2}. %for more details. Since this is a part of gluing transitive Courant algebroids, we give the result here using the language of $(2,1)$-presheaf to be self-contained.
\end{rmk}

By Proposition \ref{pro:tlo} and \ref{pro:tlm},  it is not hard to see that after the plus construction we arrive at a $(2,1)$-sheaf ${\tlp}^+$ which maps to the category of transitive Lie algebroids with connections essentially surjectively and fully faithfully.

\subsection{$(2,1)$-sheaf ${\ecalgdpp}$
 of exact Courant algebroids with connections}\label{app:ec}

The {\bf standard Courant algebroid} is  $(TM\oplus T^*M,\Courant{-,-}^E_S,\pair{-,-}^E,\pr_{TM})$, where  $\Courant{-,-}^E_S$ is the standard Dorfman bracket given by
\begin{equation}\label{eq:standardb}
  \Courant{X+\xi,Y+\eta}^E_S=[X,Y]+L_X\eta-\ii_Yd\xi,
\end{equation}
and $\pair{-,-}^E$ is the canonical symmetric bilinear form given by
  \begin{equation}\label{eq:standardpair}
  \pair{X+\xi,Y+\eta}^E=\half\big( \xi(Y)+\eta(X)\big),
\end{equation}

A Courant algebroid $C$ is called {\bf exact} if we have the following short exact sequence
$$
0\longrightarrow T^*M\stackrel{\rho^*}{\longrightarrow} C\stackrel{\rho}{\longrightarrow}TM\longrightarrow 0.
$$
A {\bf connection} of an exact Courant algebroid $C$ is an isotropic splitting\footnote{A splitting $s:TM\longrightarrow C$ is called isotropic if the image of $s$ is an isotropic subbundle, i.e. $\langle s(X),s(Y)\rangle=0$, for all $X,Y\in \Gamma(TM)$.} $s:TM\longrightarrow C$. As before, connections always exist. By choosing a connection $s:TM\longrightarrow C$, the vector bundle $C$ is isomorphic to $TM\oplus T^*M$. Then transferring the Courant algebroid structure on $C$ to that on $TM\oplus T^*M$, we obtain the Courant algebroid $(TM\oplus T^*M, \Courant{-,-}_h,\pair{-,-}^E, \pr_{TM})$, where the nondegenerate symmetric pairing  $\pair{-,-}^E$ is given by \eqref{eq:standardpair} and the bracket $\Courant{-,-}^E_h$ is given by
\begin{equation}\label{eq:hbracket}
  \Courant{X+\xi,Y+\eta}^E_h= \Courant{X+\xi,Y+\eta}^E_S+\ii_Y\ii_Xh.
  \end{equation}
Here  $h\in\Omega^3_{\rm cl}(M)$, defined by $h(X,Y)=\Courant{s(X),s(Y)}-s[X,Y]$,  is the {\bf curvature} of the connection $s$. In \cite{Severa:3-form}, the authors show that exact Courant algebroids over $M$ are classified by $H^3(M, \R)$.

Now we construct the $(2,1)$-presheaf of exact Courant algebroids with
connections over
the category of (differential) manifolds $\Mfd$. For simplicity, for an object $U\in \Mfd$, we write $\tangu:=TU\oplus T^*U.$

%\begin{pro}
There is a $(2,1)$-presheaf of exact Courant algebroids with
connections $\ecalgdp: \Mfd^{\op} \to \Gpd$, where $\Mfd^{\op}$ is the
opposite  category of $\Mfd$, and $\Gpd$ is the
$2$-category of (set theoretical) groupoids and groupoid morphisms.

For an object $U\in \Mfd$, the groupoid $\ecalgdp(U)$ is made up
by the following data:
\begin{itemize}
\item{$\ecalgdp(U)_0$}: an object is a quintuple $(\tangu, \Courant{-,-}^E_S,  \pair{-,-}^E,\pr_{TU}, B)$, where $B\in \Omega^2(U)$ is a $2$-form,  $\Courant{-,-}^E_S$ and
  $\pair{-,-}^E$ are  given by \eqref{eq:standardb} and \eqref{eq:standardpair} respectively. We will simply denote an object by $(\tangu,B)$.

\item{$\ecalgdp(U)_1$}:  a $1$-morphism from $(\tangu,  \widetilde{B})$ to $ (\tangu,B) $ is a bundle automorphism of $\tangu$ given by the
  matrix
$
\left(\begin{array}{cc}
1&0\\
\huaB&1
\end{array} \right),
$ where $\huaB\in \Omega^2(U)$ is a closed $2$-form such that $\widetilde{B}-B=\huaB$.
This matrix preserves the standard Courant bracket $\Courant{-,-}^E_S$ and the pairing
$\pair{-,-}^E$. The
composition of $1$-morphisms is simply the matrix multiplication.

\end{itemize}
Then for a morphism $\varphi: U\to V$ in $\Mfd$, as in the case of Lie algebroids, the associated
functor $\ecalgdp(\varphi): \ecalgdp(V) \to \ecalgdp(U)$ is induced by
pulling back forms.\calc{, more precisely, on the level of objects,
\[\ecalgdp(\varphi)(TV\oplus T^*V, \Courant{-,-}^E_S,  \pair{-,-}^E,\pr_{TV}, B_V)
= (\tangu, \Courant{-,-}^E_S,  \pair{-,-}^E,\pr_{TU}, \varphi^*B_V),\]
and on the level of morphisms,
\[
\ecalgdp(\varphi) \left(\begin{array}{cc}
1&0\\
\huaB_V&1
\end{array} \right) = \left(\begin{array}{cc}
1&0\\
\varphi^*\huaB_V&1
\end{array} \right) .
\]
One may verify that $\ecalgdp(\varphi)$ is indeed a functor between
desired categories.}
%  \end{pro}

\vspace{3mm}
Take an open cover $\{U_i\}$ of $M\in
\Mfd$. An object in $\holim \ecalgdp(U(M)_\bullet)$ consists of
\begin{itemize}
\item  an object $\sqcup (\tangui,B_i)$ in $\ecalgdp(\sqcup U_i)_0$,

\item  $\Lambda_{ij}=\left(\begin{array}{cc}
1&0\\
\huaB_{ij}&1
\end{array} \right) \in \ecalgdp(\sqcup U_{ij})_1$ which is a 1-morphism from
$
(\tanguij,
 B_j|_{U_{ij}})
  $
  to
   $
   (\tanguij, B_i|_{U_{ij}}),
$
therefore satisfying
$$
B_j|_{U_{ij}}-B_i|_{U_{ij}}=\huaB_{ij}.
$$
\item compatibility conditions $\Lambda_{ij}\Lambda_{jk}=\Lambda_{ik}$ on $U_{ijk}$ which automatically holds.
\end{itemize}

The plus construction gives us a $(2,1)$-sheaf $\ecalgdpp$. For a manifold $M$, an
object of $\ecalgdpp(M)$ consists of a cover $\{U_i\}$ and a
$U(M)$-equivariant object of $\ecalgdp$ described above. Naturally
we ask what the above data glues to. It turns out that the gluing result is an exact Courant
algebroid with a connection. The gluing procedure is the same as the one given in \cite{Rogers:prequan}. To be self-contained, we give the result using the language of this paper.

\begin{pro}\label{thm:eco}\label{pro:eco}
An object  $\Big(\sqcup (\tangui,
 B_i),\huaB_{ij}\Big)$ in $\holim \ecalgdp(U(M)_\bullet)$  gives rise to an exact Courant algebroid $(C,\Courant{-,-},\pair{-,-},\rho)$ with a connection $s:TM\longrightarrow C$.
\end{pro}
\pf
Given an object $\Big(\sqcup (\tangui,
  B_i),\huaB_{ij}\Big)$ in $\holim\ecalgdp(U(M)_\bullet)$, as before,  the cocycle condition $\Lambda_{ij}\Lambda_{jk}=\Lambda_{ik}$ implies that $\taof{U_i}$'s glue to a vector bundle $C$ via transition matrices $\Lambda_{ij}$'s.
Since $\Lambda_{ij}$ preserves the standard bracket $\Courant{-,-}^E_S$, we have a well-defined bracket $\Courant{-,-}$ on $\Gamma(C)$.
Furthermore, since $\Lambda_{ij}$ also preserves the standard pairing $\pair{-,-}^E$ on $\tangui$, we obtain a global nondegenerate symmetric bilinear form $\pair{-,-}$ on $C$.
Obviously, $C$ fits the following exact sequence of vector bundles,
$$
0\longrightarrow T^*M  \stackrel{\rho^*}{\longrightarrow}C \stackrel{\rho} {\longrightarrow } TM\longrightarrow 0,
$$
where $\rho$ is induced by the projection $\tangui\longrightarrow TU_i$.

Also, by the facts that $\Lambda_{ij}$ preserves  the standard bracket $\Courant{-,-}^E_S$ and the standard pairing $\pair{-,-}^E$, Axioms (i)-(iii) in Definition \ref{defi:ca} are satisfied. Therefore, $(C,\Courant{-,-},\pair{-,-},\rho)$ is an exact Courant algebroid.

The 2-forms $\{B_i\}$ induce an isotropic splitting $s:TM\longrightarrow C$ via
$$
s(X)=X-\ii_XB_i, \quad X\in U_i.
$$
Note that the definition of $s$ does not depend on choices of $U_i$. In fact, if $X\in U_i\cap U_j$, it is straightforward to see that $X-\ii_XB_i \thicksim X-\ii_XB_j$. \qed\vspace{3mm}

In $\holim \ecalgdp(U(M)_\bullet)$, a  1-morphism from  an object $\Big(\sqcup (\tangui,
   \widetilde{B}_i),\widetilde{\huaB}_{ij}\Big) $  to another object $\Big(\sqcup (\tangui
   B_i),\huaB_{ij}\Big) $ consists of
  a 1-morphism $\left(\begin{array}{cc}
1&0\\
\huaB_{i}&1
\end{array} \right)$ from $ \sqcup (\tangui,
   \widetilde{B}_i)$ to $ \sqcup (\tangui,
    B_i)$  in $\ecalgdp(\cup U_i)_1$, which satisfies
  \begin{equation}\label{eq:gluecam}
    \Lambda_{ij}\left(\begin{array}{cc}
1&0\\
\huaB_{j}&1
\end{array} \right)=\left(\begin{array}{cc}
1&0\\
\huaB_{i}&1
\end{array} \right)\widetilde{\Lambda}_{ij}.
  \end{equation}
 Then we have

  \begin{pro} \label{pro:ecm}
 A  $1$-morphism in $\holim \ecalgdp(U(M)_\bullet)$ gives rise to an exact Courant algebroid isomorphism preserving connections.
  \end{pro}
\pf
The proof is similar to that of Proposition \ref{pro:tlm}. Eq. \eqref{eq:gluecam} is the important information which implies the gluing result. The fact that the bundle map $\frkB$ also preserves the connection, namely $\frkB(\widetilde{s}(X))=s(X)$, follows from the following calculation,
$$
\frkB(\widetilde{s}(X))=X-\ii_X\widetilde{B}_i+\huaB_i=X-\ii_XB_i=s(X),\quad \forall X\in \Gamma(TM).
$$
The proof is finished. \qed\vspace{3mm}

Similar to the case of transitive Lie algebroids, after the plus construction, we arrive at the $(2,1)$-sheaf ${\ecalgdp}^+$ of exact Courant algebroids with connections.

%up to here 21.11 afternoon.

\subsection{Gluing transitive Courant algebroids via local data}\label{GluingCA}

In this subsection, we give the explicit formula for the transitive Courant algebroid glued by pieces of standard transitive Courant algebroids in a transitive Courant data. This also shows how we may obtain the bracket of a general transitive Courant algebroid from the standard one. As in Proposition \ref{pro:tco}, given a transitive Courant data $C_c$,  there is a corresponding transitive Courant algebroid $(C, \Courant{-, -}, \pair{-,-}, \rho)$. Using the two splittings $s$ and $\sigma_s$ given in \eqref{eq:tcs1} and \eqref{eq:tcs2}, we obtain an isomorphism  $\huaS:TM\oplus \huaG\oplus T^*M\longrightarrow C$ given by
\begin{equation}\label{iso:tc}
\huaS(X+a+\xi)=s(X)+\sigma_s(a)+\xi.
\end{equation}
\emptycomment{Any section $e\in\Gamma(C)$ can be written as $e=s(X)+\sigma_s(a)+\xi$ for some $X\in\frkX(M),a\in\Gamma(\huaG)$ and $\xi\in\Omega^1(M)$.
Now we consider the induced pairing and bracket.} Recall that locally, $\huaS_i=\huaS|_{U_i}:TU_i\times(U_i\times\g)\times T^*U_i\longrightarrow C|_{U_i}$ and its inverse are given by
\begin{eqnarray*}
  \huaS_i(X+a+\xi)&=&X+\theta_i(X)-(\theta_i,\theta_i(X))^\g-\ii_XB_i+a-2(\theta_i,a)^\g+\xi,\\
   \huaS_i^{-1}(X+a+\xi)&=&X-\theta_i(X)-(\theta_i,\theta_i(X))^\g+\ii_XB_i+a+2(\theta_i,a)^\g+\xi.
\end{eqnarray*}
Having $TM\oplus \huaG\oplus T^*M$ equipped with the pairing given by \eqref{eq:Tpair}, a straightforward computation shows that $\huaS$ preserves the pairing.
\emptycomment{
\begin{lem}
Local isomorphisms $\huaS_i$'s preserve the pairing, that is,  $$\pair{\huaS_i(X+a+\xi),\huaS_i(Y+b+\eta)}^T=\pair{X+a+\xi,Y+b+\eta}^T.$$
Therefore,  $(C,\pair{-,-})$ is isomorphic to $(TM\oplus \huaG\oplus T^*M,\pair{-,-}^T)$.
\end{lem}}

\begin{lem}\label{lem:tcc}
Define $\nabla^i:\Gamma(TU_i)\otimes \Gamma(U_i\times \g)\longrightarrow \Gamma(U_i\times \g)$ by
\begin{equation}\label{eq:tcconnectionlocal}
 \nabla^i_Xa =X(a)+[\theta_i(X),a]_\g, \quad \forall X\in\Gamma(TU_i), a\in\Gamma(U_i\times \g).
\end{equation}
Then, we have
  $$
  \tau_{ij}\nabla^j_Xa=\nabla^i_X\tau_{ij}a,\quad \forall X\in\Gamma(TU_{ij}), a\in\Gamma(U_{ij}\times \g).
  $$
  Thus, by gluing  $\nabla^i$, we obtain  a globally well-defined connection $\nabla:\Gamma(TM)\otimes \Gamma(\huaG)\longrightarrow \Gamma(\huaG)$.
\end{lem}
\pf By \eqref{eq:auto2} and \eqref{eq:morcon11}, we have
\begin{eqnarray*}
  \tau_{ij}\nabla^j_Xa-\nabla^i_X\tau_{ij}a&=&\tau_{ij}(X(a)+[\theta_j(X),a]_\g)-X(\tau_{ij}a)-[\theta_i(X),\tau_{ij}a]_\g\\
  &=&\tau_{ij}(X(a))+[\tau_{ij}\theta_j(X),\tau_{ij}a]_\g)-X(\tau_{ij}a)-[\tau_{ij}\theta_j(X)+\phi_{ij}(X),\tau_{ij}a]_\g\\
  &=&0.
\end{eqnarray*}
The proof is finished.  \qed\vspace{3mm}

Now we see that given a Courant data $C_c$, we have a connection $\nabla$, a curvature $R$ of the underlining Lie data $A_c$ given in Lemma \ref{lem:gglueR} and a curvature $H$ of $C_c$. Thus $TM\oplus \huaG\oplus T^*M$ may be equipped with a transitive Courant algebroid structure with the Courant bracket $\Courant{-,-}^T_{\nabla, R, H}$ given as in \eqref{eq:Tbracket}.

\begin{pro}\label{GlobalCA}
The morphism $\huaS$ in \eqref{iso:tc} is an isomorphism of Courant algebroids.
 \end{pro}
 \pf We pull back the bracket on $\Gamma(C)$ to $\Gamma (TM\oplus \huaG\oplus T^*M)$ via $\huaS$ and denote it by $\Courant{-,-}_{ind}$. The only nontrivial thing to check is that $\Courant{-,-}_{ind}=\Courant{-,-}^T_{\nabla, R, H}$.  For all $a,b\in\Gamma(U_i\times \g)$, we have
\begin{eqnarray*}
  \Courant{a,b}_{ind}&=&\huaS_i^{-1}\Courant{\huaS_i(a),\huaS_i(b)}^T_S=\huaS_i^{-1}\Courant{a-2(\theta_i,a)^\g,b-2(\theta_i,b)^\g}^T_S=S_i^{-1}([a,b]_\g+\huaP(a,b))\\
  &=&[a,b]_\g+\huaP(a,b)+2(\theta_i,[a,b]_\g)^\g.
\end{eqnarray*}
By \eqref{eq:tcconnectionlocal} and Lemma \ref{lem:tcc}, we have
\begin{eqnarray*}
  \Big(\huaP(a,b)+2(\theta_i,[a,b]_\g)^\g\Big)(Y)&=&2(b,Y(a))^\g+2(\theta_i(Y),[a,b]_\g)^\g\\
  &=&2(b,Y(a)+[\theta_i(Y),a]_\g)^\g=2(b,\nabla^i_Ya)^\g\\
  &=&P(a,b)(Y),
\end{eqnarray*}
which implies that
\begin{equation}
  \label{eq:gbr1}\Courant{a,b}_{ind}=[a,b]_\g+P(a,b).
\end{equation}
For all $X\in\Gamma(TU_i), b\in\Gamma(U_i\times \g)$, we have
\begin{eqnarray*}
  \Courant{X,b}_{ind}&=&\huaS_i^{-1}\Courant{\huaS_i(X),\huaS_i(b)}^T_S=\huaS_i^{-1}\Courant{X+\theta_i(X)-(\theta_i,\theta_i(X))^\g-\ii_XB_i,b-2(\theta_i,b)^\g}^T_S\\
  &=&\huaS_i^{-1}(X(b)-2L_X(\theta_i,b)^\g+[\theta_i(X),b]_\g+\huaP(\theta_i(X),b))\\
  &=&X(b)+[\theta_i(X),b]_\g-2L_X(\theta_i,b)^\g+\huaP(\theta_i(X),b)+2(\theta_i,X(b)+[\theta_i(X),b]_\g)^\g.
\end{eqnarray*}
By \eqref{eq:tcconnectionlocal}
and
\begin{eqnarray*}
 \Big(-2L_X(\theta_i,b)^\g+\huaP(\theta_i(X),b)+2(\theta_i,X(b)+[\theta_i(X),b]_\g)^\g\Big)(Y)=-2(R(X,Y),b)^\g=-2Q(X,b)(Y),
\end{eqnarray*}
we get
\begin{equation} \label{eq:gbr2}
\Courant{X,b}_{ind}=\nabla_Xb-2Q(X,b).
  \end{equation}
Similarly, we have
\begin{equation}\label{eq:gbr3}
\Courant{a,Y}_{ind}=2Q(Y,a)-\nabla_Ya.
  \end{equation}
 For all $X,Y\in\Gamma(TU_i)$, we have
  \begin{eqnarray*}
  \Courant{X,Y}_{ind}&=&\huaS_i^{-1}\Courant{\huaS_i(X),\huaS_i(Y)}^T_S\\
  &=&\huaS_i^{-1}\Courant{X+\theta_i(X)-(\theta_i,\theta_i(X))^\g-\ii_XB_i,Y+\theta_i(Y)-(\theta_i,\theta_i(Y))^\g-\ii_YB_i}^T_S\\
  &=&\huaS_i^{-1}\Big([X,Y]+X(\theta_i(Y))-Y(\theta_i(X))+[\theta_i(X),\theta_i(Y)]_\g\\
  &&-L_X(\theta_i,\theta_i(Y))^\g+L_Y(\theta_i,\theta_i(X))^\g+\huaP(\theta_i(X),\theta_i(Y))-L_X\ii_YB_i+\ii_Yd\ii_XB_i\Big)\\
  &=&[X,Y]-\theta_i([X,Y])+X(\theta_i(Y))-Y(\theta_i(X))+[\theta_i(X),\theta_i(Y)]_\g\\
  &&-L_X\ii_YB_i+\ii_Yd\ii_XB_i+\ii_{[X,Y]}B_i+\Xi,
  \end{eqnarray*}
  where
  \begin{eqnarray*}
\Xi&=&  -L_X(\theta_i,\theta_i(Y))^\g+L_Y(\theta_i,\theta_i(X))^\g+\huaP(\theta_i(X),\theta_i(Y))-(\theta_i,\theta_i([X,Y]))^\g\\
&&+2(\theta_i,X(\theta_i(Y))-y(\theta_i(X))+[\theta_i(X),\theta_i(Y)]_\g)^\g.
\end{eqnarray*}
 Obviously, we have
\begin{eqnarray*}
  -\theta_i([X,Y])+X(\theta_i(Y))-Y(\theta_i(X))+[\theta_i(X),\theta_i(Y)]_\g&=&d\theta_i(X,Y)+[\theta_i(X),\theta_i(Y)]_\g=R_i(X,Y),\\
  -L_X\ii_YB_i+\ii_Yd\ii_XB_i+\ii_{[X,Y]}B_i&=&-dB_i(X,Y,\cdot).
\end{eqnarray*}
Furthermore, we have
  \begin{eqnarray*}
    \Xi(Z)&=&2(\theta_i(Y),Z\theta_i(X))^\g-X(\theta_i(Z),\theta_i(Y))^\g+(\theta_i[X,Z],\theta_i(Y))^\g+d(\theta_i,\theta_i(X))^\g(y,Z)\\
    &&-(\theta_i(Z),\theta_i[X,Y])^\g+2(\theta_i(Z),X(\theta_i(Y))-Y(\theta_i(X))+[\theta_i(X),\theta_i(Y)]_\g)^\g\\
    &=&(\theta_i(X),d\theta_i[Y,Z])^\g+c.p.+2(\theta_i(Z),[\theta_i(X),\theta_i(Y)]_\g)^\g\\
    &=&\Big((\theta_i,d\theta_i)^\g+\frac{1}{3}(\theta_i,[\theta_i,\theta_i]_\g)^\g\Big)(X,Y,Z)\\
    &=&\cs_3(\theta_i)(X,Y,Z).
  \end{eqnarray*}
  Therefore, by Lemma \ref{lem:gglueR} and the fact that $-dB_i+\cs_3(\theta_i)$ can be glued to a global 3-form $H$, we have
  \begin{eqnarray}\label{eq:gbr4}
  \Courant{X,Y}_{ind}=[X,Y]+R(X,Y)+\Big(-dB_i+\cs_3(\theta_i)\Big)(X,Y,\cdot)=[X,Y]+R(X,Y)+H(X,Y,\cdot).
  \end{eqnarray}
  Furthermore, it is straightforward to obtain that
  \begin{eqnarray}
    \label{eq:gbr5}\Courant{X,\eta}_{ind}=L_X\eta,\quad\Courant{\xi,Y}_{ind}=-\ii_Yd\xi,\quad\Courant{a,\eta}_{ind}=0,\quad\Courant{\eta,b}_{ind}=0,\quad\Courant{\xi,\eta}_{ind}=0.
  \end{eqnarray}
  By \eqref{eq:gbr1}-\eqref{eq:gbr5}, we deduce that the induced bracket $\Courant{-,-}_{ind}$ is exactly given by  \eqref{eq:Tbracket}, i.e. $$\Courant{-,-}_{ind}=\Courant{-,-}^T_{\nabla,R,H}.$$
  The proof is finished. \qed

\subsection{Inner automorphisms of transitive Courant algebroids}\label{app:innerauto}
In this subsection, we prove that the automorphisms that appeared in Proposition \ref{Naturaltransformation} are inner automorphisms of the standard transitive Courant algebroid $(\tagu, \Courant{-,-}^T_S, \pair{-,-}^T, \pr_{TU})$. In his letter to Weinstein \cite{sevlet}, ${\rm\check{S}}$evera claimed that the inner automorphism group $\Inn(U)$ of the standard transitive Courant algebroid over $U$ is an extension of the group of $G$-valued function $C^\infty(U, G)$ by closed 2-forms $\Omega^2_{\cl}(U)$,
\begin{equation}
 \Omega^2_{\cl}(U) \to \Inn(U) \to  C^\infty(U, G).
\end{equation}
More precisely, an inner automorphism is a pair  $(g,\omega)$, where $g$ is a $G$-valued function and $\omega\in \Omega^2(U)$, such that
\begin{equation}\label{eq:gw}
d\omega+g^*\frkC=0,
\end{equation}
where $\frkC=\frac{1}{6}(\TM,[\TM,\TM])^\g$, or equivalently $\frkC(\hat{a},\hat{b},\hat{c})=(a,[b,c])^\g$. The group structure is given by
\begin{eqnarray*}
(g_1,\omega_1)(g_2,\omega_2)&=&(g_1g_2,\omega_1+\omega_2+(g_1^*\TM,\ad_{g_1}g_2^*\TM)^\g),\\
(g,\omega)^{-1}&=&(g^{-1},-\omega).
\end{eqnarray*}

Now we give the corresponding matrix form of an inner automorphism.
The matrix corresponding to $(g,\omega)$ is given by
\begin{equation}\label{eq:severainner}
\Psi=\left(\begin{array}{ccc}
1&0&0\\
-g^*\TM&\ad_{g}&0\\
\omega-(g^*\TM)^\star \circ g^*\TM &2(g^*\TM)^\star\circ \ad_\g&1
\end{array}\right).
\end{equation}

\begin{pro}
$\Psi$ given above is an automorphism of the standard transitive Courant algebroid.
\end{pro}
\pf
It is straightforward to see that \eqref{eq:auto1}-\eqref{eq:auto3} hold. For all $X,Y,Z\in\Gamma(TU)$, we have
$$\langle L_X\omega(Y)-\ii_Y d\omega(X)-\omega([X,Y]),Z\rangle=d\omega(X,Y,Z).$$
Denote by $\beta^{sym}=-(g^*\TM)^\star \circ g^*\TM $. By straightforward computations, we have
 $$
 \langle L_X\beta^{sym}(Y)-\ii_Y d\beta^{sym}(X)-\beta^{sym}([X,Y])+\huaP(g^*\TM,g^*\TM),Z\rangle =g^*\frkC(X,Y,Z).
 $$
By \eqref{eq:gw}, we deduce that \eqref{eq:auto4} holds. Thus $\Psi$ given above is an automorphism. \qed\vspace{3mm}

See \cite[Corollary 4.2]{GRT} for a similar result on inner automorphisms in another setting.

\begin{cor}
  $\Lambda_{01}$ given in Proposition \ref{Naturaltransformation} is an inner automorphism of the standard transitive Courant algebroid.
\end{cor}

\section*{Acknowledgement} The last author thanks much
Giovanni Felder for patiently going through complicated calculations together and pointing out correct directions for the most technical part of this article during her visit to ETH. The last two authors thank Anton Alekseev and Pavol ${\rm\check{S}}$evera for many inspiring conversations.
We warmly thank Olivier Brahic, Zhuo Chen, Domenico Fiorenza, Fei Han,
Kirill Mackenzie, Eckhard Meinrenken, Ralf Meyer, Thomas Nikolaus, Dmitry Roytenberg, Thomas Schick, Chris Schommer-Pries, Urs Schreiber, Konrad Waldorf and Marco Zambon for very helpful discussion. Y. Sheng is supported by NSFC (11471139) and NSF of Jilin Province (20170101050JC);
X. Xu was partially supported by the SNSF grants P2GEP2-165118 and NCCR SwissMAP; C. Zhu is supported by the German Research Foundation
(Deutsche Forschungsgemeinschaft (DFG)) through the Institutional
Strategy of the University of G\"ottingen.

{}

 \vspace{3mm}
Department of Mathematics, Jilin University,
 Changchun 130012,  China

Email: shengyh@jlu.edu.cn
 \vspace{3mm}

Department of Mathematics, Massachusetts Institute of Technology,
Cambridge 02139, United States

Email: xxu@mit.edu
 \vspace{3mm}

Mathematics Institute, Georg-August-University
G\"ottingen,  G\"ottingen 37073, Germany

Email: czhu@math.uni-goettingen.de
\end{document}